\documentclass{amsart}[11pt]
\usepackage[margin=4cm]{geometry}
\usepackage{amsmath, amsfonts, amsthm, amssymb, mathtools}
\usepackage{mathrsfs}
\usepackage{subfig}
\usepackage{stmaryrd}
\usepackage{verbatim}
\usepackage{hyperref}
\usepackage{xcolor}
\usepackage{graphicx}
\usepackage{algorithmic, algorithm}
\usepackage{enumerate}
\usepackage{listings}
\numberwithin{equation}{section}
\numberwithin{figure}{section}

\usepackage[textsize=tiny]{todonotes}

\makeatletter
\def\@tocline#1#2#3#4#5#6#7{\relax
  \ifnum #1>\c@tocdepth 
  \else
    \par \addpenalty\@secpenalty\addvspace{#2}%
    \begingroup \hyphenpenalty\@M
    \@ifempty{#4}{%
      \@tempdima\csname r@tocindent\number#1\endcsname\relax
    }{%
      \@tempdima#4\relax
    }%
    \parindent\z@ \leftskip#3\relax \advance\leftskip\@tempdima\relax
    \rightskip\@pnumwidth plus4em \parfillskip-\@pnumwidth
    #5\leavevmode\hskip-\@tempdima
      \ifcase #1
       \or\or \hskip 1em \or \hskip 2em \else \hskip 3em \fi%
      #6\nobreak\relax
    \dotfill\hbox to\@pnumwidth{\@tocpagenum{#7}}\par
    \nobreak
    \endgroup
  \fi}
\makeatother

\hypersetup{colorlinks,linkcolor=blue!55!black,citecolor=blue!55!black,urlcolor=blue!55!black}

\mathtoolsset{showonlyrefs=true}
\newtheorem{theorem}{Theorem}[section]
\newtheorem{corollary}[theorem]{Corollary}
\newtheorem{lemma}[theorem]{Lemma}
\newtheorem{proposition}[theorem]{Proposition}

\theoremstyle{definition}
\newtheorem{definition}[theorem]{Definition}
\newtheorem{remark}[theorem]{Remark}

\newtheorem{assumption}[theorem]{Assumption}

\newcommand{\ind}{1\hspace{-2.1mm}{1}} 
\newcommand{\half}{\frac{1}{2}}

\newcommand{\CC}{\mathbb{C}}
\newcommand{\RR}{\mathbb{R}}

\newcommand{\PP}{\mathbb{P}}

\newcommand{\EE}{\mathbb{E}}

\newcommand{\VV}{\mathbb{V}}

\newcommand{\Cc}{\mathcal{C}}

\newcommand{\Gg}{\mathcal{G}}
\newcommand{\Hh}{\mathcal{H}}
\newcommand{\Jj}{\mathcal{J}}
\newcommand{\Ff}{\mathcal{F}}
\newcommand{\Mm}{\mathcal{M}}
\newcommand{\Nn}{\mathcal{N}}
\newcommand{\Oo}{\mathcal{O}}
\newcommand{\Kk}{\mathcal{K}}
\newcommand{\Ll}{\mathcal{L}}

\newcommand{\Call}{\mathtt{C}}
\newcommand{\Put}{\mathtt{P}}
\newcommand{\D}{\mathrm{d}}
\newcommand{\E}{\mathrm{e}}
\newcommand{\I}{\mathrm{i}}

\newcommand{\BS}{\mathrm{BS}}
\newcommand{\eps}{\varepsilon}

\newcommand{\Cov}{\CC\mathrm{ov}}

\DeclareMathOperator*{\sgn}{sgn}

\newcommand{\Af}{\mathfrak{A}}
\newcommand{\Bf}{\mathfrak{B}}

\newcommand{\Df}{\mathfrak{D}}
\newcommand{\If}{\mathfrak{I}}

\newcommand{\Lf}{\mathfrak{L}}

\newcommand{\Ii}{\mathcal{I}}

\newcommand{\Dd}{\mathcal{D}}

\newcommand{\Lfbar}{\overline{\Lf}}

\newcommand{\speed}{\mathfrak{s}}

\title[LDP for integrated log-normals and variance options]{Large deviations for integrated exponentials of Gaussian Volterra processes - application to variance options}

\author{Antoine Jacquier}
\address{Department of Mathematics, Imperial College London}
\email{a.jacquier@imperial.ac.uk}
\author{Mugad Oumgari}
\address{Department of Mathematics, University College London}
\email{mugad.attahallah.20@ucl.ac.uk}
\author{Dan Pirjol}
\address{School of Business, Stevens Institute of Technology}
\email{dpirjol@stevens.edu}
\author{Lingjiong Zhu}
\address{Department of Mathematics, Florida State University}
\email{zhu@math.fsu.edu}

\date{\today}

\begin{document}

\begin{abstract}
We prove small-time large deviations for the normalised time integral of the exponential of a Gaussian Volterra process, and complement it with a central limit theorem and moderate deviations.
As an application, we obtain short-maturity asymptotics 
for out-of-the-money and at-the-money variance options in rough volatility models,
as well as implied volatility expansions which we validate numerically.
\end{abstract}

\maketitle
\setcounter{tocdepth}{4}
\tableofcontents

\section{Introduction}

Variance options, also known as options on realised variance, 
are derivative contracts whose payoffs depend on the annualised realised variance of the returns of an underlying asset.
Consider an underlying asset price whose dynamics follow a diffusion of the form
\begin{equation}\label{log:asset:dynamics}
\frac{\D S_{t}}{S_{t}}
=(r-q)\D t+\sigma_{t}\D B_{t},
\end{equation}
where $r$ is the risk-free rate, $q$ the dividend yield, 
$(\sigma_{t})_{t\geq 0}$ the instantaneous volatility and~$B$ is a standard Brownian motion.
On a uniformly spaced time grid $\{t_i = i\tau\}_{i=0,\ldots,n}$
with $\tau = \frac{T}{n}$ the time-step of the sampling period expressed in years (for daily sampling $\tau=\frac{1}{252}$),
the annualised realised variance is given by
\begin{equation}\label{RV:eqn}
\mathrm{RV}_{T,n} := \frac{1}{n\tau} \sum_{i=1}^n \log^2 \left(\frac{S_{t_{i}}}{S_{t_{i-1}}}\right).
\end{equation}
In practice, the discrete-sum in the expression of $\mathrm{RV}_{T,n}$ in~\eqref{RV:eqn} is approximated, as $n\uparrow\infty$, by the quadratic variation of the 
log-asset price
$[\log S, \log S]_T= \int_0^T \sigma^2_s \D s$,
such that, as $n\uparrow\infty$, $\mathrm{RV}_{T,n}$ in~\eqref{RV:eqn} can be approximated by
$\Sigma_{T}:=\frac{1}{T}\int_0^T \sigma^2_s \D s$.
The prices of the variance Call and Put options with strike $K>0$ 
and maturity $T>0$ are given by
$$
\Call(K,T)
:= \EE\left[\left(\Sigma_{T} - K\right)_{+}\right]
\qquad\text{and}\qquad
\Put(K,T)
:= \EE\left[\left(K-\Sigma_{T}\right)_{+}\right].
$$
We are interested in studying the short-maturity asymptotics
and implied volatility expansions for these prices
in rough volatility models of the form
\begin{align}\label{rough:vol:eqn}
\sigma_{t}:=\sigma(Y_{t}),
\qquad \text{with}\qquad
Y_{t}=y_{0}+W_{t}^{H},
\end{align}
where~$W^{H}$ is a fractional Brownian motion with Hurst index $H\in(0,1)$,
possibly correlated with the Brownian motion~$B$ in~\eqref{log:asset:dynamics}. 
When $\sigma(y)=\E^{y}$ in~\eqref{rough:vol:eqn}, 
then $\Sigma_{T}=\frac{\E^{2y_{0}}}{T}\int_{0}^{T}\E^{2W_{t}^{H}}\D t$. More generally, 
we consider $\Sigma_{T} = \frac{Z_{T}}{T}$ where
\begin{equation}
Z_{T}:=\int_{0}^{T}\exp\left\{\nu X_{s}^{(\ell)}+\varphi(s)\right\}\D s,
\end{equation}
where $\varphi:\RR_{+}\rightarrow\RR$ is a continuous function, $\nu>0$ is a fixed constant, and
$$
X_{t}^{(\ell)}=\int_{0}^{t}K_{\ell}(t-u)\D W_{u}, 
$$
for some function $K_{\ell}$
whose properties will be specified in Assumption~\ref{ass:main}, and where~$W$ is a standard Brownian possibly correlated with~$B$.
We will first study
the large deviation principle (LDP), central limit theorem (CLT)
and moderate deviation principle (MDP) for $Z_{T}/T$
in the small-time regime as $T\downarrow 0$,
and will then transfer these results to obtain short-maturity
asymptotics for out-of-the-money
and at-the-money variance options. 


Our work is related to several streams of literature:
Variance options, Short-maturity asymptotics under rough volatility, 
Small-time LDPs for realised variance in rough models,
LDPs for Gaussian and Volterra processes and Short-maturity Asian options.

\subsubsection*{Variance options}
The earliest studies on pricing variance options in stochastic volatility models such as jump-augmented Heston are due to Sepp~\cite{Sepp2008HestonJumps,Sepp2008}. 
Carr, Geman, Madan and Yor~\cite{CGMY2005} proposed a Laplace transform method for pricing options on realised variance in exponential L\'evy models. 
Drimus~\cite{Drimus2012} discussed the pricing and hedging of variance options in the $\frac{3}{2}$-volatility model using transform-based methods
and in~\cite{Drimus2012LogOU} studied variance options for a general class of Log-OU stochastic volatility models.
A slightly different angle was proposed by Bergomi, who developed forward variance models,
and derived variance options pricing in such models~\cite{Bergomi1, Bergomi2,BergomiBook}.

\subsubsection*{Short-maturity asymptotics under rough volatility}
In the context of rough volatility models,
Forde and Zhang~\cite{forde2017asymptotics} showed that 
$(T^{H-\half}\log(S_{T}))$ satisfies a large deviation principle as $T\downarrow0$ and derived asymptotics for European options, 
and the moderate deviations (near-the-money) regime is treated in~\cite{bayer2019short}. 
The short-time smile skew and curvature are given in Fukasawa~\cite{fukasawa2017short} and Al\`os and Le\'on~\cite{alos2017curvature}, 
and a second-order short-time expansion is obtained in~\cite{euch2019short} for general rough volatility models.
Pathwise large deviations for rough Bergomi are studied in~\cite{jacquier2018pathwise}, 
and pathwise large and moderate deviations for Gaussian rough volatility models in Gulisashvili~\cite{gulisashvili2020gaussian,gulisashvili2021time}.
The case of the rough Heston model is tackled in~\cite{Forde2021}.
Finally, Friz, Gassiat and Pigato~\cite{friz2021precise} introduced a new methodology to analyse a large class of rough volatility models, with further analysis in the large and moderate deviations regimes~\cite{friz2022short}.

\subsubsection*{Small-time LDPs for realised variance in rough models} 
The short-time asymptotics of the at-the-money volatility and skew of options on realised variance in rough volatility models were obtained by Malliavin calculus methods by Al\`os, Garc\'ia-Lorite and Muguruza~\cite{Alos2018}.
Lacombe, Muguruza, Stone ~\cite{lacombe2021asymptotics} gave large deviations predictions for the small-time behaviour of rough Bergomi and its multi-factor extensions. 
Theorem~\ref{thm:LDP} below is closest to their result
but we consider (a) a general regularly varying kernel rather than the Riemann-Liouville one, with the observation (Theorem~\ref{thm:pathLDP}), 
(b) the additional factor~$\varphi$, 
and (c) the Laplace/negative-moment statement,
as well as a detailed tail analysis in Section~\ref{sec:tails}.
The short maturity asymptotics of variance options in the Markovian case $H=\half$ is well studied:
in stochastic volatility models, options on realised variance are analogous to Asian options under a local volatility model, 
and the results carry over directly. 
Section~\ref{sec:summary:Asian} below provides a summary of 
large deviations results in this case. 
Short-maturity asymptotics of variance options in a Markovian local-stochastic volatility framework was obtained recently in~\cite{pirjol2025}.

\subsubsection*{LDPs for Gaussian and Volterra processes}
Large deviations for Gaussian and Volterra processes 
were developed in~\cite{forde2017asymptotics} 
for~\eqref{log:asset:dynamics}-\eqref{rough:vol:eqn},
with Gulisashvili~\cite{gulisashvili2018large} providing similar results  under weaker assumptions,
further~\cite{gulisashvili2020gaussian} establishing sample-path large and moderate deviation principles for Gaussian volatility models, 
and deriving asymptotics for exit probabilities, option prices, 
and implied volatility (and extended in~\cite{gulisashvili2021time} to more super rough Gaussian models).
The rough Bergomi case was tackled in~\cite{jacquier2018pathwise}
and randomised fractional volatility models in~\cite{horvath2019asymptotic}.
In~\cite{friz2021precise}, Friz, Gassiat and Pigato apply Laplace methods and regularity structures to obtain precise asymptotics for European options for large classes of classical and rough stochastic volatility models.
Finally, \cite{bayer2019short} studies the short-maturity near-the-money option pricing in rough fractional volatility models,
while~\cite{giorgio2023short} derives small-time asymptotics for non-self-similar stochastic volatility models
and~\cite{jacquier2022large} establishes large and moderate deviation principles for stochastic Volterra systems.

\subsubsection*{Asian options}\label{sec:summary:Asian}
For $H=\half$ and $\ell\equiv1$, $Z_{t}$ is the classical exponential functional of Brownian motion~\cite{dufresne1990distribution,matsumoto2005exponential} 
and the payoff functional of a continuously monitored Asian option. 
The short-maturity analysis of such options in local volatility models has been carried out in a series of studies in~\cite{pirjol2016short,arguin2018,pirjol2019short,pirjol2023subleading}.
In particular, \cite{pirjol2016short} obtains leading-order short-maturity asymptotics for OTM and ATM with bounded volatility function, \cite{arguin2018} derives the most-likely-path in Asian option pricing, and~\cite{pirjol2019short} carries out the detailed analysis for the CEV model.
Finer asymptotics such as sub-leading order terms, can be found in~\cite{pirjol2023subleading}.
This setup was recently extended to local-stochastic volatility models in~\cite{pirjol2025asian}.

Our contributions can be summarised as follows.

\begin{itemize}
\item
For a Gaussian Volterra process with regularly varying kernel (Assumption~\ref{ass:main}), we study the small-time behaviour of~$Z_{\eps}$.
Theorem~\ref{thm:LDP} provides an LDP for $Z_{\eps}/\eps$ with rate function analysed in Section~\ref{sec:tails}, 
fully explicit for $H=\half$ (Proposition~\ref{prop:half}).
Theorems~\ref{thm:CLT} and~\ref{thm:MDP} complement it with central limit theorem (CLT) and moderate deviation principle (MDP).
\item Section~\ref{sec:options} translates these results into short-maturity asymptotics for out-of-the-money (OTM) (Theorem~\ref{thm:OTM}) and at-the-money (ATM) (Theorem~\ref{thm:ATM}) options on realised variance, and Section~\ref{sec:expansion} expands the rate function in log-moneyness around the at-the-money strike. 
\item
Based on the theoretical results in Section~\ref{sec:options}, we obtain expansions for the implied volatility for variance options in Section~\ref{sec:impliedvol}. 
We give explicit results for the level and the slope of the short-maturity implied volatility (Proposition~\ref{prop:sigimp}), which we compare with~\cite{Alos2018}
and which we validate numerically.
\end{itemize}

The paper is organised as follows:
Section~\ref{sec:setup} introduces the model setup
and the main results (large deviations, 
central limit theorem and moderate deviations)
for integrated exponentials of Gaussian Volterra processes.
Section~\ref{sec:options} develops an application to variance options in rough volatility models:
short-maturity asymptotics for out-of-the-money
and at-the-money variance option prices. 
We then derive implied volatility expansions
and conduct numerical experiments in Section~\ref{sec:impliedvol}. 
Appendix~\ref{sec:comment} discusses the difference between two commonly used kernels, 
and additional technical proofs are gathered in Appendix~\ref{sec:proofs_all}.

\section{Setup and main results}\label{sec:setup}
\subsection{Setup}
On a filtered probability space $(\Omega,\Ff,(\Ff_{t})_{t\ge0},\PP)$ supporting a standard one-dimensional Brownian motion~$W$, introduce the process
\begin{equation}\label{eq:X}
X^{(\ell)}_{t} := \int_0^t K_{\ell}(t-u) \D W_u,
\qquad \text{for }t\ge0,
\end{equation}
which is a continuous 
$\Ff^W$-adapted
Gaussian Volterra process \cite{baudoin2003equivalence} as long as
$K_{\ell}:(0,\infty) \to (0,\infty)$ satisfies the following assumption.

\begin{assumption}\label{ass:main}
$K_{\ell}(t) = t^{H-\half}\ell(t)$ for $t>0$, with $H\in(0,\half]$ and $\ell:(0,\infty)\to(0,\infty)$ slowly varying at~$0$ and bounded on compact subsets of $(0,\infty)$.
Furthermore $K_{\ell}$ is non-increasing on $(0,t_{0})$ for some $t_{0}>0$.
\end{assumption}
We write~$K = K_{\ell}$ whenever $\ell(\cdot)\equiv 1$.
This assumption, stronger than the usual square integrability condition, covers $K(t) = t^{H-\half}$ 
and $K_{\ell}(t) = t^{H-\half}\E^{-\gamma t}$ ($\gamma>0$), 
common in mathematical finance~\cite{bayer2023rough, bennedsen2022decoupling}.
The local monotonicity condition is mild and  used to obtain increment estimates uniform in~$t$ (Lemma~\ref{lem:holder}),
and can be replaced by any hypothesis guaranteeing~\eqref{eq:incr}. 
The assumption ensures that $K\in L^2_{\mathrm{loc}}((0,\infty))$ so that~\eqref{eq:X} is well defined and, for any $t>0$,
$
\VV[X^{(\ell)}_{t}]
 = \int_{0}^{t} u^{2H-1} (\ell(u))^2 \D u$
 is finite
 and $\VV[X^{(\ell)}_{\eps}] \sim \frac{\eps^{2H}(\ell(\eps))^2}{2H}$
as $\eps\downarrow 0$
by Karamata's theorem~\cite[Proposition~1.5.8]{bingham1987regular}. 
Accordingly, let
\begin{equation}\label{eq:eps}
\speed_{\eps}:= \nu \eps^{H}\ell(\eps),
\end{equation}
for some fixed $\nu>0$ and, with $\varphi:\RR_+\to\RR$ continuous, set
\begin{equation}\label{eq:Z}
Z_{t} := 
\int_{0}^{t}\exp\left\{\nu X^{(\ell)}_{s} + \varphi(s)\right\}\D s,\qquad t>0,
\end{equation}
and
\[
\Ll(u,t) := \EE\left[Z_{t}^{-u}\right],
\qquad\text{for $u \in \Dd_{t}:=\left\{u\in\RR:\EE\left[Z_{t}^{-u}\right]<\infty\right\}$}.
\]
For each $t>0$, 
the random variable~$Z_{t}$ is the integrated variance of a lognormal (rough) stochastic volatility model
and the payoff functional of a continuously monitored arithmetic Asian option written on~$\E^{\nu X}$. 

Negative moments of~$Z_{t}$ arise for example in the mixing representation of option prices~\cite{gulisashvili2020small} in uncorrelated stochastic volatility models for at-the-money small-time expansions,
with $u\in \{\pm\half,\pm\frac{3}{2}\}$ 
and through the identity
\begin{equation}\label{eq:gamma}
z^{-u} = \frac{1}{\Gamma(u)}\int_0^\infty \theta^{u-1}\E^{-\theta z} \D\theta,\qquad z,u>0,
\end{equation}
in the analysis of the Laplace transform $u\mapsto\EE\left[\E^{-u Z_{t}}\right]$, studied in the Brownian case in~\cite{pirjol2023asymptotics} 
for pricing zero-coupon bonds in the Dothan model~\cite{dothan1978term}.

Small-time at-the-money asymptotics in Gaussian self-similar models depend on the moments of the integrated variance of orders $\half$ and~$\frac{3}{2}$~\cite{gulisashvili2020small}; these are instances of $\Ll(u,t)$ at fixed~$u$, governed by Proposition~\ref{prop:secondorder} rather than by Proposition~\ref{prop:laplace}. Extreme-strike asymptotics for the same class has been studied in \cite{gulisashvili2019extreme}.

\begin{remark}\ 
\begin{enumerate}[(i)]
\item The domain $\Dd_{t}$ is not a constraint 
as $\Dd_{t}=\RR$ for every $t>0$ 
under Assumption~\ref{ass:main} by Lemma~\ref{lem:moments}. 
Domain restrictions (moment explosions) appear in this class of models only for the asset price in the correlated case~\cite{gulisashvili2020gaussian,gassiat2019martingale}.
\item For fixed $u$ the limit is trivial: 
since $Z_{t}\sim t\E^{\varphi(0)}$ almost surely, 
then $\Ll(u,t) = t^{-u}\E^{-u\varphi(0)}(1+o(1))$,
hence
$\lim_{\eps\downarrow 0}\speed_{\eps}^{2}\log\Ll(u,\eps) = 0$. 
A non-degenerate Laplace-type limit requires some scaling as in~\cite{pirjol2023asymptotics}
for the Brownian case.
\end{enumerate}
\end{remark}

For $T>0$, 
let~$\Cc_0([0,T])$ denote the space of continuous functions $f:[0,T]\to\RR$ with $f(0)=0$ endowed 
with the supremum norm $\|\cdot\|_{[0,T]}$
and~$\Hh_{T}$ the Cameron-Martin space of absolutely continuous functions~$h:[0,T]\to\RR$ with $h(0)=0$, $\dot{h}\in L^2([0,T])$
and
$$
\|h\|_{\Hh_{T}}
:= \|\dot{h}\|_{L^2([0,T])}
= \left(\int_{0}^{T}
(\dot{h}(t))^2\D t\right)^{\half}.
$$
We define the bounded operator $\Kk_{\ell}:L^{2}([0,T])\to\Cc_{0}([0,T])$ by
\begin{equation}\label{eq:Kop}
\left(\Kk_{\ell}\dot{h}\right)(t)
:= \int_0^{t}K_{\ell}(t-u)\dot{h}(u) \D u,\qquad t\in[0,T];
\end{equation}
its range $\Hh_{T}^{X^{(\ell)}}:=\Kk_{\ell}(L^{2}([0,T]))$ is the Cameron-Martin space of 
$X^{(\ell)} = \Kk_{\ell}\dot{W}$.
This operator is injective:
for $\ell\equiv 1$, it follows from~\cite{samko1993fractional} while the general case follows from Titchmarsh's convolution theorem~\cite{titchmarsh1926zeros}.
Dropping the~$\ell$ notation for $\ell\equiv 1$, 
introduce 
the Gaussian process
$X := X^{(1)} = \Kk\dot{W}$,
with Cameron-Martin space
$\Hh_{T}^{X}:=\Kk(L^{2}([0,T]))\subset\Cc_{0}([0,T])$,
with norm $\|\Kk\dot{h}\|_{\Hh_{T}^{X}} = \|\dot{h}\|_{L^2([0,T])}$;
we further write
$\CC(t,s):=\EE[X_{t}X_s]$ and 
$(\CC f)(t) := \int_0^\infty \CC(t,s)f(s)\D s$
for the corresponding covariance function and operator.
For standard Brownian motion,
$K(\cdot)\equiv 1$, thus~$\Kk$ is the integration operator, $X=W$ and $\|\Kk\dot{h}\|_{\Hh_{T}^{X}}
= \|h\|_{\Hh_{T}}
= \|\dot{h}\|_{L^2([0,T])}$.
We finally introduce 
$\If_{t}:\Cc_{0}([0,t])\to[0,\infty]$ as
\begin{equation}\label{eq:I}
\If_{t}(f) :=
\left\{
\begin{array}{ll}
\displaystyle \half\left\|\Kk^{-1}f\right\|^{2}_{L^{2}([0,t])},
 & \text{if }f \in \Hh_{t}^{X},\\
+\infty, & \text{otherwise,}
\end{array}
\right.
\end{equation}
for $t\in [0,T]$,
with $\If:=\If_1$ 
(independent of~$\ell$), as well as the continuous functional
\begin{equation}\label{eq:Psi}
\Gg_{T}(f) := \log\left(\int_{0}^{T}\E^{f(s)}\D s\right),
\quad\text{for }f\in\Cc([0,T]),
\qquad\text{and}\qquad
\Gg := \Gg_{1}.
\end{equation}

In the next lemma, we show that the domain $\Dd_t$ is the real line for every $t>0$.

\begin{lemma}\label{lem:moments}
Under Assumption~\ref{ass:main}, $\Dd_{t}=\RR$ for every $t>0$.
\end{lemma}

\begin{proof}
The function 
$g(t):=\nu X^{(\ell)}_{t} + \varphi(t)$, 
is continuous, so that
Jensen's inequality yields
\[
Z_{t} = \int_{0}^{t}\E^{g(u)}\D u
\geq t\exp\left\{\frac{1}{t}\int_{0}^{t}g(u) \D u\right\};
\]
hence, for $v>0$,
$$
Z_{t}^{-v}\leq t^{-v}
\exp\left\{-\frac{v}{t}
\int_{0}^{t}\varphi(u)\D u\right\}
\E^{-v\nu G_{t}},
$$
with $G_{t} := \frac{1}{t}\int_{0}^{t}X^{(\ell)}_{u}\D u$ centred Gaussian with finite variance 
\[
\frac{1}{t^{2}}
\int_{0}^{t}\int_{0}^{t}\Cov\left(X^{(\ell)}_u,X^{(\ell)}_{u'}\right)\D u\D u',
\]
so that $\EE[\E^{-v\nu G_{t}}]$ is finite.
Now, for $v\geq 1$, Jensen's inequality implies that
$Z_{t}^{v}\leq t^{v-1}\int_{0}^{t}\E^{v g(u)}\D u$,
and 
$\EE\left[\int_{0}^{t}\E^{vg(u)}\D u\right]
= \int_{0}^{t}\E^{v\varphi(u)
 + \half v^2\nu^2\VV[X^{(\ell)}_{u}]}\D u$
is finite. 
 For $0<v<1$, $Z_{t}^v\le 1 + Z_{t}$.
 This completes the proof.
\end{proof}

Before we proceed, let us recall the definition
of large deviation principle (LDP). 
A family of measurable random elements $X_\eps$ in a Polish space $\mathcal{X}$
is said to satisfy a \textit{large deviation principle} (LDP) at speed $s_\eps\to\infty$ and rate
function $I:\mathcal{X}\rightarrow[0,+\infty]$ if $I$ is lower-semicontinuous
and for every open set $O$ and closed set $F$ in $\mathcal{X}$,
\begin{equation*}
 -\inf_{x\in O} I(x)\leq\liminf_{\eps\downarrow0}s_\eps^{-1}\log\mathbb{P}(X_\eps\in O),
 \qquad
 \limsup_{\eps\downarrow0}s_\eps^{-1}\log\mathbb{P}(X_\eps\in F)\leq-\inf_{x\in F} I(x).
\end{equation*}
A rate function is said to be \emph{good} if every finite sublevel set is compact. 
We refer to Dembo and Zeitouni (1998) \cite{dembo1998large} for a survey on large deviations theory and its applications.

\subsection{Main large deviations results}
\begin{definition}\label{def:lambda}
For $u\in\RR$, let 
\begin{align}\label{eq:lambda}
\underline{\Jj}(u)
& :=\inf\left\{
\half\|\dot{h}\|^2_{L^2([0,1])}: \dot{h}\in L^2([0,1]),
\Gg(\Kk\dot{h}) = u\right\},\\
\widetilde{\underline{\Jj}}(u)
& :=\inf\left\{
\half\|\dot{h}\|^2_{L^2([0,1])}: \dot{h}\in L^2([0,1]),
\Gg(\Kk\dot{h}) \leq u\right\}.
\end{align}
\end{definition}

The main result is the large deviations statement in the following theorem.

\begin{theorem}\label{thm:LDP}
Under Assumption~\ref{ass:main}, 
the family $(Z_{\eps}/\eps)_{\eps>0}$ satisfies an LDP with speed~$\speed_{\eps}^{-2}$ and good rate function $\underline{\Jj}(\log(\cdot)-\varphi(0))$ on $(0,\infty)$.
\end{theorem}

\begin{proof}
The proof of Theorem~\ref{thm:LDP} will be provided in Appendix~\ref{sec:thm:LDP}.
\end{proof}

The following proposition presents a rephrasing in terms of Laplace transforms and negative moments, given by Varadhan's lemma~\cite[Theorem~4.3.1]{dembo1998large}.
We let
$\underline{\Jj}^{*}(u) :=\sup_{x\in\RR}\{ux-\underline{\Jj}(x)\}$ 
for $u \in \RR$ denote the Fenchel-Legendre transform of~$\underline{\Jj}$.

\begin{proposition}\label{prop:laplace}
Under Assumption~\ref{ass:main},
\begin{equation}\label{eq:mainlaplace}
\lim_{\eps\downarrow 0}\speed_{\eps}^{2} 
\log\EE\left[\left(\frac{Z_{\eps}}{\eps\E^{\varphi(0)}}\right)^{\frac{u}{\speed_{\eps}^{2}}}\right]
 = \underline{\Jj}^{*}(u)\text{ is finite for all }u\in\RR.
\end{equation}
\end{proposition}

\begin{proof}
The proof of Proposition~\ref{prop:laplace} will be provided in Appendix~\ref{sec:thm:LDP}.
\end{proof}

\subsection{Properties of the rate function}\label{sec:tails}

We will repeatedly use the elementary bound
valid for any $t>0$ and $h\in\Hh_{t}$:
\begin{equation}\label{eq:CS}
|(\Kk\dot{h})(t)|\leq \|h\|_{\Hh_{t}}\left(\int_0^t u^{2H-1}\D u\right)^{\half}
 = \|h\|_{\Hh_{t}}
 \frac{t^{H}}{\sqrt{2H}},
\end{equation}
as well as the H\"older regularity inequality
$$
\left|(\Kk\dot{h})(t) - (\Kk\dot{h})(s)\right|\leq \|h\|_{\Hh_{t}}\frac{|t-s|^H}{\sqrt{H}},
$$
which follows from Lemma~\ref{lem:monotone}. 
We introduce the quantities
$$
\Af := \int_{0}^{1}X_s\D s
\qquad\text{and}\qquad \Bf := \int_{0}^{1}X_s^2\D s,
$$
and note that, by stochastic Fubini's theorem,
\begin{equation}\label{eq:AB}
\EE[\Af^2]=\frac{2}{(2H+1)^2(H+1)}
\qquad\text{and}\qquad 
\EE[\Bf]=\frac{1}{2H(2H+1)}.
\end{equation}

The following lemma, proved in Appendix~\ref{sec:lem:prop_lambda_proof}, gathers properties for~$\underline{\Jj}$ and~$\widetilde{\underline{\Jj}}$.

\begin{lemma}\label{lem:prop_lambda}
The following hold:
\begin{enumerate}[(i)]
\item $\underline{\Jj}$ is non-increasing on $(-\infty,0]$, 
non-decreasing on $[0,\infty)$, $\underline{\Jj}(0)=0$, $\underline{\Jj}(u)>0$ for $u\ne0$ and 
$\underline{\Jj}(u)
\geq \frac{u^2}{2 \EE[\Af^2]}$ for every $u\le0$;
\item $\widetilde{\underline{\Jj}}$ is convex on $\RR$, and $\widetilde{\underline{\Jj}}=\underline{\Jj}$ on $(-\infty,0]$;
\item $\underline{\Jj}(u)
= H u^2\left(1+\Oo\left(\frac{\log(u)}{u}\right)\right)$ as $u\uparrow\infty$
 and 
$\underline{\Jj}(u)
=\frac{u^2}{2 \EE[\Af^2]}+o(u^2)$
as $u\to 0$,
and $\underline{\Jj}(u)\sim\Lfbar \E^{2H|u|}$
as $u\downarrow-\infty$, 
for some strictly positive constant~$\Lfbar$.
\end{enumerate}
\end{lemma}

Next, we provide explicit upper and lower bounds for~$\Lfbar$, 
proved in Appendix~\ref{sec:prop:bounds_ch_proof}.

\begin{proposition}\label{prop:bounds_ch}
The following bounds hold for $H\in (0,\half)$:
\begin{equation}\label{eq:upper_lower}
\frac{(H+\half)^2(H+1)}{H^2\E^{2}}
\leq \Lfbar \leq \frac{\Gamma(2H)}{\sin(\pi H)\Gamma(H+\half)^{2}}
\frac{(1-H)^{2(H-1)}}{4H^{2H+1}}.
\end{equation}
\end{proposition}

The following is immediate from the proof of Lemma~\ref{lem:prop_lambda}:
\begin{remark}\label{rem:half}
For $H=\half$ and $\ell\equiv 1$, 
then $\Lfbar=2$.
\end{remark}

Next, we compute the limit of the logarithmic moment generating function:

\begin{proposition}\label{prop:secondorder}
Assume in addition that 
$\sup_{s\leq \eps}|\varphi(s)-\varphi(0)|
 = o(\speed_{\eps}^2)$.
Then 
$$
\lim_{\eps\downarrow 0} \frac{1}{\speed_{\eps}^{2}} \log\EE\left[\left(\frac{Z_{\eps}}{\eps\E^{\varphi(0)}}\right)^{-u}\right]
 = \frac{u(u+1) \EE[\Af^2]-u \EE[\Bf]}{2},
\qquad\text{for any }u\in\RR.
$$
\end{proposition}

\begin{proof}
Write 
$Y_{\eps} 
:= \int_{0}^{1}\E^{\speed_{\eps}X^{\eps}_s + \delta_{\eps}(s)}\D s$,
with $X^{\eps}$ in~\eqref{eq:hatXT}, $\delta_{\eps}:=\varphi_{\eps}$ in~\eqref{eq:phi_notations} then $\|\delta_{\eps}\|_\infty
=\overline{\varphi}_{\eps}
 = o(\speed_{\eps}^2)$.
Expanding in~$\speed_{\eps}$, we obtain
$Y_{\eps}=1+\speed_{\eps}\Af_{\eps} + \half\speed_{\eps}^2\Bf_{\eps} + \Oo(\speed_{\eps}^3)$,
with $\Af_{\eps}:=\int_{0}^{1}X_{s}^{\eps}\D s$ and $\Bf_{\eps}:=\int_{0}^{1}(X_{s}^{\eps})^2\D s$. 
Hence
$\log(Y_{\eps}) = \speed_{\eps}\Af_{\eps} + \frac{\speed_{\eps}^2}{2}(\Bf_{\eps}-\Af_{\eps}^2)
 + \Oo\left(\speed_{\eps}^3\right)$,
and, since $\EE[\Af_{\eps}]=0$, then
$$
\EE\left[Y_{\eps}^{-u}\right]
 = 1 + \frac{\speed_{\eps}^2}{2}
\left(u^2\EE\left[\Af_{\eps}^2\right]
 - u\left(\EE\left[\Bf_{\eps}\right]-\EE\left[\Af_{\eps}^2\right]\right)\right)
 + \Oo\left(\speed_{\eps}^3\right).
$$
Uniform integrability of the remainders follows from $\sup_{\eps}\EE[\exp\{c\|X^{\eps}\|_\infty\}]<\infty$ by Borell-TIS~\cite[Theorem~2.1.1]{adler2007random}, and $\EE[\Af_{\eps}^2]\to\EE[\Af^2]$, $\EE[\Bf_{\eps}]\to\EE[\Bf]$ by Lemma~\ref{lem:coupling}.
\end{proof}


The following corollary is a direct consequence of Theorem~\ref{thm:LDP} and~\eqref{eq:exact}.

\begin{corollary}\label{cor:lefttail}
For every $z\in(0,1)$,
$$
\lim_{\eps\downarrow 0} \speed_{\eps}^2\log\PP\left(\frac{Z_{\eps}}{\eps\E^{\varphi(0)}}\leq z\right)=-\underline{\Jj}(\log(z))=-z^{-2H}\Lf\left(\frac{1}{z}\right).
$$
\end{corollary}

\subsection{Consistency at $H=\half$}

For $H=\half$ the constrained problems defining $\underline{\Jj}$ reduce to the classical Liouville ODE and can be solved explicitly.
\begin{proposition}\label{prop:half}
Let $H=\half$. Then $\underline{\Jj}\in \Cc^{\infty}(\RR)$ is strictly convex and reads
\begin{equation*}
\underline{\Jj}(u) 
 =
 \left\{
\begin{array}{ll}
2u_{+} (u_{+}-\tanh(u_{+})),
& \text{if }u\geq 0,\\
2u_{-} (\tan(u_{-})-u_{-}),
& \text{if }u< 0,
\end{array}
 \right.
\end{equation*}
where~$u_{-}$ and~$u_{+}$ are respectively the unique solutions in $[0,\frac{\pi}{2})$ and $[0,\infty)$ of the equations
$$
u  = \log\left(\frac{\sinh(2u_{+})}{2u_{+}}\right)
\qquad\text{and}\qquad
u = \log\left(\frac{\sin(2u_{-})}{2u_{-}}\right).
$$
Furthermore,
\begin{equation*}
\underline{\Jj}(u)=
\begin{cases}
\frac{3}{2} u^{2}(1 + o(1)), & \text{as }u\to 0,\\
\frac{u^{2}}{2}+u(\log(2u)-1) + \Oo\left(\log(u)^{2}\right), & \text{as } u\uparrow\infty,\\
2 \E^{|u|}-\frac{\pi^{2}}{2}+\Oo\left(\E^{-|u|}\right), & \text{as }u\downarrow-\infty,
\end{cases}
\end{equation*}
which recovers $\Lfbar=2$ from Remark~\ref{rem:half}.
\end{proposition}

\begin{proof}
Let $H=\half$, so that 
$(\Kk\dot{h})(t) = \int_0^t\dot{h}(s)\D s$,
and
$\If(f) = \half\int_{0}^{1}(f'(s))^2\D s$,
for $f \in \Hh_{1}$ (and $\If(f)=+\infty$ otherwise), and the adjoint of~$\Kk$ on $L^2([0,1])$ is 
$$
(\Kk^{*} g)(v) = \int_v^1 g(s)\D s.
$$
For $u>0$, 
the infimum in the definition of~$\underline{\Jj}$ is attained: 
$\If$ is a good rate function, $\Gg$ is continuous on $\Cc_0[0,1]$ 
and $\{\Gg\geq u\}$ is closed. 
At a minimiser, $\Gg_{1}(f)=u$, 
and the constraint gradient 
\[
\nabla\Gg_{1}(f)=\frac{\E^{f}}{\int_{0}^{1}\E^{f(t)}\D t} 
\]
does not vanish.
By Karush-Kuhn-Tucker, there exists $a>0$ such that $\dot{h} = a \Kk^{*}\E^{f}$.
Any KKT point satisfies
$f'' = -a \E^{f}$ on $[0,1]$
with $f(0)=0$ and $f'(1)=0$.
The energy
\[
\half f'^2 + a \E^{f} = a \E^{f(1)} 
\]
yields $f'>0$ on $[0,1)$.
The terminal-value problem at $t=1$ has the unique solution
$$
f(t) = f(1) - 2\log\cosh(\gamma(1-t)),
\qquad\text{with }
\gamma:=\sqrt{\frac{a}{2}\E^{f(1)}}.
$$
The boundary condition $f(0)=0$ forces $\E^{f(1)} = \cosh(\gamma)^{2}$. 
Writing $u_+:=\gamma$, the KKT points are therefore exactly the family
$$
f_{u_+}(t) = 2\log\left(\frac{\cosh(u_+)}{\cosh(u_+(1-t))}\right),
\qquad\text{for } u_+\in(0,\infty)
\qquad\text{and}\qquad
a = \frac{2u_+^{2}}{\cosh^{2}(u_+)}.
$$
Direct computation gives
$$
\int_{0}^{1}\E^{f_{u_+}(x)}\D x
= \cosh^{2}(u_+)\frac{\tanh(u_+)}{u_+}
= \frac{\sinh(2u_+)}{2u_+},
$$
and 
$\If(f_{u_+})
= 2u_+\int_0^{u_+}\tanh^{2}(x)\D x
 = 2u_+(u_+-\tanh(u_+))$.
The rest of the proof follows by direct computations.
For $u<0$ the multiplier changes sign: any KKT point satisfies $f''=a\E^{f}$ on $[0,1]$ with $a>0$, $f(0)=0$ and $f'(1)=0$, so that, with $u_{-}\in(0,\frac{\pi}{2})$,
$$
f_{u_{-}}(t) = 2\log\left(\frac{\cos(u_{-})}{\cos(u_{-}(1-t))}\right),
\qquad
\int_{0}^{1}\E^{f_{u_{-}}(x)}\D x = \frac{\sin(2u_{-})}{2u_{-}},
$$
and $\If(f_{u_{-}}) = 2u_{-}\left(\tan(u_{-})-u_{-}\right)$.
\end{proof}

\begin{remark}\label{rem:PZ}
Proposition~\ref{prop:half} agrees with the solution for the rate function for Asian options in the Black-Scholes model obtained in~\cite{pirjol2016short}. 
More precisely one has $\underline{\Jj}(u) = \Jj_{\mathrm{\BS}}(\E^{u})$ where $\Jj_{\mathrm{\BS}}(x)$ is the rate function defined in~\cite[Proposition 12]{pirjol2016short}. 
The asymptotics $u \to \pm\infty$ also agree with the leading-order terms in~\cite[Proposition~13]{pirjol2016short}. 
The higher-order terms in the expansion of $\underline{\Jj}(u)$ are 
\[
\underline{\Jj}(u) = \frac{3}{2} u^2-\frac{3}{10}u^3+\frac{109}{1400} u^4 + \Oo\left(u^5\right).
\]
\end{remark}


\subsection{Central limit theorem and moderate deviations}
Moderate deviations for option prices were introduced in~\cite{friz2018option} and studied in~\cite{bayer2019short,friz2022short} for rough volatility models
and in~\cite{jacquier2022large,li2017moderate}.
\cite{jacquier2020pathwise} studied them for realised variance in Markovian volatility models.
In our setting, they take the following form.
First, we have the following central limit theorem.

\begin{theorem}\label{thm:CLT}
Under Assumption~\ref{ass:main}, with $\Phi_{\eps}:=\int_{0}^{1}\E^{\varphi(\eps t)-\varphi(0)}\D t$, \begin{equation}\label{eq:CLT}
\frac{1}{\speed_{\eps}}\left(\frac{Z_\eps}{\eps} - \E^{\varphi(0)}\Phi_{\eps}\right)
\text{ converges in distribution to } 
\Nn\left(0,  \E^{2\varphi(0)} \EE\left[\Af^2\right]\right),
\text{ as }\eps\downarrow 0.
\end{equation}
\end{theorem}

\begin{proof}
The proof of Theorem~\ref{thm:CLT} will be given in Appendix~\ref{sec:proofs}.
\end{proof}

Next, we have the following moderate deviations result.

\begin{theorem}\label{thm:MDP}
Let $(a_{\eps})_{\eps}$ be a positive sequence 
such that $\lim_{\eps\downarrow 0}a_{\eps} = \lim_{\eps\downarrow 0} \frac{\speed_{\eps}}{a_{\eps}} = 0$. 
Under Assumption~\ref{ass:main}, the family
$(\frac{1}{a_{\eps}}(\frac{Z_\eps}{\eps} - \E^{\varphi(0)}\Phi_{\eps}))_{\eps>0}$
satisfies, as $\eps\downarrow 0$, a large deviation principle on~$\RR$ with speed $a_{\eps}^2\speed_{\eps}^{-2}$ and good rate function
\begin{equation}\label{eq:J}
\Jj(u) := \frac{u^2 \E^{-2\varphi(0)}}{2 \EE[\Af^2]}
 = \left(H+\half\right)^2(H+1)\E^{-2\varphi(0)} u^2,
\qquad\text{for all }u\in\RR.
\end{equation}
\end{theorem}

\begin{proof}
The proof of Theorem~\ref{thm:MDP} will be given in Appendix~\ref{sec:proofs}.
\end{proof}

\begin{remark}\label{rem:consistency}\ 
\begin{enumerate}[(i)]
\item The slowly varying factor~$\ell$ is invisible in the limit and only enters through~$\speed_{\eps}$.
\item By the Delta method~\cite[Theorem~3.1]{van2000asymptotic},~\eqref{eq:CLT} is equivalent to the weak convergence
$$
\frac{1}{\speed_{\eps}}\left(\log\left(\frac{Z_\eps}{\eps}\right) - \log\left(\E^{\varphi(0)}\Phi_{\eps}\right)\right) \Rightarrow\Nn\left(0,\EE\left[\Af^2\right]\right),
$$
that is $\frac{1}{\speed_{\eps}}(R_\eps - \log(\Phi_{\eps}))\Rightarrow\Nn\left(0,\EE\left[\Af^2\right]\right)$, 
with~$R_{\eps}$ in Section~\ref{sec:thm:LDP}.
\item Typical of moderate deviations~\cite{bryc1997large,djellout1999large,jacquier2020small}, 
$\Jj$ is the quadratic approximation at its minimum of the rate function of Theorem~\ref{thm:LDP};
indeed, for~$z$ close to~$\E^{\varphi(0)}$, Lemma~\ref{lem:prop_lambda}(iii) gives
$$
\underline{\Jj}\left(\log\left(z \E^{-\varphi(0)}\right)\right)
= \frac{\log\left(z \E^{-\varphi(0)}\right)^2}{2\EE[\Af^2]}(1+o(1))\\
= \Jj\left(z-\E^{\varphi(0)}\right)(1+o(1)).
$$
\item Formally at $H=\half$, $\ell\equiv1$, $\varphi\equiv 0$, so that
$\EE[\Af^2]=\frac{1}{3} = \VV[\int_{0}^{1} W_s\D s]$ and~\eqref{eq:CLT} reads 
$$
\frac{1}{\nu\sqrt{\eps}}\left(\frac{1}{\eps}\int_0^\eps \E^{\nu W_t}\D t - 1\right)\Rightarrow\Nn\left(0,\frac{1}{3}\right).
$$
\item This is consistent with~\cite{pirjol2016short}.
Taking $\varphi(t)=\log(S_{0})-\half\nu^{2}t$ 
then~$\frac{Z_{\eps}}{\eps}$ is the arithmetic average over $[0,\eps]$ of the Black-Scholes price $S_{t}=S_{0}\exp\{\nu W_{t}-\half\nu^{2}t\}$.
Since $\sup_{s\leq\eps}|\varphi(s)-\varphi(0)|=\half\nu^{2}\eps=o(\speed_{\eps})$, Corollary~\ref{cor:ATMlognormal} yields
$$
\EE\left[\left(\frac{1}{\eps}\int_{0}^{\eps}S_{t}\D t-S_{0}\right)_{+}\right]
= \frac{\nu S_{0}\sqrt{\eps}}{\sqrt{6\pi}}(1+o(1)),
$$
which is the at-the-money asymptotics of~\cite[Theorem~6]{pirjol2016short} 
(with $\sigma(S_{0})=\nu$).
\end{enumerate}
\end{remark}

\begin{remark}\label{rk:literature}
Small-time CLT for stochastic Volterra equations with H\"older coefficients have been recently obtained by Friesen, Gerhold, Wiedermann~\cite{friesen2026small}, including a functional CLT~\cite[Corollary~2.4]{friesen2026small} for a kernel $K(t) = \Ll(t)t^{H-\half}$ with $\Ll$ locally Lipschitz~\cite[Example~2.8]{friesen2026small} and small-time Gaussian limit of $\frac{1}{\eps}\int_0^\eps v_t\D t$~\cite[Proposition~2.11(iii), Example~2.10]{friesen2026small}, where~$(v_{t})_{t\geq 0}$ is the instantaneous variance process. 
For~$\Ll$ locally Lipschitz and $\varphi\equiv0$, Theorem~\ref{thm:CLT} 
thus follows from~\cite{friesen2026small} applied to $v = \exp(\nu X)$ 
(after localisation of the exponential as in~\cite[Remark~2.6]{friesen2026small}).
For $K(t)=t^{H-\half}$, \cite[Example~2.10]{friesen2026small} gives the limiting kernel $\overline{K}(t)=\sqrt{2H}t^{H-\half}$ and the normalisation $\sqrt{\lambda(n)}\sim\sqrt{2H}n^{H}$ so that~\cite[Remark~2.6, Proposition~2.11(iii)]{friesen2026small} with $f(x)=\E^{\nu x}$ give that
$\frac{1}{\eps^{H}}(\frac{1}{\eps}\int_{0}^{\eps}\E^{\nu X_{t}}\D t-1)$ is asymptotically centred Gaussian with variance
$$
\frac{\nu^{2}}{2H}\int_{0}^{1}\left(\int_{0}^{s}\overline{K}(t)\D t\right)^{2}\D s
= \nu^{2}\EE[\Af^{2}],
$$
using~\eqref{eq:AB}, which is exactly~\eqref{eq:CLT} with $\speed_{\eps}=\nu\eps^{H}$ and $\varphi\equiv 0$.
Compared with~\cite{friesen2026small}, Theorem~\ref{thm:CLT} 
allows for a general slowly varying~$\ell$, 
for instance the log-modulated kernels in~\cite{bayer2021log}, 
for which the tightness argument of~\cite{friesen2026small} is not available 
(as explained in~\cite[Example~2.8]{friesen2026small});
this is precisely the role of the uniform estimates in Lemmas~\ref{lem:holder}-\ref{lem:coupling},
which in our Gaussian case allow to bypass the tightness arguments.
\end{remark}

\section{Short-maturity variance options in rough volatility models}\label{sec:options}


We now consider the following rough volatility model:
\begin{equation}\label{eq:modelS2}
\frac{\D S_{t}}{S_{t}} = \sigma(Y_{t})\D B_{t},
\qquad Y_{t} = y_{0} + W_{t}^{H},
\end{equation}
where~$B$ is a standard Brownian motion and $\sigma:\RR\to[0,\infty)$ a continuous function.
The process~$W^{H}$, 
possibly correlated with~$B$, 
is a fractional Brownian motion with Hurst index $H\in(0,1)$, that is a zero-mean Gaussian process
with covariance function
\begin{equation}
\mathrm{Cov}\left(W_{t}^{H},W_{s}^{H}\right)
 = \half\left(|t|^{2H}+|s|^{2H}-|t-s|^{2H}\right).
\end{equation}
and Volterra representation
$W_{t}^{H} = \int_{0}^{t}K_{H}(s,t)\D W_{s}$,
for some Brownian motion~$W$, with
\begin{align*}
& K_{H}(s, t)
\\
& = \begin{cases}c_{+} s^{\half-H} \int_s^t(u-s)^{H-\frac{3}{2}} u^{H-\half} \D u, & H \in\left(\half, 1\right), \\
c_{H}\left[\left(\frac{t}{s}\right)^{H-\half}(t-s)^{H-\half}-\left(H-\half\right) s^{\half-H} \int_s^t u^{H-\frac{3}{2}}(u-s)^{H-\half} \D u\right], & H \in\left(0, \half\right),
\end{cases}\
\end{align*}
where $c_{+}=\left(\frac{H(2 H-1)}{\mathrm{B}\left(2-2 H, H-\half\right)}\right)^{\half}$, $c_{H}=\left(\frac{2 H}{(1-2 H) \mathrm{B}\left(1-2 H, H+\half\right)}\right)^{\half}$ and $\mathrm{B}(\cdot,\cdot)$ the Beta function, while $K_{\half}\equiv 1$.
Denoting
$$
\Sigma_{T} := \frac{1}{T}\int_{0}^{T}\D\left\langle \log(S),\log(S)\right\rangle_t
= \frac{1}{T}\int_{0}^{T}\sigma^{2}\left(Y_{t}\right)\D t,
$$
the normalised quadratic variation over $[0,T]$, 
the prices of the variance Call and Put options with strike $K>0$ at inception read
$$
\Call(K,T)
= \EE\left[\left(\Sigma_{T} - K\right)_{+}\right]
\qquad\text{and}\qquad
\Put(K,T) = \EE\left[\left(K-\Sigma_{T}\right)_{+}\right].
$$

\begin{remark}
The kernel~$K_{H}$ plays for the fractional Brownian motion the role played by the convolution kernel~$K_{\ell}$ of Assumption~\ref{ass:main} for the Volterra process~$X^{(\ell)}$ of Section~\ref{sec:setup}.
The statements below are formulated for~$K_{H}$; 
in the lognormal case $\sigma^{2}(y)=\E^{\nu y}$, the corresponding results for the Volterra driver of Section~\ref{sec:setup} follow directly from Theorem~\ref{thm:LDP} and Theorem~\ref{thm:CLT}, with explicit rate functions and constants (Corollaries~\ref{cor:OTMlognormal}-\ref{cor:ATMlognormal}).
\end{remark}

Mimicking~\eqref{eq:Kop}-\eqref{eq:I},
we introduce, whenever finite,
\begin{equation}
(\Kk_{H}h)(t)
:= \int_0^{t}K_{H}(u,t)h(u) \D u,\qquad t\in[0,T].
\end{equation}

\subsection{Out-of-the-money asymptotics}\label{sec:OTM}

We first consider the short-maturity asymptotics
for out-of-the-money variance Call ($K>\sigma_{0}^{2}$)
and Put ($K<\sigma_{0}^{2}$) options, where $\sigma_{0}^{2}:=\sigma^{2}(y_{0})$ denotes the spot variance.
We make the following assumption on the volatility function~$\sigma(\cdot)$.

\begin{assumption}\label{assump:1}
The map $\sigma:\RR\to[0,\infty)$
is continuous and there exist $\theta_{1},\theta_{2}>0$ such that
$\sigma^{2}(y)\leq\theta_{1}\E^{\theta_{2}|y|}$
for all $y\in\RR$.
\end{assumption}

We present the leading-order asymptotics for the out-of-the-money variance Call and Put options
in the next theorem.

\begin{theorem}\label{thm:OTM}
For any $K>0$, let
\begin{equation}\label{eq:Lambda_var}
\Lambda(K) := \inf\left\{\half\|h\|_{\Hh_{1}}^2:
\int_{0}^{1}\sigma^{2}\left(y_{0} + (\Kk_{H}\dot{h})(t)\right)\D t=K\right\}.
\end{equation}
Under Assumption~\ref{assump:1}, 
if~$\Lambda$ is right-continuous on $(\sigma_{0}^{2},\infty)$ (for Calls), respectively left-continuous on $(0,\sigma_{0}^{2})$ (for Puts), then
\begin{equation}\label{CP:OTM}
\begin{array}{rll}
\lim_{T\downarrow 0}T^{2H}\log \Call(K,T)
 & = -\Lambda(K), & \text{ if }K>\sigma_{0}^{2},\\
\lim_{T\downarrow 0}T^{2H}\log \Put(K,T)
 & = -\Lambda(K), & \text{ if }K \in (0,\sigma_{0}^{2}).
\end{array}
\end{equation}
\end{theorem}

\begin{proof}
The proof of Theorem~\ref{thm:OTM} is given in Appendix~\ref{sec:proofs:options}.
\end{proof}

The continuity assumption of~$\Lambda$ is guaranteed for example whenever~$\sigma^{2}$ is strictly monotone or continuously differentiable, as shown in Lemma~\ref{lem:Lambda_cont}.
In the lognormal specification, the variational problem~\eqref{eq:Lambda_var} is precisely the one defining the rate function~$\underline{\Jj}$ of Definition~\ref{def:lambda}, and Theorem~\ref{thm:LDP} yields explicit out-of-the-money asymptotics for the Volterra driver of Section~\ref{sec:setup}.
Recall~$Z_{\eps}$ in~\eqref{eq:Z} and~$\speed_{\eps}$ in~\eqref{eq:eps}, and that $Z_{\eps}/\eps$ converges almost surely to $\sigma_{0}^{2}=\E^{\varphi(0)}$, the spot variance.

\begin{corollary}\label{cor:OTMlognormal}
Under Assumption~\ref{ass:main},
\begin{align}
\lim_{\eps\downarrow 0}\speed_{\eps}^{2}\log \EE\left[\left(\frac{Z_{\eps}}{\eps}-K\right)_{+}\right]
&= -\underline{\Jj}\left(\log\left(\frac{K}{\sigma_{0}^{2}}\right)\right),
\qquad\text{for every } K>\sigma_{0}^{2},
\nonumber
\\
\lim_{\eps\downarrow 0}\speed_{\eps}^{2}\log \EE\left[\left(K-\frac{Z_{\eps}}{\eps}\right)_{+}\right]
&= -\underline{\Jj}\left(\log\left(\frac{K}{\sigma_{0}^{2}}\right)\right),
\qquad\text{for every } K \in (0,\sigma_{0}^{2}).
\nonumber
\end{align}
\end{corollary}

\begin{proof}
The proof of Corollary~\ref{cor:OTMlognormal} is given in Appendix~\ref{sec:proofs:options}.
\end{proof}

\begin{remark}\label{rem:OTMlognormal}\
\begin{enumerate}[(i)]
\item Lemma~\ref{lem:prop_lambda}(iii)
and Proposition~\ref{prop:half} translate into option-price asymptotics: writing $u:=\log(\frac{K}{\sigma_{0}^{2}})$ for the log-moneyness of the variance, the decay rate behaves as $\frac{u^{2}}{2\EE[\Af^{2}]}$ near the money 
(consistently with the moderate deviations of Theorem~\ref{thm:MDP}), as $Hu^{2}$ for deep out-of-the-money Calls, and as $\Lfbar\E^{2H|u|}$ for deep out-of-the-money Puts; at $H=\half$ the rate function is fully explicit.
\item For $\ell\equiv1$, $\varphi\equiv0$ and~\eqref{eq:modelS2} with $Y=y_{0}+X^{(1)}$ and $\sigma^{2}(y)=\E^{\nu y}$, $\nu>0$ then
$\int_{0}^{T}\sigma^{2}(Y_{t})\D t=\sigma^{2}(y_{0})Z_{T}$,
so that Corollary~\ref{cor:OTMlognormal} reads
$$
\lim_{T\downarrow 0}T^{2H}\log \EE[(\Sigma_{T}-K)_{+}]
= -\frac{1}{\nu^{2}}\underline{\Jj}\left(\log\left(\frac{K}{\sigma_{0}^{2}}\right)\right),
\qquad \text{for } K>\sigma_{0}^{2},
$$
the exact analogue of Theorem~\ref{thm:OTM} for the Riemann-Liouville kernel, 
now with explicit rate function; out-of-the-money Puts are analogous.
\end{enumerate}
\end{remark}

For $g\in\Cc_{0}([0,1])$, let $I_{H}(g):=\half\|\dot{h}\|^{2}_{L^{2}([0,1])}$ if $g=\Kk_{H}\dot{h}$ for some $\dot{h}\in L^{2}([0,1])$ and $I_{H}(g):=+\infty$ otherwise, and
$F(g) := \int_{0}^{1}\sigma^{2}\left(y_{0}+g(t)\right)\D t$,
so that $F$ is continuous on $\Cc_{0}([0,1])$, $F(0)=\sigma_{0}^{2}$ and
$\Lambda(K) = \inf\{I_{H}(g): F(g)=K\}$ in~\eqref{eq:Lambda_var}, with $\inf\emptyset=+\infty$.

\begin{lemma}\label{lem:Lambda_cont}
Let $\sigma$ be continuous.
\begin{enumerate}[(i)]
\item For $K>\sigma_{0}^{2}$, $\Lambda(K) = \inf\{I_{H}(g): F(g)\geq K\}$, and for $K<\sigma_{0}^{2}$, $\Lambda(K) = \inf\{I_{H}(g): F(g)\leq K\}$.
In particular, $\Lambda$ is non-decreasing on $(\sigma_{0}^{2},\infty)$ and non-increasing on $(0,\sigma_{0}^{2})$.
\item If~$\sigma^{2}$ is strictly monotone or continuously differentiable, then~$\Lambda$ is right-continuous on $(\sigma_{0}^{2},\infty)$ and left-continuous on $(0,\sigma_{0}^{2})$.
\end{enumerate}
\end{lemma}

\begin{proof}
(i) Let $K>\sigma_{0}^{2}$ and~$g$ such that $F(g)\geq K$ and $I_{H}(g)<\infty$.
The map $\theta\mapsto F(\theta g)$ is continuous on $[0,1]$, equal to $\sigma_{0}^{2}<K$ at $\theta=0$ and to $F(g)\geq K$ at $\theta=1$, so that $F(\theta g)=K$ for some $\theta\in(0,1]$, and $I_{H}(\theta g)=\theta^{2}I_{H}(g)\leq I_{H}(g)$.
Hence 
\[
\inf\left\{I_{H}(g):F(g)\geq K\right\}\geq\Lambda(K),
\]
and the reverse inequality is trivial.
The case $K<\sigma_{0}^{2}$ is identical, and monotonicity follows since the sets $\{F\geq K\}$ and $\{F\leq K\}$ are nested.

(ii) Assume that~$\sigma^{2}$ is strictly increasing (otherwise replace~$\mathfrak{c}$ by~$-\mathfrak{c}$), 
and let
$$
\mathfrak{c}(t) := \EE\left[W^{H}_{t}\int_{0}^{1}W^{H}_{s}\D s\right]
= \half\int_{0}^{1}\left(t^{2H}+s^{2H}-|t-s|^{2H}\right)\D s,
\qquad \text{for all } t\in[0,1],
$$
strictly positive on $(0,1]$,
which belongs to the Cameron-Martin space of~$W^{H}$, with
$$
I_{H}(\mathfrak{c})
= \half\VV\left[\int_{0}^{1}W^{H}_{s}\D s\right]
= \frac{1}{4(H+1)}.
$$
Let $K>\sigma_{0}^{2}$ with $\Lambda(K)<\infty$ (otherwise there is nothing to prove by~(i)), $\varsigma>0$, and~$g$ such that $F(g)=K$ and $I_{H}(g)\leq\Lambda(K)+\varsigma$.
For $\delta>0$, strict monotonicity yields $F(g+\delta \mathfrak{c})>K$, and $\delta'\mapsto F(g+\delta' \mathfrak{c})$ is continuous.
Hence for every $K'\in(K,F(g+\delta \mathfrak{c})]$ there exists $\delta'\in(0,\delta]$ such that $F(g+\delta' \mathfrak{c})=K'$, and therefore
$$
\Lambda(K')\leq I_{H}(g+\delta' \mathfrak{c})
\leq \left(\sqrt{\Lambda(K)+\varsigma}+\delta\sqrt{I_{H}(\mathfrak{c})}\right)^{2}.
$$
Letting~$K'\downarrow K$, 
then $\delta,\varsigma\downarrow 0$, 
hence $\limsup_{K'\downarrow K}\Lambda(K')\leq\Lambda(K)$ and right-continuity follows from~(i).
For $K<\sigma_{0}^{2}$, the same argument with $g-\delta \mathfrak{c}$ gives left-continuity.

Assume instead that~$\sigma^{2}$ is continuously differentiable and let $K>\sigma_{0}^{2}$ with $\Lambda(K)<\infty$.
Since~$I_{H}$ is a good rate function and $\{F\geq K\}$ is closed, the infimum in~(i) is attained at some~$g$, with $I_{H}(g)=\Lambda(K)>0$ because $F(0)=\sigma_{0}^{2}<K$, and $F(g)=K$ (otherwise $F(\theta g)=K$ for some $\theta\in(0,1)$, with $I_{H}(\theta g)<I_{H}(g)$).
For~$\widehat{\mathfrak{c}}$ in the Cameron-Martin space of~$W^{H}$, dominated convergence theorem gives
$$
\lim_{\delta\downarrow 0}\frac{F(g+\delta\widehat{\mathfrak{c}})-F(g)}{\delta}
= \int_{0}^{1}(\sigma^{2})'\left(y_{0}+g(t)\right)\widehat{\mathfrak{c}}(t)\D t.
$$
Since this space contains all smooth functions with compact support in $(0,1)$~\cite{decreusefond1999stochastic}, if the right-hand side vanishes for all of them, then $(\sigma^{2})'(y_{0}+g(t))=0$ for all $t\in[0,1]$, so that~$\sigma^{2}$ is constant on the interval $y_{0}+g([0,1])$ and $F(g)=\sigma_{0}^{2}\neq K$.
Hence there exists~$\widehat{\mathfrak{c}}$ such that $F(g+\delta\widehat{\mathfrak{c}})>K$ for all $\delta>0$ small enough, and the argument above, replacing~$\mathfrak{c}$ by~$\widehat{\mathfrak{c}}$ 
and taking $\varsigma=0$ gives right-continuity at~$K$; 
left-continuity for $K<\sigma_{0}^{2}$ follows by analogous arguments.
\end{proof}

\subsection{The variational problem}\label{sec:variational}

We now solve the variational problem~\eqref{eq:Lambda_var} 
for OTM variance Call and Put options in Theorem~\ref{thm:OTM},
assuming $\sigma(\cdot)$ continuous and almost everywhere differentiable.
The following is proved in Appendix~\ref{sec:proofs:options}.

\begin{proposition}\label{prop:EL}
The optimal~$h$ for~\eqref{eq:Lambda_var} 
satisfies the integral equation
\begin{equation}\label{eq:Integ_equation}
\dot{h}(t) + \lambda\int_{t}^{1}(\sigma^{2})'\left(y_{0}+(\Kk_{H}\dot{h})(s)\right)K_{H}(t,s)\D s=0,
\qquad\text{for any }0\leq t\leq 1,
\end{equation}
with the constraint
$\int_{0}^{1}\sigma^{2}(y_{0}+(\Kk_{H}\dot{h})(t))\D t=K$,
where~$\lambda\in\RR$ is a Lagrange multiplier, which exists as soon as $(\sigma^{2})'(y_{0}+(\Kk_{H}\dot{h})(\cdot))$ does not vanish identically on~$[0,1]$.
\end{proposition}

\subsection{Special cases}

\subsubsection{Special Case: $\sigma(y)=\xi\sqrt{|y|}$.}

\begin{corollary}\label{cor:special_1}
If $\sigma(y)=\xi\sqrt{|y|}$ for  $\xi>0$,
then for all $y_{0}\ne 0$,
the rate function reads
\begin{equation}\label{eq:special_1_rate}
\half\int_{0}^{1}\dot{h}(t)^{2}\D t
 = \frac{\left(\xi^{2}|y_{0}| - K\right)^{2}\int_{0}^{1}
 \left|\int_{t}^{1}K_{H}(t,s)\D s\right|^{2}\D t}{2\xi^{4}
 \left|\int_{0}^{1}\int_{0}^{t}\int_{s}^{1}K_{H}(s,t)K_{H}(s,u)\D u \D s \D t\right|^{2}}.
\end{equation}
\end{corollary}

\begin{remark}\label{rem:special}
If $\sgn(y_{0}+\Kk_{H}\dot{h}) = \sgn(y_{0})$ on~$[0,1]$, 
then
the triple integral is equal to 
$\VV[\int_{0}^{1}W^{H}_{t}\D t] = \frac{1}{2(H+1)}$ 
(as in the proof of Theorem~\ref{thm:ATM}), 
so that
$$
\half\int_{0}^{1}\dot{h}(t)^{2}\D t = \frac{(H+1)\left(K-\xi^{2}|y_{0}|\right)^{2}}{\xi^{4}}.
$$
Note that
\[
(\Kk_{H}\dot{h})(t) = 2(H+1) \mathrm{sgn}(y_{0})\left(\frac{K}{\xi^{2}}-|y_{0}|\right)\Cc(t),
\]
with $\Cc(t):=\int_{0}^{1}\mathrm{Cov}(W^{H}_{t},W^{H}_{s})\D s\geq 0$,
so this always holds for Calls ($K>\xi^{2}|y_{0}|$), whereas for Puts it requires
$2(H+1)(|y_{0}|-\frac{K}{\xi^{2}})\max_{[0,1]}\Cc<|y_{0}|$ and hence
deeper out-of-the-money Puts are not covered here.
Although $\sigma^{2}(y)=\xi^{2}|y|$ is neither monotone nor continuously differentiable, the continuity condition of Theorem~\ref{thm:OTM} holds 
nevertheless: 
on $(\sigma_{0}^{2},\infty)$, $\Lambda$ is given by the continuous formula above below, and on $(0,\sigma_{0}^{2})$ it is the value function $K\mapsto\inf\{I_{H}(g):F(g)\leq K\}$ of a convex problem (Lemma~\ref{lem:Lambda_cont}(i), with~$F$ convex here), hence convex, finite and thus continuous.

For Call options, let $g:=\Kk_{H}\dot{h}$, $\flat:=\mathrm{sgn}(g)$. 
Then
\[
\int_{0}^{1}|g(t)|\D t = \left\langle\Kk_{H}^{*}\flat,\dot{h}\right\rangle\leq\|\Kk_{H}^{*}\flat\|\|\dot{h}\|, 
\quad\text{in $L^{2}([0,1])$}, 
\]
with
\begin{align*}
\|\Kk_{H}^{*}\flat\|^{2}
&= \int_{0}^{1}\int_{0}^{1}\flat(s)\flat(t)\mathrm{Cov}\left(W^{H}_{s},W^{H}_{t}\right)\D s\D t
\\
&\leq \int_{0}^{1}\int_{0}^{1}\mathrm{Cov}\left(W^{H}_{s},W^{H}_{t}\right)\D s\D t = \frac{1}{2(H+1)}.
\end{align*}
Since $F(g)\leq\xi^{2}(|y_{0}|+\int_{0}^{1}|g(t)|\D t)$, any~$g$ with $F(g)=K>\xi^{2}|y_{0}|$ satisfies
\[
I_{H}(g)\geq(H+1)\left(\int_{0}^{1}|g(t)|\D t\right)^{2}\geq(H+1)\left(\frac{K}{\xi^{2}}-|y_{0}|\right)^{2},
\]
and this lower bound is attained by the above candidate.
Hence $\Lambda(K)=\frac{(H+1)(K-\xi^{2}|y_{0}|)^{2}}{\xi^{4}}$ for all $K>\xi^{2}|y_{0}|$ 
with no assumption on the sign of the optimal path.

For Put options, $F$ is convex and $\Lambda(K)=\inf\{I_{H}(g):F(g)\leq K\}$ by Lemma~\ref{lem:Lambda_cont}(i).
When the candidate~$g^{*}$ above keeps the sign of~$y_{0}$ on~$[0,1]$, $F$ is differentiable at~$g^{*}$ and~$g^{*}$ satisfies the first-order condition~\eqref{eq:Integ_equation} with $\lambda>0$, so that, by convexity, it is the global minimiser and the formula above is equal to~$\Lambda(K)$.
Explicitly,
\[
\Cc(t)=\half\left(t^{2H}+\frac{1-t^{2H+1}-(1-t)^{2H+1}}{2H+1}\right);
\]
for $H\geq\half$, $\Cc$ is non-decreasing with $\max_{t\in[0,1]}\Cc(t)=\Cc(1)=\half$, and the condition above reads $K>\frac{H}{H+1}\xi^{2}|y_{0}|$, whereas for $H<\half$ the maximum is attained in the interior of~$[0,1]$.
\end{remark}

\begin{proof}[Proof of Corollary~\ref{cor:special_1}]
In this case,
the optimal~$h$ should (Remark~\ref{rem:special})
satisfy $y_{0}+\int_{0}^{t}K_{H}(s,t)\dot{h}(s)\D s>0$ for any $0\leq t\leq 1$
if $y_{0}>0$ and 
$y_{0}+\int_{0}^{t}K_{H}(s,t)\dot{h}(s)\D s<0$ for any $0\leq t\leq 1$
if $y_{0}<0$. 
For $y_{0}>0$,~\eqref{eq:Integ_equation} reads
$\dot{h}(t)+\lambda\xi^{2}\int_{t}^{1}K_{H}(t,s)\D s = 0$,
with the constraint
$\xi^{2}\int_{0}^{1}\left(y_{0}+(\Kk_{H}\dot{h})(t)\right)\D t = K$,
so that
\begin{equation}
\xi^{2}y_{0}
-\xi^{2}\int_{0}^{1}\int_{0}^{t}K_{H}(s,t)\lambda\xi^{2}\int_{s}^{1}K_{H}(s,u) \D u \D s \D t = K,
\end{equation}
which implies that
\begin{equation}\label{eq:special_1_lambda}
\lambda
=\frac{\xi^{2}|y_{0}| - K}{\xi^{4}\int_{0}^{1}\int_{0}^{t}\int_{s}^{1}K_{H}(s,t)K_{H}(s,u)\D u \D s \D t},
\end{equation}
and a direct computation gives the rate function in~\eqref{eq:special_1_rate}.

In the case where $y_{0}<0$, then~\eqref{eq:Integ_equation} reads
$\dot{h}(t) = \lambda\xi^{2}\int_{t}^{1}K_{H}(t,s)\D s$
with the constraint
$\xi^{2}\int_{0}^{1}\left(y_{0}+(\Kk_{H}\dot{h})(t)\right)\D t=-K$,
so that
\[
\xi^{2}y_{0}
 + \xi^{2}\int_{0}^{1}\int_{0}^{t}K_{H}(s,t)\lambda\xi^{2}\int_{s}^{1}K_{H}(s,u)\D u \D s \D t=-K,
\]
and the corollary follows.
\end{proof}

\subsubsection{Special Case: $H=\half$.}

\begin{lemma}\label{lem:monotone_path}
Let $H=\half$ and assume~$\sigma^{2}$ monotone.
Then, for every $K\neq\sigma_{0}^{2}$ with~$\Lambda(K)$ finite, the infimum in~\eqref{eq:Lambda_var} is attained at a monotone path 
and the rate function reads
\begin{align*}
\Lambda(K) 
&= \half\int_{0}^{1} \dot{h}(t)^{2}\D t
 \\
&= \lambda\left(K-\sigma(y_{0}+h(1))^{2}\right)
 =  \half\left(\int_{0}^{h(1)}\frac{\D x}{\sqrt{\Delta(x)}}\right)
\left(\int_{0}^{h(1)}\sqrt{\Delta(x)}\D x\right),
\end{align*}
with $\Delta(x):=|\sigma(y_{0}+x)^{2}-\sigma(y_{0}+h(1))^{2}|$, 
and $h(1)$ is the solution of the equation
\begin{equation}\label{h1:eq:main}
\frac{\int_0^{h(1)} \sigma^2(y_{0}+h) \frac{\D h}{\sqrt{\Delta(h)}}}
{\int_0^{h(1)}\frac{\D h}{\sqrt{\Delta(h)}}} = K.
\end{equation}
\end{lemma}

\begin{proof}
Assume that~$\sigma^{2}$ is non-decreasing and that $K>\sigma_{0}^{2}$; the other cases are identical, replacing the running maximum below by the running minimum, or~$h$ by~$-h$.
Let $h\in\Hh_{1}$ with $F(h)=K$ and $M(t):=\max_{s\leq t}h(s)$.
Then~$M$ is non-decreasing, $M(0)=0$, $M\geq h$, and $0\leq M(t)-M(s)\leq\int_{s}^{t}|\dot{h}(u)|\D u$ for $s\leq t$, so that $M\in\Hh_{1}$ with $\|M\|_{\Hh_{1}}\leq\|h\|_{\Hh_{1}}$, and $F(M)\geq F(h)=K$ since~$\sigma^{2}$ is non-decreasing.
As in the proof of Lemma~\ref{lem:Lambda_cont}(i), $F(\theta M)=K$ for some $\theta\in(0,1]$, and~$\theta M$ is non-decreasing with $I_{\half}(\theta M)\leq I_{\half}(h)$.
The infimum may therefore be restricted to non-decreasing paths, which form a closed subset of~$\Cc_{0}([0,1])$, and it is attained there since~$I_{\half}$ is a good rate function and~$F$ is continuous.

We now prove the form of the rate function. Since it is characterised as
\begin{equation}
\half\inf\left\{\|h\|_{\Hh_{1}}^{2}, h(0) = 0, \int_{0}^{1}\sigma(y_{0}+h(t))^{2}\D t=K\right\},
\end{equation}
define the Lagrange functional
\begin{equation}
\Lambda[h]:=\half \|h\|_{\Hh_{1}}^{2} + \lambda\left(\int_{0}^{1}\sigma(y_{0}+h(t))^{2}\D t-K\right),
\end{equation}
where $\lambda$ is a Lagrange multiplier.
The Euler-Lagrange equation is given by
\begin{equation}
\frac{\D}{\D t}\dot{h}(t)
=\lambda(\sigma^{2})'(y_{0}+h(t)),
\end{equation}
with boundary conditions
$\{h(0)=0,\dot{h}(1)=0\}$,
that is
\begin{equation}
\half\frac{\D}{\D t}\dot{h}(t)^{2}
=\frac{\D}{\D t}\lambda\sigma^{2}(y_{0}+h(t)).
\end{equation}
Integrating from~$1$ to~$t$ and applying
the transversality condition $\dot{h}(1)=0$, we obtain
\begin{equation}
\half \dot{h}(t)^{2}
=\lambda\sigma^{2}(y_{0}+h(t)) - \lambda\sigma^{2}(y_{0}+h(1)),
\end{equation}
and therefore
\begin{equation}
\dot{h}(t) = \pm \sqrt{2|\lambda|}\sqrt{|\sigma(y_{0}+h(t))^{2} - \sigma(y_{0}+h(1))^{2}|},
\end{equation}
where the sign~$\pm$ is the sign of~$\dot{h}$, namely~$+$ for an increasing optimal path and~$-$ for a decreasing one:
when~$\sigma^{2}$ is increasing, $+$ corresponds to Calls and~$-$ to Puts, and conversely when~$\sigma^{2}$ is decreasing. 

In the~$+$ case then, with $h(0)=0$,
\begin{equation}
\frac{1}{\sqrt{2|\lambda|}}\int_{0}^{h(t)}
\frac{\D x}{\sqrt{|\sigma(y_{0}+x)^{2} - \sigma(y_{0}+h(1))^{2}|}} = t,
\end{equation}
and in particular
\begin{equation}\label{lambda:eq}
\frac{1}{\sqrt{2|\lambda|}}\int_{0}^{h(1)}\frac{\D x}{\sqrt{\sigma(y_{0}+x)^{2} - \sigma(y_{0}+h(1))^{2}}}=1.
\end{equation}
We assume here and in the following $(\sigma^{2})'(y_{0}+h(1))\neq 0$, which ensures convergence of the integral.
The function~$h(t)$ satisfies the constraint 
$\int_{0}^{1}\sigma(y_{0}+h(t))^{2}\D t=K$. 
This can be used to determine $h(1)$ by writing
\begin{align}\label{h1:eq}
\int_{0}^{1}\sigma(y_{0}+h(t))^{2}\D t 
&= \int_0^{h(1)} \sigma(y_0 + h(t))^2 \frac{\D h}{h'(t)} \\
& = \frac{1}{\sqrt{2|\lambda |}}
\int_0^{h(1)} \sigma(y_0 + x)^2 \frac{\D x}{\sqrt{\Delta(x)}} \\
& = \left( \int_0^{h(1)} \frac{\D x}{\sqrt{\Delta(x)}}\right)^{-1}
\int_0^{h(1)} \sigma(y_0 + x)^2 \frac{\D x}{\sqrt{\Delta(x)}} = K,
\end{align}
where we used~\eqref{lambda:eq} in the last line. This gives an equation for $h(1)$, see~\eqref{h1:eq:main}. 

Finally, the rate function is given by
\begin{equation}\label{LamK}
\Lambda(K)=\half\int_{0}^{1}\dot{h}(t)^{2}\D t
=\lambda\int_{0}^{1}\sigma(y_{0}+h(t))^{2}\D t-\lambda\sigma(y_{0}+h(1))^{2}
=\lambda\left(K-\sigma(y_{0}+h(1))^{2}\right).
\end{equation}
In the~$-$ case, analogous computations yield the rate function.
For both cases, from~\eqref{lambda:eq},
$$
\lambda = \frac{\sgn\left(K-\sigma(y_{0}+h(1))^{2}\right)}{2}\left(\int_{0}^{h(1)}\frac{\D x}{\sqrt{\Delta(x)}}\right)^{2} .
$$
Substituting this and $K$ from \eqref{h1:eq} into \eqref{LamK} gives the stated result for the rate function.
\end{proof}

\begin{remark}\label{rem:halfconsistency}
For the lognormal specification $\sigma(y)^{2} = \E^{\nu y}$, the Euler-Lagrange equation $\ddot{h}(t)=\lambda(\sigma^{2})'(y_{0}+h(t))$ above becomes, after the substitution $f=\nu h$, a Liouville equation solved explicitly in the proof of Proposition~\ref{prop:half}, and the two computations agree.
\end{remark}

\subsection{At-the-money asymptotics}\label{sec:ATM}

We finally consider the short-maturity asymptotics of at-the-money variance Call and Put options, at $K=\sigma_{0}^{2} = \sigma(y_{0})^{2}$.
The following assumption on~$\sigma$ is similar to Assumption~\ref{assump:1} but requires additional smoothness.

\begin{assumption}\label{assump:2}
$\sigma\in\Cc^2(\RR\to[0,\infty))$ and there exist $\theta_{1},\theta_{2}>0$ such that
$\left|(\sigma^{2})''(y)\right|
\leq\theta_{1}\E^{\theta_{2}|y|}$
for all $y\in\RR$.
\end{assumption}

In the following theorem,
proved in Appendix~\ref{sec:proofs:options},
we present the leading-order asymptotics for at-the-money variance Call and Put options.

\begin{theorem}\label{thm:ATM}
Under Assumption~\ref{assump:2},  
\begin{equation}
\lim_{T\downarrow 0}\frac{\Call(\sigma_{0}^{2},T)}{T^{H}}
=\lim_{T\downarrow 0}\frac{\Put(\sigma_{0}^{2},T)}{T^{H}}
=\frac{\left|(\sigma^{2})'(y_{0})\right|}{2\sqrt{\pi(H+1)}}.
\end{equation}
\end{theorem}

\begin{remark}\label{rem:ATMdegenerate}
If $(\sigma^{2})'(y_{0})=0$, the limit in Theorem~\ref{thm:ATM} vanishes and the leading order changes.
In this case, under Assumption~\ref{assump:2},
$$
\lim_{T\downarrow 0}\frac{\Call(\sigma_{0}^{2},T)}{T^{2H}} = \frac{\left((\sigma^{2})''(y_{0})\right)_{+}}{2(2H+1)}
\qquad\text{and}\qquad
\lim_{T\downarrow 0}\frac{\Put(\sigma_{0}^{2},T)}{T^{2H}} = \frac{\left((\sigma^{2})''(y_{0})\right)_{-}}{2(2H+1)}.
$$
Indeed, $\Sigma_{T}-\sigma_{0}^{2}=\frac{1}{2T}\int_{0}^{T}(\sigma^{2})''(R_{t})(W^{H}_{t})^{2}\D t$ with $|R_{t}-y_{0}|\leq|W^{H}_{t}|$.
Since~$(\sigma^{2})''$ is continuous, Cauchy-Schwarz inequality, and the exponential bound in Assumption~\ref{assump:2}  give
$\EE[\frac{1}{T}\int_{0}^{T}|(\sigma^{2})''(R_{t})-(\sigma^{2})''(y_{0})|(W^{H}_{t})^{2}\D t] = o(T^{2H})$,
while $\EE[\frac{1}{T}\int_{0}^{T}(W^{H}_{t})^{2}\D t] = \frac{T^{2H}}{2H+1}$.
\end{remark}

\begin{corollary}\label{cor:ATMlognormal}
Under Assumption~\ref{ass:main},
with~$\Phi_{\eps}$ as in Theorem~\ref{thm:CLT},
\begin{align}
\lim_{\eps\downarrow 0}\frac{1}{\speed_{\eps}}
\EE\left[\left(\frac{Z_{\eps}}{\eps}-\E^{\varphi(0)}\Phi_{\eps}\right)_{+}\right]
&=\lim_{\eps\downarrow 0}\frac{1}{\speed_{\eps}}
\EE\left[\left(\E^{\varphi(0)}\Phi_{\eps}-\frac{Z_{\eps}}{\eps}\right)_{+}\right]
\nonumber
\\
&=\E^{\varphi(0)}\sqrt{\frac{\EE[\Af^{2}]}{2\pi}}
=\frac{\E^{\varphi(0)}}{(2H+1)\sqrt{\pi(H+1)}}.
\nonumber
\end{align}
If moreover $\sup_{s\leq\eps}|\varphi(s)-\varphi(0)|=o(\speed_{\eps})$, then $\E^{\varphi(0)}\Phi_{\eps}$ may be replaced by~$\E^{\varphi(0)}$.
\end{corollary}

\begin{proof}
With  
$Y_{\eps}
 := Z_\eps / (\eps\E^{\varphi(0)})$ then
$$
\EE\left[\left(\frac{Z_{\eps}}{\eps}-\E^{\varphi(0)}\Phi_{\eps}\right)_{+}\right]=\E^{\varphi(0)}\EE[(Y_{\eps}-\Phi_{\eps})_{+}], 
$$
and $x_{+}=\half(x+|x|)$ for all $x\in\RR$.
The proof of Theorem~\ref{thm:CLT} shows that $\speed_{\eps}^{-1}(Y_{\eps}-\Phi_{\eps})$ converges in distribution to $\zeta\sim\Nn(0,\EE[\Af^{2}])$, and that the family $\{|\speed_{\eps}^{-1}(Y_{\eps}-\Phi_{\eps})|^{p}\}_{\eps\leq\eps_{0}}$ is uniformly integrable for every~$p\geq 1$; in particular both $\EE[\speed_{\eps}^{-1}(Y_{\eps}-\Phi_{\eps})]\to\EE[\zeta]=0$ and $\EE[|\speed_{\eps}^{-1}(Y_{\eps}-\Phi_{\eps})|]\to\EE[|\zeta|]$.
Therefore
$$
\lim_{\eps\downarrow 0}
\frac{1}{\speed_{\eps}}\EE\left[(Y_{\eps}-\Phi_{\eps})_{+}\right]
 = \half\lim_{\eps\downarrow 0}
\left\{\EE\left[\frac{Y_{\eps}-\Phi_{\eps}}{\speed_{\eps}}\right]+\EE\left[\left|\frac{Y_{\eps}-\Phi_{\eps}}{\speed_{\eps}}\right|\right]\right\}
= \frac{\EE[|\zeta|]}{2}
 = \sqrt{\frac{\EE[\Af^{2}]}{2\pi}},
$$
from~\eqref{eq:AB}.
The Put statement is identical since 
$(-x)_{+}=\frac{|x|-x}{2}$.
Finally, if $\sup_{s\leq\eps}|\varphi(s)-\varphi(0)|=o(\speed_{\eps})$, then $|\Phi_{\eps}-1|\leq \E^{\overline{\varphi}_{\eps}}-1=o(\speed_{\eps})$ and, since $x\mapsto x_{+}$ is $1$-Lipschitz,
$$
\left|\EE\left[\left(\frac{Z_{\eps}}{\eps}-\E^{\varphi(0)}\Phi_{\eps}\right)_{+}\right]
-\EE\left[\left(\frac{Z_{\eps}}{\eps}-\E^{\varphi(0)}\right)_{+}\right]\right|
\leq \E^{\varphi(0)}\left|\Phi_{\eps}-1\right| = o(\speed_{\eps}),
$$
so the centring may indeed be replaced by~$\E^{\varphi(0)}$.
\end{proof}

\begin{remark}\label{rem:ATMcompare}
At $H=\half$ the constants of Theorem~\ref{thm:ATM} (for lognormal~$\sigma(\cdot)$) and Corollary~\ref{cor:ATMlognormal} coincide; for $H\neq\half$ they differ (a factor~$2$ versus~$2H+1$), reflecting the different kernels (Molchan-Golosov versus Riemann-Liouville) through the variance of the time-integrated driver, $\frac{1}{2(H+1)}$ versus~$\EE[\Af^{2}]$ in~\eqref{eq:AB}.
Part of this discrepancy is a matter of normalisation, since $\VV[X_{1}]=\frac{1}{2H}$ whereas $\VV[W^{H}_{1}]=1$: for the unit-variance Riemann-Liouville kernel $\sqrt{2H}(t-s)^{H-\half}$ of Section~\ref{sec:expansion}, the ratio of the two constants is $\frac{\sqrt{8H}}{2H+1}$, which is equal to one only at $H=\half$.
In both cases, at-the-money variance options vanish at the rate~$T^{H}$ (up to the slowly varying factor~$\ell$), in contrast with the~$\sqrt{T}$ decay of at-the-money European options written on the asset itself.
\end{remark}

The following proposition shows that the results hold with very minor modifications when changing the kernel from~$K_{H}$ to~$K$.

\begin{proposition}\label{prop:unified}
Let Assumption~\ref{ass:main} hold and take~\eqref{eq:modelS2} with $Y=y_{0}+X^{(\ell)}$ instead of $Y=y_{0}+W^{H}$.
Let $\widetilde{\speed}_{T}:=T^{H}\ell(T)$ and~$\Lambda_{X}$ defined as~$\Lambda$ in~\eqref{eq:Lambda_var} with~$\Kk$ instead of~$\Kk_{H}$.
\begin{enumerate}[(i)]
\item Under Assumption~\ref{assump:1}, if~$\sigma^{2}$ is strictly monotone or continuously differentiable, then~\eqref{CP:OTM} holds with $\widetilde{\speed}_{T}^{2}$ in place of~$T^{2H}$ and~$\Lambda_{X}$ in place of~$\Lambda$.
\item Under Assumption~\ref{assump:2},
$$
\lim_{T\downarrow 0}\frac{\Call(\sigma_{0}^{2},T)}{\widetilde{\speed}_{T}}
=\lim_{T\downarrow 0}\frac{\Put(\sigma_{0}^{2},T)}{\widetilde{\speed}_{T}}
=\left|(\sigma^{2})'(y_{0})\right|\sqrt{\frac{\EE[\Af^{2}]}{2\pi}}
=\frac{\left|(\sigma^{2})'(y_{0})\right|}{(2H+1)\sqrt{\pi(H+1)}}.
$$
\end{enumerate}
\end{proposition}
\begin{proof}
(i) With~$X^{T}$ as in~\eqref{eq:hatXT}, 
taking $\eps=T$, $Y_{tT}=y_{0} + \widetilde{\speed}_{T}X^{T}_{t}$ for $t\in[0,1]$ and Theorem~\ref{thm:pathLDP} 
(with $\nu=1$) provides an LDP for $(\widetilde{\speed}_{T}X^{T})_{T>0}$ with speed~$\widetilde{\speed}_{T}^{-2}$ and good rate function~$\If$.
The moment bounds remain valid since $\sup_{t\leq 1}\VV[X^{(\ell)}_{t}]=\VV[X^{(\ell)}_{1}]<\infty$, and so does Lemma~\ref{lem:Lambda_cont}, with $\mathfrak{c}:=\Kk\Kk^{*}\ind=\EE[X_{\cdot}\Af]$, which is positive on~$(0,1]$.
Regarding~(ii), in the proof of Theorem~\ref{thm:ATM}, the first-order term becomes
$$
(\sigma^{2})'(y_{0})\frac{1}{T}\int_{0}^{T}X^{(\ell)}_{t}\D t=(\sigma^{2})'(y_{0})\widetilde{\speed}_{T}\int_{0}^{1}X^{T}_{s}\D s,
$$
which is centred Gaussian, and $\VV[\int_{0}^{1}X^{T}_{s}\D s]$ converges to~$\EE[\Af^{2}]$ as $T\downarrow 0$ by Lemma~\ref{lem:coupling}.
The remainder is bounded by a constant times
$\sup_{t\leq T}\VV[X^{(\ell)}_{t}]=\VV[X^{(\ell)}_{T}]\sim\frac{\widetilde{\speed}_{T}^{2}}{2H}=o(\widetilde{\speed}_{T})$.
The last equality is exactly~\eqref{eq:AB}.
\end{proof}
For $\sigma^{2}(y)=\E^{\varphi(0)+\nu y}$ and $y_{0}=0$, Proposition~\ref{prop:unified} is consistent with Corollaries~\ref{cor:OTMlognormal} and~\ref{cor:ATMlognormal} (with~$\varphi$ constant), since then $\speed_{T}=\nu\widetilde{\speed}_{T}$ and $(\sigma^{2})'(y_{0})=\nu\E^{\varphi(0)}$.

\subsection{Rate function expansion}
\label{sec:expansion}

We now restrict ourselves to rough volatility models driven by a Riemann-Liouville process $X_t  = \int_0^t K_{H}(s,t) \D Z_s$ with kernel 
$K_{H}(s,t)=\sqrt{2H} (t-s)^{H-\half}$ and such that $\VV[X_t]=t^{2H}$.
(Here and until the end of Section~\ref{sec:impliedvol}, $K_{H}$ and~$\Kk_{H}$ refer to this Riemann-Liouville kernel and not to the Molchan-Golosov kernel of Section~\ref{sec:options}, and~$Z$ is a standard Brownian motion and not the process in~\eqref{eq:Z}; Theorem~\ref{thm:OTM} and Proposition~\ref{prop:EL} hold for this driver, with the same proofs.)
We consider the rough stochastic volatility model
$\D S_t = \sigma(X_t)S_t \D W_t$, and assume 
that the instantaneous variance can be expanded around $y=0$ as
\begin{equation}\label{sig2exp}
(\sigma(y))^2 = v_0 + v_1 y + v_2 y^2 + \Oo\left(y^3\right).
\end{equation}
Under these assumptions we can compute a more explicit analytical expansion for the rate function for variance options. 

\begin{remark}
This class includes for example the mixed rough Bergomi model with
\begin{equation}\label{eq:mixrB}
\sigma^2(y) = v_0 \sum_{i=1}^N w_i  \E^{\eta_i y},
\end{equation}
where the weights $(w_{i})$
are positive and sum up to~$1$.
Expanding in powers of~$y$ gives
\begin{equation}\label{coeffrB}
v_1 = v_0 \sum_{i=1}^N w_i \eta_i,
\qquad v_2 = \half v_0 \sum_{i=1}^N w_i \eta_i^2,\qquad \ldots.
\end{equation}
\end{remark}

Denote the rate function for OTM variance options appearing in the short maturity 
asymptotics of Theorem~\ref{thm:OTM} with the Riemann-Liouville process as
\begin{equation}
J_H(K) := \inf \left\{\half \|h\|^2_{\Hh_1} :
    \int_{0}^{1} \sigma\left((\Kk_H \dot{h})(t)\right)^2 \D t = K \right\},
\end{equation}
with $y_0=0$ since the driver~$X$ starts at zero.
The following proposition, proved in Appendix~\ref{sec:proof_prop:Jexpansion}, provides an expansion of~$J_H(K)$ 
in powers of the log-moneyness of the variance option $y = \log(\frac{K}{v_0})$:

\begin{proposition}\label{prop:Jexpansion}
Assume that~$\sigma^2(\cdot)$ admits the expansion~\eqref{sig2exp} with $v_{1}\neq 0$.
Then 
\begin{equation}\label{JH0}
J_H(K) =  j_2 y^2 + j_3 y^3 + \Oo\left(y^4\right),
\end{equation}
where the coefficients read
\begin{equation}\label{eq:coeffs_j2_j3}
j_2:= \frac{(H+\half)^2 (H+1)}{2H} \frac{v_0^2}{v_1^2} 
\qquad\text{and}\qquad
j_3:= 
j_2 \left( 1 - 8(H+1)^2 \frac{v_0 v_2}{v_1^2} R(H)\right),
\end{equation}
with 
\begin{equation}\label{Rdef}
R(H) := \int_{0}^{1} s^{2H+1} {}_2F_1\left(1,-H-\half,H+\frac{3}{2},s\right)^2\D s.
\end{equation}
\end{proposition}

The observation in the following lemma will be useful for computational observations and will be used later to constrain the sign of the ATM skew of options on realised variance in the rough Bergomi models.

\begin{lemma}\label{lem:boundR}
The inequality 
$R(H)
> \frac{1}{4(H+1)^2}$
holds for all $H \in [0,\half]$.
\end{lemma}
\begin{proof}
By Cauchy-Schwarz inequality,
for any $H\in[0,\half]$,
\begin{equation}\label{eq:RH_lowerbound}
R(H)
:= \int_{0}^{1} s^{2\left(H+\half\right)} (\mathrm{F}_{1}(s))^2\D s
\geq \left(\int_{0}^{1} s^{H+\half} \mathrm{F}_{1}(s)\D s\right)^2,
\end{equation}
where
$$
\mathrm{F}_{1}(s) := {}_2F_1\left(1,-H-\half,H+\frac{3}{2},s\right)
 = \left(H+\half\right)\int_{0}^{1}(1-t)^{H-\half}(1-st)^{H+\half}\D t.
$$
Let $u=st$; By Fubini's theorem, the lower bound~\eqref{eq:RH_lowerbound} (without the square) reads
\begin{align*}
\int_{0}^{1} s^{H+\half} \mathrm{F}_{1}(s)\D s 
 & = \left(H+\half\right)\int_{0}^{1} s^{H+\half} \left(\int_{0}^{1}(1-t)^{H-\half}(1-st)^{H+\half}\D t\right)\D s\\
 & = \left(H+\half\right)\int_{0}^{1} \left(\int_{0}^{s}(s-u)^{H-\half}(1-u)^{H+\half}\D u\right)\D s\\
 & = \left(H+\half\right)\int_{0}^{1} (1-u)^{H+\half}\int_{u}^{1}(s-u)^{H-\half}\D s\D u
\\
 & = \int_{0}^{1} (1-u)^{2\left(H+\half\right)}\D u
 = \frac{1}{2(H+1)},
 \end{align*}
hence 
$R(H)
\geq \frac{1}{4(H+1)^2}$,
with strict inequality
since $s^{H+\half}\mathrm{F}_{1}(s)|_{s=0}=0$ and $\mathrm{F}_{1}(1)>0$.
\end{proof}
\begin{remark}
To our knowledge, the equality
$\int_0^1 s^\alpha {}_2F_1(1,-\alpha,\alpha+1,s) \D s
= \frac{1}{1+2\alpha}$
holds in fact for any
$\alpha > -\half$
(the proof of the lemma is unchanged) 
but is not available from standard textbooks such as~\cite{abramowitz1970handbook}.
\end{remark}

In the rough Bergomi model, we can compute $j_2$ and $j_3$ explicitly 
as an immediate corollary of Proposition~\ref{prop:Jexpansion}.

\begin{corollary}\label{corr:rBergomi}
In the rough Bergomi model,
$\sigma^2(y) = v_0 \E^{\eta y}$,
so that $v_1=v_0\eta$ and $v_2=\half  v_0 \eta^2$ from~\eqref{coeffrB} and
the coefficients~\eqref{eq:coeffs_j2_j3} read
\begin{equation}
j_2 = \frac{(H+\half)^2 (H+1)}{2\eta^2 H} 
\qquad\text{and}\qquad
j_3 = j_2 \left( 1 - 4(H+1)^2 R(H)\right).
\end{equation}
and $j_{3}\leq 0$ from Lemma~\ref{lem:boundR}. 
\end{corollary}


%

\begin{remark}
In the limiting Markovian case $H=\half$, Corollary~\ref{corr:rBergomi} reads
\begin{equation}
J_{H=\half}(K) = \frac{3}{2\eta^2} y^2 - \frac{3}{10\eta^2} y^3 + \Oo\left(y^4\right),
\end{equation}
which agrees with the first two terms in the expansion of the rate function for an Asian option in Black-Scholes~\cite[Equation~(35)]{pirjol2016short}.
In this case 
$R(\half) = \frac{2}{15}$.
\end{remark}

\section{Implied volatility of variance options}\label{sec:impliedvol}

The short-maturity asymptotics for OTM options on realised variance in Theorem~\ref{thm:OTM} can be translated into asymptotics for the implied volatility for options on realised variance.
As in Section~\ref{sec:expansion}, we only consider rough volatility models driven by a Riemann-Liouville process with Hurst index~$H$.
We denote $\sigma_{\BS}(K,T)$ the implied volatility of an option on realised variance with strike~$K$ and maturity~$T$,
defined
as the unique non-negative solution to 
\begin{equation}
\Call(K,T) = \Call_{\BS}\left(K,T;\frac{1}{T}
\EE\left[\int_0^T \sigma^2(Y_t)\D t\right],\sigma_{\BS}(K,T)\right),
\end{equation}
and with the Black-Scholes Call option price~\cite{black1973pricing} 
$$
\Call_{\BS}(K,T;F,\sigma) := F \Nn(d_{+}) - K \Nn(d_{-}),
\quad\text{with}\quad
d_{\pm} := \frac{\log\left(\frac{F}{K}\right)}{\sigma\sqrt{T}} \pm \half \sigma\sqrt{T}.
$$
Since $\log \Call_{\BS}(K,T) \sim -\frac{1}{2\sigma^2 T} \log\left(\frac{K}{F(0)}\right)^2$ as $T\downarrow 0$, standard arguments~\cite{Forde2011,pirjol2016short}
yield
\begin{equation} 
\lim_{T\downarrow 0} \frac{\sigma_{\BS}^2(K,T)}{T^{2(H-\half)}} = \frac{y^2}{2J_H(K)}.
\end{equation}
Expanding the implied volatility in log-moneyness $y=\log\left(\frac{K}{K_{\mathrm{ATM}}}\right)
 = \log\left(\frac{K}{v_{0}}\right)$, we obtain
\begin{equation}
\sigma_{\BS}(K,T) = \sigma_{\mathrm{ATM}}(T) + s(T) y + \Oo\left(y^2\right).
\end{equation}
The next result,
which follows directly from Proposition~\ref{prop:Jexpansion},
gives the short-maturity asymptotics of these coefficients.
 
\begin{proposition}\label{prop:sigimp}
The short-maturity asymptotics for the ATM implied volatility and skew of an option on realised variance in a rough stochastic volatility model driven by a Riemann-Liouville process, and volatility $\sigma^2(y) = v_0 + v_1 y + v_2 y^2 + \Oo\left(y^3\right)$ are given by
\begin{equation}\label{ATMlimit}
\lim_{T\downarrow 0} \frac{\sigma_{\mathrm{ATM}}(T)}{T^{H-\half}} = \frac{1}{\sqrt{2j_2}} = \frac{1}{H+\half} \sqrt{\frac{H}{H+1}} \frac{|v_1|}{v_0},
\end{equation}
and
\begin{equation}\label{skewlimit}
\lim_{T\downarrow 0} \frac{s(T)}{T^{H-\half}} = -\frac{j_3}{(2j_2)^{\frac{3}{2}}} = -\frac{1}{2(H+\half)} \sqrt{\frac{H}{H+1}}
\frac{|v_1|}{v_0} \left( 1- 4(H+1)^2 \frac{2v_0 v_2}{v_1^2} R(H) \right),
\end{equation}
with~$R(H)$ defined in~\eqref{Rdef}. 
\end{proposition}

\begin{remark}
In the case of the Molchan-Golosov kernel, the same computation can be extracted from Theorem~\ref{thm:ATM} and yields
\begin{equation}
\lim_{T\downarrow 0} \frac{\sigma_{\mathrm{ATM}}(T)}{T^{H-\half}} = 
\frac{1}{\sqrt{2(H+1)}} \frac{|v_1|}{v_0},
\end{equation}
with the same $T$-scaling as~\eqref{ATMlimit} but different normalisation factor, as explained in Remark~\ref{rem:ATMcompare}. 
At $H=\half$ the two results agree.
\end{remark}

In the log-normal SABR model
with $H=\half$ and 
$\sigma^2(y) = \sigma_{0}^{2} \E^{\eta y}$, with $\eta>0$,
\eqref{ATMlimit} recovers the short-maturity asymptotics of the ATM Asian implied volatility in Black-Scholes~\cite[Equation~(55)]{pirjol2016short}:
\begin{equation}
\lim_{T\downarrow 0} \sigma_{\BS}(K_{\mathrm{ATM}},T) = \frac{\eta}{\sqrt3}.
\end{equation}
For $H=\half$, the short-maturity asymptotics of the ATM skew~\eqref{skewlimit} also reproduces the known result for the Asian skew in Black-Scholes $s(0) = \frac{\eta}{10\sqrt3}$ as in~\cite[Equation~(35)]{Alos2022} and~\cite[Remark~20 and Equation~(55)]{pirjol2016short}, 
given more precisely as 
\begin{equation}
\lim_{T\downarrow 0}\sigma_{\BS}^{\textrm{Asian},\BS}\left(\frac{K}{S_0},T\right)
= \frac{\sigma}{\sqrt3}\left(
1 + \frac{1}{10} x - \frac{23}{2100} x^2 + \Oo\left(x^3\right) \right),\quad x = \log\left(\frac{K}{S_0}\right).
\end{equation}
The first two terms in this expansion were obtained
using Malliavin calculus techniques in~\cite{Alos2022} for a general stochastic volatility model. 
The application to the Black-Scholes model is presented in~\cite[Section~5.1]{Alos2022}.


\subsection{ATM skew sign} 
For $H\in (0,\half]$ the ATM skew~\eqref{skewlimit} of the variance option implied volatility in rough Bergomi is always positive
by Corollary~\ref{corr:rBergomi} and  Lemma~\ref{lem:boundR}.
We show that this holds more generally for the mixed rough Bergomi model.

\begin{corollary}
Consider the mixed rough Bergomi model~\eqref{eq:mixrB}
driven by a Riemann-Liouville process with $H\in (0,\half]$ and $\sum_{i=1}^{N}w_{i}\eta_{i}\neq 0$. 
The ATM skew of the variance option is positive for sufficiently small maturity.
\end{corollary}

\begin{proof}
In the mixed rough Bergomi model, $v_1 = v_{0}\sum_{i=1}^N w_i \eta_i$ and $v_2 = \half v_{0}\sum_{i=1}^N w_i \eta_i^2$,
so that Jensen's inequality gives
$$
\frac{2v_0 v_2}{v_1^2} = \frac{\sum_{i=1}^N w_i \eta_i^2}
{\left(\sum_{i=1}^N w_i \eta_i\right)^2} \geq 1,
$$
and therefore
$1- 4(H+1)^2 \frac{2v_0 v_2}{v_1^2} R(H) \leq 1- 4(H+1)^2 R(H) < 0$
by Lemma~\ref{lem:boundR}.
Substituting into~\eqref{skewlimit} the conclusion follows.
\end{proof}


\subsection{Comparison with the literature}\label{sec:literature}
Al\`os, Garc\'ia-Lorite and Muguruza~\cite{Alos2018} study the short-maturity asymptotics of variance options in rough volatility models, including the rough Bergomi model. 
In the mixed rough Bergomi model with $\sigma^2(y) = v_0(\gamma \E^{\nu y} + (1-\gamma) \E^{\eta y})$,
the first two coefficients can be read off from~\eqref{coeffrB}:
$v_1=v_0(\gamma\nu + (1-\gamma)\eta) $ and $v_2 = \half v_0(\gamma\nu^2 + (1-\gamma)\eta^2)$.
The ATM volatility and skew are given in~\cite[Example~14,  Example~21]{Alos2018}:
\begin{align*}
&\lim_{T\downarrow 0} \frac{\sigma_{\mathrm{ATM}}(T)}{T^{H-\half}} = \frac{1}{H+\half} \sqrt{\frac{H}{H+1}} (\gamma \nu + (1-\gamma)\eta),\\
&\lim_{T\downarrow 0} \frac{s(T)}{T^{H-\half}} 
= - \frac{(\gamma\nu + (1-\gamma) \eta)\sqrt{H}}{(2H+1)\sqrt{H+1}}
+
\frac{\gamma \nu^2 + (1-\gamma) \eta^2}{\gamma\nu + (1-\gamma)\eta} 
\sqrt{2H}(2H+2)^{\frac{3}{2}} \left(H+\half\right) \Ii(H),\nonumber
\end{align*}
as in~\eqref{ATMlimit}
with $\Ii(H)$ given in~\cite[Equation~(58)]{Alos2018}:
\begin{equation}\label{Icaldef}
\Ii(H) = \lim_{T\downarrow 0} \frac{1}{T^{4H+3}}
\int_0^T \frac{(T-s)^{H+\half}}{H+\half} 
\int_s^T (T-u)^{2H+1} (u-s)^{H-\half} \mathrm{F}_{2}\left(\frac{T-u}{s-u}\right)\D u\D s,
\end{equation}
and $\mathrm{F}_{2}(x) := {}_2F_1\left(\half-H,H+\half,H+\frac{3}{2},x\right)$.
Our prediction from~\eqref{skewlimit} is very similar:
\begin{equation}
\lim_{T\downarrow 0} \frac{s(T)}{T^{H-\half}}
=  
-\frac{\sqrt{H}}{(2H+1)\sqrt{H+1}}
\left\{
(\gamma\nu + (1-\gamma) \eta)
-
4(H+1)^{2} R(H)\frac{\gamma \nu^2 + (1-\gamma) \eta^2}{\gamma\nu + (1-\gamma)\eta} 
\right\}.
\end{equation}
Comparing the two expressions, the first (negative) term is the same. 
The positive terms would be also identical if 
\begin{equation}\label{eq:identity_I_R}
\Ii(H)
= \frac{R(H)}{2(H+\half)^2}.
\end{equation}
The following lemma is reassuring:
\begin{lemma}
The identity~\eqref{eq:identity_I_R} holds for all $H\in [0,\half]$.
\end{lemma}
\begin{proof}
The first observation is that the double integral in~$\Ii(H)$ is homogeneous in~$T$
and the limit~$\Ii(H)$ is in fact superfluous.
Since
${}_2F_1(a,b,c,z) 
= (1-z)^{-a}{}_2F_1\left(a,c-b,c,\frac{z}{z-1}\right)$
from~\cite[15.3.4]{abramowitz1970handbook},
setting $z=\frac{b}{b-a}$,
then $\frac{z}{z-1} = \frac{b}{a}$ and
\begin{equation}\label{eq:ident_2F1}
\mathrm{F}_2\left(\frac{b}{b-a}\right) 
= \left(1-\frac{b}{a}\right)^{\half-H}{}_2F_1\left(\half-H,1,H+\frac{3}{2},\frac{b}{a}\right).
\end{equation}
Therefore, with the change of variables $s=T(1-a)$
and $u=T(1-b)$ with $0<b<a<1$, then
$u-s = (a-b)T$, $T-s = aT$ and $T-u = bT$, so that  
\begin{align*}
\mathcal{I}(H)& = \frac{1}{T^{4H+3}}
\int_{0}^{T} \frac{(T-s)^{H+\half}}{H+\half} 
\int_{s}^{T} (T-u)^{2H+1} (u-s)^{H-\half} \mathrm{F}_{2}\left(\frac{T-u}{s-u}\right)\D u\D s\\
& = \frac{T^2}{T^{4H+3}}
\int_{0}^{1} \frac{(aT)^{H+\half}}{H+\half} 
\int_{0}^{1-\frac{s}{T}} (bT)^{2H+1} ((a-b)T)^{H-\half} \mathrm{F}_{2}\left(\frac{b}{b-a}\right)\D b\D a\\
& = \frac{1}{H+\half} 
\int_{0}^{1} a^{H+\half}
\int_{0}^{a} b^{2H+1} (a-b)^{H-\half} \mathrm{F}_{2}\left(\frac{b}{b-a}\right)\D b\D a\\
& = \frac{1}{H+\half} 
\int_{0}^{1} a^{2H}
\int_{0}^{a} b^{2H+1}{}_2F_1\left(\half-H,1,H+\frac{3}{2},\frac{b}{a}\right)\D b\D a,
\end{align*}
using~\eqref{eq:ident_2F1} and noting that the inner upper bound is equal to
$1-\frac{s}{T} = a$.
Finally, letting $b=at$ in the inner integral, 
we obtain (without the un-necessary limit)
\begin{align}\label{eq:I_H_simplified}
\Ii(H)
 & = \frac{1}{H+\half} 
\int_{0}^{1} a^{4H+2}
\int_{0}^{1} t^{2H+1}{}_2F_1\left(\half-H,1,H+\frac{3}{2},t\right)\D t\D a\nonumber\\
 & = \frac{1}{H+\half} 
\left(\int_{0}^{1} a^{4H+2}\D a\right)
\left(\int_{0}^{1} t^{2H+1}{}_2F_1\left(\half-H,1,H+\frac{3}{2},t\right)\D t\right)\nonumber\\
 & = \frac{1}{(H+\half)(4H+3)}
 \int_{0}^{1} t^{2H+1}{}_2F_1\left(\half-H,1,H+\frac{3}{2},t\right)\D t.
\end{align}
Comparing with~\eqref{Rdef},
there does not seem to be an immediate link between the two forms of~${}_2F_1$ in~\cite{abramowitz1970handbook}.
However, write
$$
\mathfrak{i}(t):= {}_2F_1\left(\half-H,1,H+\frac{3}{2},t\right)
\quad\text{and}\quad
\mathfrak{r}(t):= \left({}_2F_1\left(1,-H-\half,H+\frac{3}{2},t\right)\right)^{2},
$$
corresponding to the integrands in~$\Ii(H)$ and~$R(H)$ respectively.
Then, using~\cite[15.2.1]{abramowitz1970handbook},
$$
\frac{\D}{\D t}
\left\{t^{2H+1}(1-t)(\mathfrak{i}(t) - \mathfrak{r}(t))\right\}
 = (4H+3) t^{2H+1}\mathfrak{r}(t)
 -(2H+1) t^{2H+1}\mathfrak{i}(t);
$$
integrating between~$0$ and~$1$,
since the boundary terms vanish, then
$$
\int_{0}^{1}t^{2H+1}\mathfrak{i}(t)\D t = \frac{4H+3}{2H+1}\int_{0}^{1}t^{2H+1}\mathfrak{r}(t)\D t,
$$
and therefore, from~\eqref{eq:I_H_simplified} and~\eqref{Rdef}, we can write
$$
\Ii(H)
 = \frac{1}{(H+\half)(4H+3)}
 \int_{0}^{1} t^{2H+1}\mathfrak{i}(t)\D t
= \frac{1}{2(H+\half)^2}
\int_{0}^{1}t^{2H+1}\mathfrak{r}(t)\D t
= \frac{R(H)}{2(H+\half)^2}.
$$
This completes the proof.
\end{proof}

The results in~\cite{Alos2018} have been tested numerically for arbitrary strikes (without expanding in log-moneyness) in~\cite{lacombe2021asymptotics} by a direct numerical solution of the variational problem for the rate function $J_H(K)$.
We present an independent test in the next section.

\subsection{Numerical tests in the rough Bergomi model}

We compare here the analytical results for the rate function and implied volatility with a numerical evaluation of the variational problem in Theorem~\ref{thm:OTM} for the rough Bergomi model with $\sigma^2(y) = v_0 \E^{\eta y}$.
This variational problem can be approximated with an optimisation problem by discretising the function $h:[0,1]\to \RR$ 
on the grid $\{0,\frac{1}{n},\frac{2}{n},\cdots,1\}$.
Denote $g_j := h(\frac{j}{n}) - h(\frac{j-1}{n})$ with $j=1,2,\cdots,n$ the increments of~$h(\cdot)$ 
on each interval of step $\frac{1}{n}$ and recall that $h(0)=0$.
The variational problem is reduced to an optimisation over the~$n$ variables $\{g_j\}_{j=1}^n$.
This is different from the approximation in~\cite{lacombe2021asymptotics} 
who approximate the optimal path with an expansion in a polynomial basis, 
with polynomials up to 8-th degree. 
Our results therefore give an independent test of the analytical computation of the rate function by expansion around the ATM point.

Denote by~$J_{H,n}(y)$ the approximation of the rate function following from this optimisation problem with~$n$ points, given by
\begin{equation}\label{JHn}
J_{H,n}(y) = \inf_{g_i} \frac{n}{2} \sum_{i=1}^n (g_i)^2,
\end{equation}
where the infimum is taken subject to a non-linear constraint 
$\int_0^1 \sigma^2( (\mathcal{K}_H \dot{h})(t)) \D t = K$. 
In the rough Bergomi model $\sigma^2(y)=v_0 \E^{\eta y}$ the constraint is
expressed as the double sum
\begin{equation}\label{constraint-discrete}
\frac{1}{n}\sum_{i=1}^{n}\exp\left(\frac{\eta\sqrt{2H}}{(2H+1)(2n)^{H-\half}}
\left\{\sum_{j=1}^{i-1}g_{j}\left[(2i-2j+1)^{H+\half}-(2i-2j-1)^{H+\half}\right]+g_{i}\right\}\right)=\E^{y},
\end{equation}

This was obtained by approximating $h$ with a piecewise linear function with derivative
$\dot{h} = n g_{j}$ on $(\frac{j-1}{n},\frac{j}{n}]$ and evaluating the integral in the exponent of the constraint exactly.
A similar approach is used to obtain the objective function in~\eqref{JHn}.
For $t\in(\frac{i-1}{n},\frac{i}{n}]$,
$$
(\Kk_{H}\dot{h})(t)
=\frac{\sqrt{2H}n}{H+\half}\left\{\sum_{j=1}^{i-1}g_{j}\left[\left(t-\frac{j-1}{n}\right)^{H+\half}-\left(t-\frac{j}{n}\right)^{H+\half}\right]
+g_{i}\left(t-\frac{i-1}{n}\right)^{H+\half}\right\}.
$$

\begin{figure}[h!]
\centering
\includegraphics[width=1.75in]{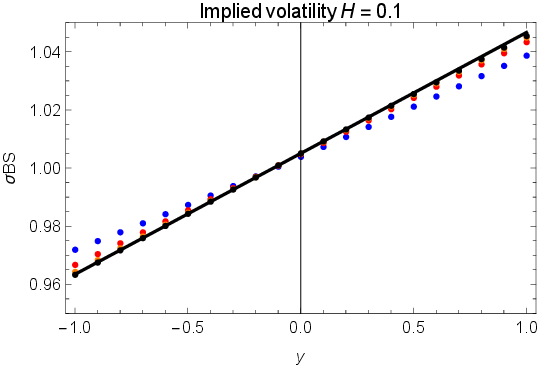}
\includegraphics[width=1.75in]{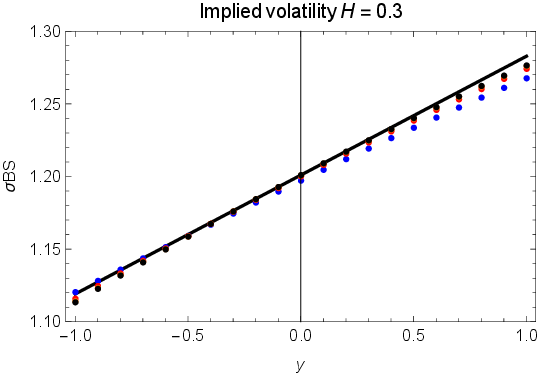}
\includegraphics[width=1.75in]{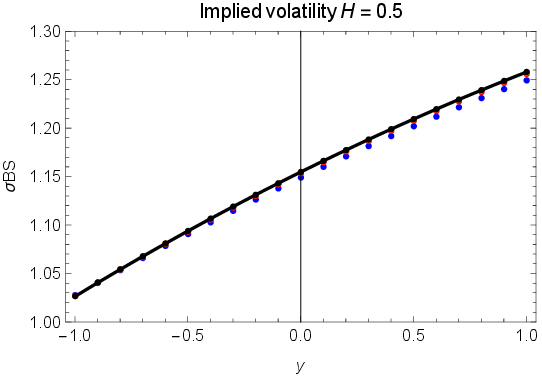}
\caption{
Asymptotic implied volatility in rough Bergomi 
(rescaled by a factor of $T^{-H+\half }$).
The black lines for $H\in\{0.1, 0.3\}$ 
show the analytical result from Proposition~\ref{prop:sigimp}.
The dots show the numerical evaluation by discretisation as described in text
with $n=5$ (blue), 10 (red), 20 (orange) and 40 (black) points. 
For the Markovian case $H=\half$ the black curve shows the full asymptotic result 
known from the Asian options in the Black Scholes model. }
\label{Fig:numimpvol}
\end{figure}

Finally, the outer integral in the constraint is approximated using the mid-point rule:
\begin{equation}
\frac{1}{n} \sum_{i=1}^n \sigma^2\left((\Kk_{H}\dot{h})\left(\frac{2i-1}{2n}\right)\right) = K.
\end{equation}
Substituting here $\sigma^2(y)=v_0 \E^{\eta y}$ gives the constraint
\eqref{constraint-discrete}.
For the numerical tests we take $\eta=2$ and $H \in \{ 0.1, 0.3, 0.5\}$. 
Figure~\ref{Fig:numimpvol} shows the rescaled asymptotic implied volatility 
$\lim_{T\downarrow 0} \sigma_{\mathrm{BS}}(K,T) T^{-H+\half }$. The solid black curve shows the linear approximation from Proposition~\ref{prop:sigimp} (for the Markovian case $H=\half$ we used the exact result). 
The coloured dots show the numerical evaluation from discretisation with~$n=5$ (blue), 10 (red), 20 (orange) and 40 points (black). 
As~$n$ increases, the numerical evaluation approaches the theoretical value from the~$y$ expansion keeping the first two terms in the rate function. 

The numerical test confirms that keeping only the linear term in the asymptotic implied volatility is a good approximation for the exact result. 
This is also confirmed analytically for the Markovian case $H=\half$ where the convexity is known $\frac{23}{2100} = 0.01$, 
and is a factor of~$10$ smaller than the coefficient of the linear term~$0.1$.
The numerical tests suggest that a similar suppression of the convexity term holds also for $H\neq \half $, which explains the approximate linearity of the smile for all values of~$H$. 


\section{Conclusion}\label{sec:conclusion}

We established small-time large deviations, central limit theorem 
and moderate deviations for the normalised time integral of the exponential of a Gaussian Volterra process. 
As an application we present short maturity asymptotics for options on realised variance in rough volatility models, 
for both regimes of out-of-the-money and at-the-money variance options.  
We consider in detail two types of rough volatility models, 
driven by a fractional Brownian motion 
and by a Gaussian Volterra process with Riemann-Liouville kernel. 
For the latter, we solve the variational problem for the rate function in an expansion around the at-the-money point,
which is then used to obtain an expansion of the implied volatility in log-moneyness. 
The first two terms in this expansion are obtained in closed form, confirming and simplifying~\cite{Alos2018}. 
We also prove the positivity of the ATM skew of the implied volatility of variance options for all $H\in (0,\half]$ in the rough Bergomi model.
These results are numerically validated through a series of tests.


\appendix

\section{Comparison of the two kernels}\label{sec:comment}

For $H\in(0,\half)$, consider the kernels
$K_{H}$ from Section~\ref{sec:options} 
with 
$K(u) = u^{H-\half}$ for $u>0$, 
following Assumption~\ref{ass:main} with $\ell\equiv 1$.
Therefore
$W^{H}_t=\int_0^t K_{H}(s,t) \D W_s$ defines
a standard fBm, 
with $\VV[W^{H}_1] = 1$ and covariance
$R_H(t,s)=\half\left(t^{2H}+s^{2H}-|t-s|^{2H}\right)$,
while $X^H_t=\frac{1}{\Gamma(H+\half)}\int_0^t K(t-s)\D W_s$ defines the Riemann-Liouville version.

\begin{proposition}\label{prop:fact}
For $0<s<t$,
$K_{H}(s,t) = K(t-s)\chi(\frac{t}{s})$ with
\begin{equation}\label{eq:phi}
\chi(x) := c_H\left[x^{H-\half}
+\left(\half-H\right)(x-1)^{\half-H}
\int_1^x v^{H-\frac{3}{2}}(v-1)^{H-\half} \D v\right],
\qquad x\geq 1.
\end{equation}
\end{proposition}

\begin{proof}
Clearly,
$K_{H}(as,at)=a^{H-\half}K_{H}(s,t)$ for every $a>0$.
Since $K$ has the same homogeneity, 
write $K_{H}(s,t) = K(t-s)L(s,t)$, 
so~$L$ is homogeneous of degree zero and thus a
function of $x=\frac{t}{s}$ alone. 
Taking $s=1$, $t=x$ and substituting $u=v$ gives~\eqref{eq:phi}.
\end{proof}

\begin{lemma}\label{lem:ode}
The function $\chi\in \Cc^\infty((1,\infty))\cap \Cc([1,\infty))$ and solves
\begin{equation}\label{eq:ode}
(x-1) \chi'(x)+\left(H-\half\right)\chi(x)
=c_H\left(H-\half\right)x^{H-\half},
\end{equation}
for all $x>1$, with boundary condition
$\chi(1)=c_H$.
\end{lemma}

\begin{proof}
Differentiate the expression for~$K_{H}$ with respect to~$t$ and obtain~\cite[\S5.1.3]{nualart2006malliavin}
\begin{equation}\label{eq:dK}
\frac{\partial K_{H}}{\partial t}(s,t)
=c_H\left(H-\half\right)\left(\frac{s}{t}\right)^{\half-H}
(t-s)^{H-\frac{3}{2}}.
\end{equation}
Now, differentiating
$K(t-s)\chi(\frac{t}{s})$ in $t$ gives
$$
\left(H-\half\right)(t-s)^{H-\frac{3}{2}}\chi\left(\frac{t}{s}\right)
+\frac{(t-s)^{H-\half}}{s}\chi'\left(\frac{t}{s}\right). 
$$
Equating with
\eqref{eq:dK} and using $t-s=s(x-1)$ yields~\eqref{eq:ode}. 
The boundary condition
$\chi(1)=c_H$ and the continuity at~$x=1$ are read off~\eqref{eq:phi}, since the integral there is
$\Oo((x-1)^{H+\half})$ as $x\downarrow1$.
\end{proof}

\begin{theorem}\label{thm:main}
For~$\chi$ in~\eqref{eq:phi}, the following hold:
\begin{enumerate}[(i)]
\item\label{it:pos} $\chi\in\Cc([1,\infty)\to[0,\infty))\cap\Cc^{\infty}((1,\infty)\to[0,\infty))$ and strictly increasing.
\item\label{it:diag} As $x\downarrow 1$,
\begin{equation}\label{eq:diag}
\chi(x)=c_H\left[1+\frac{(\half-H)^2}{H+\half} (x-1)
 + \Oo\left((x-1)^2\right)\right].
\end{equation}
\item\label{it:rv} As $x\uparrow\infty$,
$\chi(x) \sim \frac{H}{c_H} x^{\half-H}$, so~$\chi$ is regularly varying at $+\infty$ with index
$\half-H$.
\item\label{it:bounds} There exist
$0<\underline{\kappa}_H\le \overline{\kappa}_H<\infty$, depending only on $H$, such that
\begin{equation}\label{eq:twosided}
\underline{\kappa}_H x^{\half-H}\le\chi(x)\le \overline{\kappa}_H x^{\half-H}
\qquad \text{for all }x\geq 1.
\end{equation}
\end{enumerate}
\end{theorem}

\begin{proof}
\eqref{it:pos} Positivity, continuity and smoothness are immediate. 
Regarding monotonicity, rewrite~\eqref{eq:ode} as
$$
(x-1)\chi'(x)=\left(\half-H\right)\left(\chi(x)-c_Hx^{H-\half}\right);
$$
by~\eqref{eq:phi} the bracket equals
$$
c_H\left(\half-H\right)(x-1)^{\half-H}\int_1^x v^{H-\frac{3}{2}}(v-1)^{H-\half}\D v>0, 
$$
for
$x>1$, so $\chi'>0$ on $(1,\infty)$.

\eqref{it:diag} Expand the two summands of~\eqref{eq:phi} at $x=1$:
$$
x^{H-\half}=1+\left(H-\half\right)(x-1)+\Oo\left((x-1)^2\right), 
$$
and
$$
\int_1^x v^{H-\frac{3}{2}}(v-1)^{H-\half}\D v
=\frac{(x-1)^{H+\half}}{H+\half} \left(1+\Oo(x-1)\right), 
$$
so the second
summand is 
$\frac{\half-H}{H+\half}(x-1) (1+\Oo(x-1))$.
Adding the linear
coefficients,
$$
\left(H-\half\right)+\frac{\half-H}{H+\half}
=\left(\half-H\right)\left(\frac{1}{H+\half}-1\right)
=\frac{(\half-H)^2}{H+\half}, 
$$
which gives~\eqref{eq:diag}. 

\eqref{it:rv} As $x\uparrow\infty$ the first summand in~\eqref{eq:phi} tends
to $0$, $(x-1)^{\half-H}\sim x^{\half-H}$, and by monotone convergence theorem
\begin{equation*}
\int_1^\infty v^{H-\frac{3}{2}}(v-1)^{H-\half} \D v
\overset{v=\frac{1}{u}}{=}\int_{0}^{1} u^{-2H}(1-u)^{H-\half} \D u
=\mathrm{B}\left(1-2H, H+\half\right)<\infty.
\end{equation*}
Hence
$$
\chi(x)\sim c_H\left(\half-H\right)\mathrm{B}\left(1-2H,H+\half\right) x^{\half-H}, 
$$
and 
$$
c_H^2\left(\half-H\right)\mathrm{B}\left(1-2H,H+\half\right)=H.
$$
A measurable function asymptotic to a positive multiple of~$x^{\rho}$ is regularly varying of index~$\rho$ in the sense of
Karamata~\cite{bingham1987regular}.

\eqref{it:bounds} The map $x\mapsto\chi(x)x^{H-\half}$ is
continuous and positive on $[1,\infty)$,
converges to $\frac{H}{c_H}>0$ at
infinity by~\eqref{it:rv}, 
hence attains a positive minimum and a
finite maximum.
\end{proof}

\begin{remark}\label{rem:hyper}
It follows from Decreusefond-\"Ust\"unel~\cite{decreusefond1999stochastic} that
\begin{equation}\label{eq:DU}
K^{\mathrm{DU}}_H(s,t)=\frac{(t-s)^{H-\half}}{\Gamma(H+\half)} 
{}_2\mathrm{F}_1\left(H-\half, \half-H; H+\half; 
1-\frac{t}{s}\right),\qquad 0<s<t,
\end{equation}
so that the factor $\chi(x)$ is, up
to a multiplicative constant, ${}_2\mathrm{F}_1$ evaluated at
$1-\frac{t}{s}\in(-\infty,0]$. 
Concretely, with the right normalisation,
for all $x\geq 1$,
$$
\chi(x)=c_H {}_2\mathrm{F}_1\left(H-\half, \half-H; H+\half; 1-x\right).
$$
\end{remark}

\section{Additional Technical Proofs}\label{sec:proofs_all}

\subsection{Large deviations proofs}
\label{sec:LDP:proofs}

\subsubsection{Rescaling and uniform estimates}

For $\eps>0$, introduce
\begin{equation}\label{eq:phi_notations}
\varphi_{\eps}(t) := \varphi(\eps t)-\varphi(0),
\qquad
\overline{\varphi}_{\eps} := \sup_{t\in[0,1]}|\varphi_{\eps}(t)|,
\qquad
\Phi_{\eps} := \int_{0}^{1} \E^{\varphi_{\eps}(t)}\D t.
\end{equation}

We now define the rescaling
\begin{equation}\label{eq:hatXT}
X^{\eps}_{t}
:= \frac{X^{(\ell)}_{\eps t}}{\eps^H\ell(\eps)},\qquad \text{for all }t\in[0,1],
\end{equation}
and define
$W^{\eps}_{t} := \frac{1}{\sqrt{\eps}}W_{\eps t}$ and 
$\psi_{\eps}(t) := t^{H-\half}\ell_{\eps}(t)$
with $\ell_{\eps}(t):=\frac{\ell(\eps t)}{\ell(\eps)}$,
so that
\begin{equation}\label{eq:rep_Xhat_T}
X^{\eps}_{t}
 = \int_{0}^{t}\psi_{\eps}(t-u) \D W^{\eps}_{u}
\qquad\text{with}\qquad
\VV\left[X^{\eps}_{t}\right]
=\frac{\VV[X^{(\ell)}_{\eps t}]}{\eps^{2H}\ell(\eps)^2},
\qquad
\text{for any } t\in[0,1],
\end{equation}

We shall use the following elementary inequality:

\begin{lemma}\label{lem:monotone}
For $f:(0,\infty)\to[0,\infty)$ non-increasing with $f\in L^2((0,c))$ for any $c>0$,
$$
\int_0^{t}
\left(f(u) - f(u+h)\right)^2\D u \leq \int_0^h (f(u))^2 \D u,
\qquad\text{for any }h,t>0.
$$
\end{lemma}

\begin{proof}
For $0\leq b\leq a$, $(a-b)^2\leq a^2-b^2$; 
letting $a=f(u)$ and $b=f(u+h)$,
we obtain
\begin{align*}
\int_0^t|f(u) - f(u+h)|^2\D u
&\leq \int_0^t f^2(u)\D u- \int_h^{t+h}f^2(u)\D u
\\
&=\int_0^h f^2(u)\D u- \int_{t}^{t+h}f^2(u)\D u
 \leq \int_0^h f^2(u)\D u.
\end{align*}
This completes the proof.
\end{proof}

\begin{lemma}\label{lem:holder}
Fix $\eta\in(0,H)$. 
There exist $C\geq 1$ and $\eps_{0}>0$, depending only on $\ell,\eta$, such that for all $\eps\in(0,\eps_{0}]$, $w\in(0,1]$ and $0\leq s\leq t\leq 1$,
\begin{equation}\label{eq:incr}
\ell_{\eps}(w)\leq C w^{-\eta}
\qquad\text{and}\qquad
\EE\left[\left|X^{\eps}_{t} - X^{\eps}_s\right|^2\right]
\leq \frac{C^2}{H-\eta}|t-s|^{2(H-\eta)}.
\end{equation}
The same bounds hold for $X$ 
(with $\ell\equiv1$, $C=1$, $\eta=0$). 
Thus 
$$
\vartheta^2 := \sup_{\eps\leq \eps_{0}}\sup_{t\leq 1}
\VV\left[X^{\eps}_{t}\right] < \infty
\qquad\text{and}\qquad
\overline{m} := \sup_{\eps\leq \eps_{0}}\EE\left[\|X^{\eps}\|_\infty\right] < \infty.
$$
\end{lemma}

\begin{proof}
The first bound in~\eqref{eq:incr} follows from~\cite[Theorem~1.5.6]{bingham1987regular} applied to the slowly varying function $x\mapsto \ell(\frac{1}{x})$ at infinity, together with $w\leq 1$. 
For the second bound in~\eqref{eq:incr}, Lemma~\ref{lem:monotone} 
(with Assumption~\ref{ass:main}, as~$\psi_{\eps}$ is non-increasing for $\eps\leq \eps_0$) gives, with $h:=t-s$,
$$
\EE\left[\left|X^{\eps}_{t}-X^{\eps}_s\right|^2\right]
 = \int_0^{h}(\psi_{\eps}(t))^2\D t
 + \int_0^s
 \left(\psi_{\eps}(t)-\psi_{\eps}(t+h)\right)^2\D t
 \leq 2\int_0^{h}\psi_{\eps}^2(t)\D t,
$$
and the bound follows by integration.
The uniform bound on $\EE[\|X^{\eps}\|_\infty]$ follows from Dudley's entropy bound~\cite[Chapter~11]{ledoux1991probability}: the metric $d_{\eps}(s,t) = \EE[(X^{\eps}_{t}-X^{\eps}_s)^2]^{\half}$ is dominated by $c|t-s|^{H-\eta}$ uniformly in~$\eps$, so the covering numbers satisfy $\Nn(d_{\eps},\delta)\leq 1+c'\delta^{-\frac{1}{H-\eta}}$ and $\int_0^{\infty}\sqrt{\log(\Nn(d_{\eps},\delta))} \D\delta$ is finite and bounded uniformly in~$\eps$.
\end{proof}

\begin{lemma}\label{lem:coupling}
Let  $\widetilde{X}^{\eps}_{t}:=\Kk\dot{W}^{\eps}$ and 
$\Df^{\eps} := X^{\eps} - \widetilde{X}^{\eps}$. Then
$$
\varrho_{\eps}^2
:=\sup_{t\in[0,1]}\EE\left[|\Df^{\eps}_{t}|^2\right]
\leq \int_{0}^{1} t^{2H-1}
(\ell_{\eps}(t)-1)^2\D t
\xrightarrow[\eps\downarrow 0]{}0,
\quad\text{and}\quad
\lim_{\eps\downarrow 0}
\EE[\|\Df^{\eps}\|_\infty] = 0.
$$
\end{lemma}

\begin{proof}
Since 
$\Df^{\eps}_{t} = \int_0^t(t-u)^{H-\half}(\ell_{\eps}(t-u)-1)\D W^{\eps}_u$,
we have
$$
\EE\left[|\Df^{\eps}_{t}|^2\right]
= \int_0^t u^{2H-1}(\ell_{\eps}(u)-1)^2\D u.
$$
By uniform convergence for slowly varying functions~\cite[Theorem~1.2.1]{bingham1987regular}, $\lim_{\eps\downarrow 0}\ell_{\eps}(u)= 1$ uniformly on $[\delta,1]$ for each $\delta>0$; 
by~\eqref{eq:incr}, since the integrand is dominated by 
$$
u^{2H-1}\left(Cu^{-\eta}+1\right)^2
\leq C u^{2(H-\half-\eta)}\in L^1([0,1]),
$$
then dominated convergence theorem yields $\lim_{\eps\downarrow 0}\varrho_{\eps}= 0$. 
Finally~$\Df^{\eps}$ is centred Gaussian with canonical metric dominated as in Lemma~\ref{lem:holder} and with diameter at most~$2\varrho_{\eps}$
and Dudley's bound~\cite[\S4]{dudley2016vn} thus gives 
$$
\EE[\|\Df^{\eps}\|_\infty]
\leq c \int_0^{2\varrho_{\eps}}\log\left(1+c'\delta^{-\frac{1}{H-\eta}}\right)^{\half} \D\delta\downarrow 0.
$$
\end{proof}

\subsubsection{Large deviations for $(\speed_{\eps}X^{\eps})_{\eps>0}$}

Before proving the main LDP for $(Z_{\eps}/\eps)_{\eps>0}$, we prove one for the family
$(\speed_{\eps}X^{\eps})_{\eps>0}$:

\begin{theorem}\label{thm:pathLDP}
Under Assumption~\ref{ass:main}, as $\eps\downarrow 0$,
the family $(\speed_{\eps}X^{\eps})_{\eps>0}$ satisfies an LDP on $\Cc_0([0,1])$ with speed~$\speed_{\eps}^{-2}$ and good rate function~$\If$ in~\eqref{eq:I}.
\end{theorem}

\begin{proof}
By Lemma~\ref{lem:holder} with $\ell\equiv1$ and the Kolmogorov-Centsov criterion, $X$ has a modification with paths in $\Cc_0([0,1])$, Hölder of any order strictly smaller than~$H$; its law~$\mu$ is a centred Gaussian measure on the separable Banach space $\Cc_0([0,1])$. Since $\widetilde{X}^{\eps} = \Kk\dot{W}^{\eps}$ with~$\Kk$ injective and bounded from $L^2([0,1])$ into $\Cc_0([0,1])$, the Cameron-Martin space of~$\mu$ is $\Kk(L^2([0,1]))$ with $\|\Kk\dot{h}\|_{\Hh_{1}^{X}}=\|\dot{h}\|_{L^2}$, as explained in~\cite[Chapter~2]{lifshits2013gaussian} or~\cite{decreusefond1999stochastic}.
Therefore the LDP for $(\speed_{\eps}\widetilde{X}^{\eps})$ with speed $\speed_{\eps}^{-2}$ and good rate function $\half\|\cdot\|_{\Hh_{1}^{X}}^2=\If$ 
follows from Schilder's theorem for centred Gaussian measures on separable Banach spaces~\cite[\S3.4]{deuschel2001large} and~\cite[\S5.2]{dembo1998large}.
Using Lemma~\ref{lem:coupling}, Borell-TIS inequality~\cite[Theorem~2.1.1]{adler2007random} yields, for every $\delta>0$ and all~$\eps$ small enough such that $\speed_{\eps}\EE[\|\Df^{\eps}\|_\infty]
 < \frac{\delta}{2}$,
\begin{align*}
\speed_{\eps}^2\log\PP\left(\speed_{\eps}\|\Df^{\eps}\|_\infty>\delta\right)
 & \leq \speed_{\eps}^2\log\exp\left\{-\frac{\left(\frac{\delta}{\speed_{\eps}} - \EE[\|\Df^{\eps}\|_\infty]\right)^2}{2\varrho_{\eps}^2}\right\}\\
& \leq -\frac{(\delta-\speed_{\eps}\EE[\|\Df^{\eps}\|_\infty])^2}{2\varrho_{\eps}^2}
\le-\frac{\delta^2}{8\varrho_{\eps}^2},
\end{align*}
which tends to $-\infty$ since $\varrho_{\eps}\to 0$. Hence $\left(\speed_{\eps}X^{\eps}\right)$ and $\left(\speed_{\eps}\widetilde{X}^{\eps}\right)$ are exponentially equivalent and the LDP follows from~\cite[Theorem~4.2.13]{dembo1998large}.
\end{proof}


\subsubsection{Proof of Theorem~\ref{thm:LDP}}
\label{sec:thm:LDP}
Introduce $R_{\eps}:=\log\left(\frac{Z_{\eps}}{\eps}\right) - \varphi(0)$. 
Changing variables $s = \eps t$ in~\eqref{eq:Z} and using~\eqref{eq:rep_Xhat_T},
\begin{equation}\label{eq:Z_T_int}
\frac{Z_{\eps}}{\eps\E^{\varphi(0)}}
 = \int_{0}^{1}\exp\left\{\speed_{\eps}X^{\eps}_{t}\right\} \E^{\varphi(\eps t)-\varphi(0)} \D t.
\end{equation}
With~$\overline{\varphi}_{\eps}$
in~\eqref{eq:phi_notations},
we obtain the deterministic bound
\begin{equation}\label{eq:phi-bound}
\left|R_{\eps}-\Gg_{1}\left(\speed_{\eps}X^{\eps}\right)\right|\leq \overline{\varphi}_{\eps},
\end{equation}
so~$R_{\eps}$ and $\Gg_{1}(\speed_{\eps}X^{\eps})$ are exponentially equivalent at any speed. 
The map $\Gg:\Cc_0([0,1])\to\RR$ is continuous in the supremum norm, so the contraction principle~\cite[Theorem~4.2.1]{dembo1998large} applied to Theorem~\ref{thm:pathLDP} gives the LDP for $\Gg_{1}(\speed_{\eps}X^{\eps})$ with good rate function
$$
u\mapsto\inf_{f\in\Cc_{0}([0,1])}\left\{\If(f):\Gg_{1}(f)=u\right\}=\underline{\Jj}(u).
$$
Then~\cite[Theorem~4.2.13]{dembo1998large} yields an LDP for $(R_{\eps})_{\eps>0}$ on $\RR$ as $\eps\downarrow 0$, with speed $\speed_{\eps}^{-2}=(\nu^2\eps^{2H}\ell(\eps)^2)^{-1}$ and good rate function~$\underline{\Jj}$. 
The LDP for $(Z_{\eps}/\eps)_{\eps>0}$ follows by contraction.

\begin{proof}[Proof of Proposition~\ref{prop:laplace}]
By Theorem~\ref{thm:LDP} it suffices to check the conditions of Varadhan's lemma~\cite[Theorem~4.3.1]{dembo1998large}, that is, 
$$
\limsup_{\eps\downarrow 0}\speed_{\eps}^2\log\EE\left[\exp\left(\gamma|u|\speed_{\eps}^{-2}|R_{\eps}|\right)\right]<\infty, 
$$
for some $\gamma>1$. 
From~\eqref{eq:phi-bound} and $|\Gg_{1}(f)|\leq \|f\|_\infty$, 
then $|R_{\eps}|\leq \speed_{\eps}\|X^{\eps}\|_\infty+\overline{\varphi}_{\eps}$, so
$$
\speed_{\eps}^2\log\EE\left[\E^{c\speed_{\eps}^{-2}|R_{\eps}|}\right]\leq c \overline{\varphi}_{\eps}+\speed_{\eps}^2\log\EE\left[\E^{\frac{c}{\speed_{\eps}}\|X^{\eps}\|_\infty}\right]
\leq c \overline{\varphi}_{\eps}+c \speed_{\eps} \EE[\|X^{\eps}\|_\infty] + \frac{c^2\vartheta^2}{2},
$$
for any $c>0$,
where the last step is the Gaussian concentration bound 
$$
\EE\left[\E^{\beta\|Y\|_\infty}\right]\leq \exp\left\{\beta\EE[\|Y\|_\infty] + \frac{\beta^2\vartheta^2}{2}\right\}, 
$$
a consequence of Borell-TIS inequality~\cite[Theorem~2.1.1]{adler2007random} with $\beta=\frac{c}{\speed_{\eps}}$. 
The right-hand side is bounded uniformly in $\eps\leq \eps_{0}$ by Lemma~\ref{lem:holder}. 
Varadhan's lemma~\cite[Theorem~4.3.1]{dembo1998large} applies for every~$u$, and finiteness of~$\underline{\Jj}^{*}$ follows.
\end{proof}


\subsubsection{Proof of Lemma~\ref{lem:prop_lambda}}\label{sec:lem:prop_lambda_proof}
(i) $f\equiv0$ gives $\Gg_{1}(0)=0$ and $\If(0)=0$; positivity for $u\ne0$ follows since $\If$ is a good rate function and $\Gg^{-1}(\{u\})$ a closed set not containing~$0$.
Now, let $\Gg_{1}(f)=u'<u\le0$ with $\If(f)<\infty$. The map $\Gg_{1}(\cdot f)$ is continuous on $[0,1]$, so there exists $\theta\in(0,1]$ with $\Gg_{1}(\theta f)=u$; since $\If(\theta f)=\theta^2I(f)\leq I(f)$ then $\underline{\Jj}(u)\leq I(f)$, and taking the infimum over such $f$ gives 
$$
\underline{\Jj}(u)\leq \underline{\Jj}(u'). 
$$
The case $u\ge0$ is symmetric.
Finally, by Jensen's inequality, $\Gg_{1}(f)\geq \int_{0}^{1}f(s)\D s$, so $\Gg_{1}(f)=u\le0$ forces $\int_{0}^{1}f(s)\D s\leq u$; minimising $\If$ over the half-space $\{\langle\Kk^{*}\ind,\dot{h}\rangle\leq u\}$ gives $u^2/(2\|\Kk^{*}\ind\|^2)$ and $\|\Kk^{*}\ind\|^2=\EE[\Af^2]$.

(ii) $\Gg$ is convex by H\"older's inequality and so is~$\If$, 
hence $(f,u)\mapsto I(f)+\ind_{\{\Gg_{1}(f)\leq u\}}$ is jointly convex, and minimisation in~$f$ preserves convexity, yielding~$\widetilde{\underline{\Jj}}$. 
Moreover, 
$$
\widetilde{\underline{\Jj}}(u)=\inf_{u'\leq u}\underline{\Jj}(u'), 
$$
which equals $\underline{\Jj}(u)$ for $u\le0$ by~(i).
We now move on to~(iii) and first prove the behaviour at infinity.
For the lower bound, 
if $\Gg_{1}(f)=u$ then 
$$
u\leq \max_{w\in[0,1]}f(w)\leq \frac{1}{\sqrt{2H}}
\|\dot{h}\|_{L^2}, 
$$
by~\eqref{eq:CS}, and hence 
$\If(f)=\half\|\dot{h}\|^2\geq H u^2$.
For the upper bound, 
take $\dot{h}_\theta:=\theta\sqrt{2H}(1-v)^{H-\half}$, so that $\|\dot{h}_\theta\|^2=\theta^2$, $\If(\Kk\dot{h}_\theta)=\frac{\theta^2}{2}$, and $g:=\Kk\dot{h}_1$ satisfies
$$
g(t) = \sqrt{2H} t^{H+\half}\int_{0}^{1}(1-u)^{H-\half}(1-ut)^{H-\half}\D u,
$$
and is strictly increasing with $g(1) = \frac{1}{\sqrt{2H}}$. 
With the H\"older inequality following~\eqref{eq:CS}, $g(1)-g(t)\leq \frac{(1-t)^H}{\sqrt{H}}$, 
so for $\theta\geq 1$,
$$
\int_{0}^{1}\E^{\theta g(u)}\D u 
\geq \E^{\theta g(1)}\int_{0}^{1}
\exp\left\{-\frac{\theta u^{H}}{\sqrt{H}}\right\}\D u 
\geq \E^{\theta g(1)} c_H \theta^{-\frac{1}{H}},
$$
hence 
\begin{equation}\label{eq:G1_lowerbound}
\Gg_{1}(\theta g)\geq \frac{\theta}{\sqrt{2H}} - \frac{\log(\theta)}{H} - C_H,
\end{equation}
with $C_{H}:=-\log(c_{H})$. 
Since $\Gg_{1}(\cdot g)$ is continuous, increasing and unbounded
then for each large~$u$ there is $\theta_{u}$ 
with $\Gg_{1}(\theta_{u} g)=u$, and~\eqref{eq:G1_lowerbound} forces 
$\theta_{u}\leq \sqrt{2H} \left(u + \frac{\log(\theta_{u})}{H}+C_H
\right)$,
that is, 
$\theta_{u}=\sqrt{2H} u\left(1+\Oo\left(\frac{\log(u)}{u}\right)\right)$.
Then 
$$
\underline{\Jj}(u)\leq \frac{\theta_u^2}{2} = H u^2\left(1+\Oo\left(u^{-1}\log(u)\right)\right).
$$
Consider now the behaviour at the origin.
For the upper bound, take $f_\kappa:=\kappa \Kk\Kk^{*}\ind$, that is $\dot{h}_\kappa=\kappa\Kk^{*}\ind$, for which $\int_{0}^{1}f_\kappa(s)\D s=\kappa\EE[\Af^2]$ and $\If(f_\kappa)=\half\kappa^2\EE[\Af^2]$; since $\Gg_{1}(f)=\int_{0}^{1}f(s)\D s+\Oo(\|f\|_\infty^2)$, calibrating $\Gg_{1}(f_\kappa)=u$ gives $\kappa = \frac{u}{\EE[\Af^2]} + \Oo(u^2)$ and the stated upper bound. 
For the lower bound, Jensen's inequality reads 
$\Gg_{1}(f)\geq \int_{0}^{1}f(w)\D w$
and 
$$
\underline{\Jj}(u)\geq \inf\left\{\If(f):\int_{0}^{1}f(w)\D w\leq u\right\} = \frac{u^2}{2\EE[\Af^2]},
\qquad\text{for }u<0; 
$$
the case $u>0$ requires in addition the a priori bound $\|f\|_\infty=\Oo(|u|)$ on near-optimisers, which follows from~\eqref{eq:CS} and $\If(f)=\Oo(u^2)$.

We now prove the asymptotic behaviours and start by introducing the function
\begin{equation}\label{eq:cT}
\Lf(t):=\inf_{h\in\Hh_{t}}\left\{\If_{t}(\Kk\dot{h}): \Gg_{t}(\Kk \dot{h})\leq 0\right\},
\end{equation}
for $t\geq 1$, with~$\If_{t}$ in~\eqref{eq:I}
and~$\Gg_{t}$ in~\eqref{eq:Psi}.
Note that, from Lemma~\ref{lem:prop_lambda}(ii) and scaling,
\begin{equation}\label{eq:exact}
\underline{\Jj}(\log(z))
= \widetilde{\underline{\Jj}}(\log(z))
= \frac{1}{z^{2H}} \Lf\left(\frac{1}{z}\right),
\qquad\text{for every } z\in(0,1].
\end{equation}
The claim $\underline{\Jj}(u)\sim\Lfbar \E^{2H|u|}$
as $u\downarrow-\infty$ then follows immediately from the following:

\noindent\textbf{Claim:}
\emph{The map $\Lf(\cdot)$ is non-decreasing on $[1,\infty)$,
$\Lfbar := \lim_{t\uparrow\infty}\Lf(t) = \Lf(\infty)\in(0,\infty)$,
the infimum for~$\Lf(\infty)$
in~\eqref{eq:cT}
attained uniquely.}

We now prove the claim.
Let $1\leq t<t'$ and $\dot{h}$ admissible for~$t'$.
Then $\Kk(\dot{h}\ind_{[0,t]}) = \Kk\dot{h}$ on $[0,t]$ so that
$$
\int_{0}^{t}\E^{\Kk(\dot{h} \ind_{[0,t]})(s)}\D s
= \int_{0}^{t}\E^{\Kk(\dot{h})(s)}\D s
\leq \int_{0}^{t'}\E^{\Kk(\dot{h})(s)}\D s \leq 1,
$$
hence $\dot{h}\ind_{[0,t]}$ is admissible for~$\Lf(t)$
and
$\|\dot{h}\ind_{[0,t]}\|_{L^{2}([0,t])}
\leq \|\dot{h}\|_{L^{2}([0,t'])}$.
Taking infima yields 
$\Lf(t)\leq \Lf(t')\leq \Lf(\infty)$
and 
\[
\Lfbar
= \lim_{t\uparrow\infty}\Lf(t)
= \sup_{t\geq 1}\Lf(t)\leq \Lf(\infty),
\]
since~$\Lf$ is non-decreasing.
Following the proof of Proposition~\ref{prop:bounds_ch},  
$g_{\beta,c}(\cdot):=-\beta\log(1+c\cdot)$
with $\beta=\frac{1}{1-H}$
and $c=\frac{1}{\beta-1}$
satisfies $g_{\beta,c} = \Kk\dot{h}_{\beta,c}$
for some $\dot{h}_{\beta,c}\in L^{2}(0,\infty)$ with
$\|\dot{h}_{\beta,c}\|^{2}_{L^{2}(0,\infty)} = \beta^2 c^{2H}\Gamma_{H}$
(with~$\Gamma_{H}$ in Proposition~\ref{prop:bounds_ch})
and $\int_{0}^{\infty}\E^{g_{\beta,c}(s)}\D s=1$.
Therefore~$\dot{h}_{\beta,c}$ belongs to the admissible set (thus non empty) in the infimum defining~$\Lf(\infty)$.
Restricting this~$\dot{h}_{\beta,c}$ to $[0,t]$ shows that this also holds for the admissible sets 
in~\eqref{eq:cT}, and hence~$\Lf$ is finite.
Positivity follows from the proof of Proposition~\ref{prop:bounds_ch}.
Consider now a sequence $(t_{n})_{n\geq 1}\subset [1,\infty)$ with $\lim_{n\uparrow\infty}t_{n}=\infty$ and, for each~$n$,
let $\dot{h}_{n} \in L^{2}([0,t_n])$ be admissible for~$\Lf(t_n)$ with
$$
\half\|\dot{h}_{n}\|^2_{L^2([0,t_n])}\leq \Lf(t_n)+\frac{1}{t_n} \leq \Lfbar+1.
$$
We now extend each~$\dot{h}_{n}$ to $(0,\infty)$ by zero, 
hence the new sequence 
(with the same notation) 
$(\dot{h}_{n})_{n\geq 1}$ is bounded in $L^2(0,\infty)$,
hence converges weakly along a subsequence to some 
$\dot{h}^{*}\in L^2(0,\infty)$.
For each fixed $u$, $(u-\cdot)^{H-\half}\ind_{[0,u]}\in L^2(0,\infty)$, so $\Kk\dot{h}_{n}$ converges pointwise to~$\Kk\dot{h}^{*}$; 
For any $T>0$
and large enough~$n$ so that $t_n\geq T$, 
then
\[
\int_{0}^{T}\E^{(\Kk\dot{h}_{n})(u)}\D u
\leq \int_{0}^{t_{n}}\E^{(\Kk\dot{h}_{n})(u)}\D u\leq 1;
\]
Fatou's lemma on compacts then give 
$$
\int_{0}^{T}\E^{(\Kk\dot{h}^{*})(u)}\D u
\leq \liminf_{n\uparrow\infty}
\int_{0}^{T}\E^{(\Kk\dot{h}_{n})(u)}\D u\leq 1,
$$
and monotone convergence theorem yields
$
\int_{0}^{\infty}\E^{(\Kk\dot{h}^{*})(u)}\D u\leq 1$,
so that~$\dot{h}^{*}$ is admissible.
Weak lower semicontinuity of the norm finally gives 
$$
\Lf(\infty)
\leq \half\|\dot{h}^{*}\|^2_{L^2([0,\infty))}
\leq \liminf_{n\uparrow\infty}
\half\|\dot{h}_{n}\|^2_{L^2([0,\infty))}
\leq \liminf_{n\uparrow\infty}
\left(\Lf(t_{n})+\frac{1}{t_{n}}\right) = \Lfbar,
$$
and therefore
$\Lfbar=\Lf(\infty)$,
with the infimum attained at~$\dot{h}^{*}$. 
Uniqueness of the minimiser follows since the map 
$\dot{h}\mapsto \int_{0}^{\infty}\E^{(\Kk\dot{h})(u)}\D u$ is convex on $L^{2}(0,\infty)$.
Assume there exist two minimisers~$\dot{h}_{1}$
and~$\dot{h}_{2}$.
Then
$$
\half\left\|\frac{\dot{h}_{1}+\dot{h}_{2}}{2}\right\|^2
= \half\frac{\|\dot{h}_{1}\|^2 + \|\dot{h}_{2}\|^2}{2}
- \half\left\|\frac{\dot{h}_{1}-\dot{h}_{2}}{2}\right\|^2
 = \Lf(\infty) - \frac{1}{8}\left\|
\dot{h}_{1}-\dot{h}_{2}\right\|^2 <\Lf(\infty),
$$
which is a contradiction and the minimiser is unique.


\subsubsection{Proof of Proposition~\ref{prop:bounds_ch}}\label{sec:prop:bounds_ch_proof}
We start with the upper bound.
Recall the Riemann-Liouville fractional integral of order $H+\half>0$:
$$
(\Ii^{H+\half}f)(t) := \frac{1}{\Gamma(H+\half)}
\int_{0}^{t}(t-s)^{H-\half}f(s)\D s.
$$
Since $\Kk = \Gamma(H+\half)\Ii^{H+\half}$ and $\phi := \log(1+\cdot)$ is absolutely continuous with $\phi(0)=0$,
$$
\Kk^{-1}\phi = \Gamma\left(H+\half\right)^{-1}\Ii^{\half-H}\phi', 
$$
which is~\eqref{eq:GammaH} because $\phi'(u)=(1+u)^{-1}$; 
furthermore $\Ii^{H+\half}\Ii^{\half-H}\phi' = \Ii^{1}\phi' = \phi$ from~\cite[Chapter~2]{samko1993fractional}. 
As $u\downarrow 0$, 
$$
(\Kk^{-1}\phi)(u)\sim\frac{u^{\half-H}}{\Gamma(H+\half)\Gamma(\frac{3}{2}-H)}, 
$$
which is square integrable near~$0$; as $u\uparrow\infty$, 
$$
(\Kk^{-1}\phi)(u)=\Oo\left(u^{-H-\half}\log(u)\right), 
$$
which is square integrable at infinity since $H>0$. 
Hence $\Gamma_{H}:=\|\Kk^{-1}\log(1+\cdot)\|^2_{L^2(0,\infty)}$ is finite, with
\begin{equation}\label{eq:GammaH}
\left(\Kk^{-1}\log(1+\cdot)\right)(u)=\frac{1}{\Gamma(H+\half)\Gamma(\half-H)}\int_0^u\frac{(u-v)^{-H-\half}}{1+v} \D v.
\end{equation}

Fix $\beta>1$, $c>0$ and  $g_{\beta,c}(u):=-\beta\log(1+cu)$. 
Then
$$
\int_0^\infty\E^{g_{\beta,c}(u)}\D u = \int_0^\infty(1+cu)^{-\beta}\D u = \frac{1}{c(\beta-1)}, 
$$
so the choice $c=\frac{1}{\beta-1}$ 
makes~$g_{\beta,c}$ admissible. 
By scaling, $\|\Kk^{-1}g_{\beta,c}\|^2 = \beta^2c^{2H}\Gamma_{H}$, hence
$$
\Lfbar\leq \frac{\beta^2\Gamma_{H}}{2(\beta-1)^{2H}}.
$$
Minimising $2\log(\beta)-2H\log(\beta-1)$ gives $\beta^{*}=(1-H)^{-1}$ and the upper bound follows.

For the lower bound,
with $g=\Kk\dot{k}$ admissible and $t>1$,
Jensen's inequality yields
$$
\exp\left\{\frac{1}{t}
\int_0^t g(u)\D u\right\}
\leq \frac{1}{t}\int_0^t \E^{g(u)}\D u
\leq \frac{1}{t},
$$
so
$\langle\varpi,g\rangle_{L^2}\geq \log(t)>0$
with 
$\varpi:=-\frac{1}{t}\ind_{[0,t]}$.
On the other hand, 
$$
\left\langle\varpi,\Kk\dot{k}\right\rangle=\left\langle\Kk^{*}\varpi,\dot{k}\right\rangle\leq \|\Kk^{*}\varpi\|_{L^2}\|\dot{k}\|_{L^2}, 
$$
with 
$\|\Kk^{*}\varpi\|^2=\left\langle\varpi,\CC\varpi\right\rangle
= \VV\left[\frac{1}{t}\int_0^t X_u\D u\right]
 = t^{2H}\EE\left[\Af^2\right]$.
Therefore,
$$
\half\|\dot{k}\|^2 \geq \frac{\langle\varpi,g\rangle^2}{2\langle\varpi,\CC\varpi\rangle} \geq \frac{\log(t)^2}{2 \EE[\Af^2] t^{2H}}.
$$
Maximising over $t>1$ gives $\log(t)=\frac{1}{H}$ and the result follows from~\eqref{eq:AB}.

We now compute~$\Gamma_{H}$, providing the two equivalent representations for any $H\in(0,\half)$:
$$
\Gamma_{H}
  = \frac{1}{\Gamma(H+\half)^2}\int_{\RR}
\frac{|\Gamma(H+\I x)|^{2} \cosh(\pi x)}{\cosh(2\pi x)-\cos(2\pi H)} \D x\\
 = \frac{\Gamma(2H)}{2H\sin(\pi H)\Gamma(H+\half)^{2}}.
$$

Write $\psi:=\Kk^{-1}\phi$. Since by Riemann-Liouville fractional integral, 
$$
\Kk\dot{h} = \Gamma\left(H+\half\right)\Ii^{H+\half}\dot{h},
$$
and $\phi(0)=0$ with 
$\phi'(v)=(1+v)^{-1}$, then $\psi = \frac{1}{\Gamma(H+\half)}\Ii^{\half-H}\phi'$, that is
$$
\psi(u) = \frac{\cos(\pi H)}{\pi}\int_0^u\frac{(u-v)^{-H-\half}}{1+v} \D v
=\frac{\cos(\pi H)}{\pi} \frac{u^{\half-H}}{\half-H} {}_2F_1\left(1,1;\frac{3}{2}-H;-u\right),
\quad H<\half,
$$
where we used the reflection formula 
$$
\Gamma\left(H+\half\right)\Gamma\left(\half-H\right)
= \frac{\pi}{\sin(\pi(H+\half))}
 = \frac{\pi}{\cos(\pi H)},
$$ 
and Euler's integral. 
 Thus $\psi(u)\sim u^{\half-H}$ at $0$ and $\psi(u)\sim c u^{-H-\half}\log(u)$ as $u\uparrow\infty$, 
 so $\psi\in L^2((0,\infty))$ and its Mellin transform 
$$ 
(\Mm\psi)(t) := \int_0^\infty\psi(u)u^{t-1}\D u 
$$
is analytic in the strip $
 \Re(t) \in (H-\half,H+\half)$.
Now, 
$\Mm[(1+\cdot)^{-1}](t) = \frac{\pi}{\sin(\pi t)}$ for $0<\Re(t)<1$ and, using Fubini's theorem and the Beta integral,
$$
\Mm\left[\Ii^{\beta}\gamma\right](t) = \frac{\Gamma(1-t-\beta)}{\Gamma(1-t)} M\gamma(t+\beta),
\qquad \Re(t)<1-\beta,
$$
applied with $\beta = \half-H \in [0,\half)$ gives, using the reflection formula,
\begin{equation}\label{eq:Mpsi}
\Mm\psi(t) = \frac{\Gamma(H+\half-t)^{2} \Gamma(\half-H+t)}{\Gamma(H+\half) \Gamma(1-t)}.
\end{equation}
By Mellin-Plancherel theorem, 
$$
\Gamma_{H} = \frac{1}{2\pi}\int_{\RR}\left|\Mm\psi\left(\half+\I t\right)\right|^2\D t.
$$
On the critical line, $H+\half-t = H-\I t$ and $\half-H+t=1-H+\I t$, we can compute
$$
\left|\Gamma\left(\half-\I t\right)\right|^{2}
 = \frac{\pi}{\cosh(\pi t)}
 \quad\text{and}\quad
\left|\Gamma(H-\I t) \Gamma(1-H+\I t)\right|^{2}
=\frac{\pi^{2}}{\sin^{2}(\pi H)+\sinh^{2}(\pi t)}.
$$
Substituting~\eqref{eq:Mpsi} and simplifying yields the claimed expression for~$\Gamma_{H}$.
Since the integrand decays like $|t|^{2H-1}\E^{-2\pi|t|}$, the integral converges.
The cases $H=\half$ 
and $H\downarrow 0$ follow by direct computations.
Let now $\gamma(v):=(1+v)^{-1}$.
By~\eqref{eq:GammaH}, 
since all the integrands are non-negative, Tonelli's theorem yields
$$
\Gamma_{H}
= \frac{1}{\Gamma(H+\half)^{2}\Gamma(\half-H)^{2}}
\int_{0}^{\infty}\int_{0}^{\infty}\gamma(v)\gamma(v')
\left(\int_{v\vee v'}^{\infty}
\frac{\D u}{(u-v)^{H+\half}(u-v')^{H+\half}}
\right)\D v\D v'.
$$
For $v<v'$, the substitution $u=v'+(v'-v)w$ gives
$$
\int_{v'}^{\infty}
\frac{\D u}{(u-v)^{H+\half}(u-v')^{H+\half}}
= (v'-v)^{-2H}\int_{0}^{\infty}
\frac{\D w}{w^{H+\half}(1+w)^{H+\half}}
= \frac{\mathrm{B}(\beta,2H)}{(v'-v)^{2H}},
$$
with $\mathrm{B}(\half-H,2H) = \frac{\Gamma(\half-H)\Gamma(2H)}{\Gamma(H+\half)}$, so that, by symmetry,
$$
\Gamma_{H}
= \frac{2\Gamma(2H)}{\Gamma(H+\half)^{3}\Gamma(\half-H)}
\int_{0}^{\infty}\gamma(v)\left(\int_{v}^{\infty}
\frac{\gamma(v')}{(v'-v)^{2H}}\D v'\right)\D v.
$$
With $v'=v+(1+v)x$, the inner integral is equal to
$\frac{1}{(1+v)^{2H}}\int_{0}^{\infty}\frac{x^{-2H}}{1+x}\D x
= \frac{1}{(1+v)^{2H}}
\frac{\pi}{\sin(2\pi H)}$,
and $\int_{0}^{\infty}\frac{\D v}{(1+v)^{1+2H}} = \frac{1}{2H}$, which proves the claim.
The equivalent form follows from the reflection formula $\Gamma(H+\half)\Gamma(\half-H)=\frac{\pi}{\cos(\pi H)}$.

\subsection{Proofs of the CLT and MDP}\label{sec:proofs}

With~$\eps_{0}$ and~$\overline{m}$ from Lemma~\ref{lem:holder} and~$\varrho_{\eps}^2$
from Lemma~\ref{lem:coupling},
together with $\vartheta_{0}^{2} := \sup_{t\leq 1}\VV[X_t] = \frac{1}{2H}$, 
and $\overline{\speed} := \sup_{\eps\leq \eps_{0}}\speed_{\eps} < \infty$.
The following lemma is straightforward.

\begin{lemma}\label{lem:tails}
Let $Y$ be a centred Gaussian process with continuous paths on $[0,1]$, 
$\mathfrak{v} := \sup_{t\leq 1}\VV[Y_t]$ and $m := \EE[\|Y\|_{\infty}]$, both finite. 
Then 
$$
\EE\left[\|Y\|_{\infty}^p \E^{c\|Y\|_{\infty}}\right] \leq \Psi\left(p,c,m,\sqrt{\mathfrak{v}}\right)<\infty,
$$
for all $p,c\geq 0$, with~$\Psi$ non-decreasing in~$m$ and~$\mathfrak{v}$ (in~$m$ when $p=0$ or $p\geq 1$).
\end{lemma}

\begin{proof}
Since $\|Y\|_{\infty} \leq m + (\|Y\|_{\infty} - m)_{+}$, then
$\PP\left((\|Y\|_{\infty} - m)_{+}\geq x\right)\leq \exp\left\{-\frac{x^2}{2\mathfrak{v}}\right\}$
by Borell-TIS~\cite[Theorem~2.1.1]{adler2007random}. 
Since $f:x\mapsto(m+x)^p \E^{c(m+x)}$ is non-decreasing, then
\begin{align*}
\EE\left[\|Y\|_{\infty}^p \E^{c\|Y\|_{\infty}}\right]
 &\leq \EE[f((\|Y\|_{\infty} - m)_{+})]
 \\
 & = f(0) + \int_0^\infty f'(x) \PP((\|Y\|_{\infty} - m)_{+}\geq x)\D x\\
 & \leq f(0) + \int_0^\infty f'(x) \exp\left\{-\frac{x^2}{2\mathfrak{v}}\right\}\D x
  =: \Psi\left(p,c,m,\sqrt{\mathfrak{v}}\right).
\end{align*}
The function~$\Psi$ is non-decreasing in~$\sqrt{\mathfrak{v}}$ for all $p,c\geq 0$, and in~$m$ when $p=0$ or $p\geq 1$, since then both~$f(0)$ and~$f'$ are non-decreasing in~$m$; for $p\in(0,1)$ this may fail, but $\|Y\|_{\infty}^{p}\leq 1+\|Y\|_{\infty}$ reduces this case to $p\in\{0,1\}$.
\end{proof}
Applying Lemma~\ref{lem:tails} to~$X^{\eps}$, 
$\widetilde{X}^{\eps}$ and $\Df^{\eps} := X^{\eps}-\widetilde{X}^{\eps}$ yield
\begin{equation}\label{eq:unifexp}
\sup_{\eps\leq \eps_{0}}\EE\left[\|X^{\eps}\|_\infty^{p}  \E^{c\|X^{\eps}\|_\infty}\right]
+ \sup_{\eps\leq \eps_{0}}\EE\left[\|\widetilde{X}^{\eps}\|_\infty^p\right]
+ \sup_{\eps\leq \eps_{0}}\EE\left[\|\Df^{\eps}\|_\infty^p\right] < \infty,
\end{equation}
for all $p,c\geq 0$ since $\|\Df^{\eps}\|_\infty\leq \|X^\eps\|_\infty+\|\widetilde{X}^{\eps}\|_\infty$. 
Since $\varrho_{\eps}\to 0$ and $\EE[\|\Df^{\eps}\|_\infty]\to 0$, Lemma~\ref{lem:tails} applied to~$\Df^{\eps}$ 
with $(m,\mathfrak{v}) = \left(\EE[\|\Df^{\eps}\|_\infty],\varrho_{\eps}^2\right)$ further gives
\begin{equation}\label{eq:Dtozero}
\lim_{\eps\downarrow 0}\EE\left[\|\Df^{\eps}\|_\infty^p\right] = 0, 
\quad \text{for all }p>0.
\end{equation}

Recall from Section~\ref{sec:thm:LDP} and~\eqref{eq:Z_T_int} that
\begin{equation}\label{eq:YT}
Y_{\eps}
 := \frac{Z_\eps}{\eps\E^{\varphi(0)}}
 = \int_{0}^{1} \exp\left\{\speed_{\eps}X^{\eps}_{t} + \varphi_{\eps}(t)\right\}\D t.
\end{equation}
Set
\begin{equation}\label{eq:AT}
\Af^{(\eps)} := \int_{0}^{1} \widetilde{X}^{\eps}_s\D s  = \Af\sim \Nn\left(0,\EE[\Af^2]\right),
\end{equation}
where the equality holds in distribution since~$\widetilde{X}^{\eps}$ and~$X$ 
have the same law for every~$\eps$, and $\EE[\Af^2]
= \|\Kk^{*}\ind\|_{L^2([0,1])}^2$ 
with 
$$
(\Kk^{*}\ind)(u)
= \int_u^1 (t-u)^{H-\half}\D t
= \frac{(1-u)^{H+\half}}{H+\half}.
$$

\begin{lemma}\label{lem:decomp}
For every $\eps>0$,
\begin{equation}\label{eq:decomp}
Y_\eps = \Phi_{\eps} + \speed_{\eps}\left(\Af^{(\eps)} + r^{(1)}_\eps + r^{(2)}_\eps + r^{(3)}_\eps\right),
\end{equation}
with the remainders 
\begin{equation}\label{r:3:defn}
r^{(1)}_{\eps}
:= \int_{0}^{1} \E^{\varphi_{\eps}(t)} \Df^{\eps}_{t}\D t,
\quad
r^{(2)}_{\eps}
:= \int_{0}^{1}\left(\E^{\varphi_{\eps}(t)}-1\right) \widetilde{X}^{\eps}_{t}\D t,
\quad r^{(3)}_{\eps}
:= \frac{1}{\speed_{\eps}}\int_{0}^{1} \E^{\varphi_{\eps}(t)}\rho_{\eps}(t)\D t,
\end{equation}
satisfying, almost surely, 
\begin{equation}\label{eq:rbounds}
\left|r^{(1)}_{\eps}\right| \leq \E^{\overline{\varphi}_{\eps}}\|\Df^{\eps}\|_\infty,
\quad
\left|r^{(2)}_{\eps}\right| \leq \left|\E^{\overline{\varphi}_{\eps}}-1\right|\|\widetilde{X}^{\eps}\|_\infty,
\quad
\left|r^{(3)}_{\eps}\right|
\leq \frac{\speed_{\eps} \E^{\overline{\varphi}_{\eps}}}{2}
\|X^{\eps}\|_\infty^2  \E^{\speed_{\eps}\|X^{\eps}\|_\infty}.
\end{equation}
\end{lemma}

\begin{proof}
Since $|\E^x-1-x|\leq \frac{x^2}{2}\E^{|x|}$ for all $x\in\RR$, 
then 
$\E^{\speed_{\eps}X^{\eps}_{t}} = 1 + \speed_{\eps}X^{\eps}_{t} + \rho_{\eps}(t)$,
with $|\rho_{\eps}(t)| \leq \frac{\speed_{\eps}^2}{2}\|X^{\eps}\|_\infty^2 \E^{\speed_{\eps}\|X^{\eps}\|_\infty}$. 
Inserting this into~\eqref{eq:YT} yields
$Y_{\eps} - \Phi_{\eps} = \speed_{\eps}\int_{0}^{1} \E^{\varphi_{\eps}(t)}X^{\eps}_t\D t + \int_{0}^{1} \E^{\varphi_{\eps}(t)}\rho_{\eps}(t)\D t$.
By definition, $X^{\eps} = \widetilde{X}^{\eps} + \Df^{\eps}$, and therefore
\begin{align*}
\int_{0}^{1} \E^{\varphi_{\eps}(t)}X^{\eps}_{t}\D t 
&= \int_{0}^{1}\widetilde{X}^{\eps}_{t}\D t
 + \int_{0}^{1} \E^{\varphi_{\eps}(t)} \Df^{\eps}_{t}\D t
 + \int_{0}^{1}\left(\E^{\varphi_{\eps}(t)}-1\right) \widetilde{X}^{\eps}_{t}\D t
=: \Af^{(\eps)} + r^{(1)}_{\eps} + r^{(2)}_{\eps}.
\end{align*}
The bounds~\eqref{eq:rbounds} then follow from $\|\varphi_{\eps}\|_\infty = \overline{\varphi}_{\eps}$ and $|\E^{\varphi_{\eps}}-1|\leq \E^{\overline{\varphi}_{\eps}}-1$.
\end{proof}

\subsubsection{Proof of Theorem~\ref{thm:CLT}}

Fix $p\ge1$. By~\eqref{eq:unifexp}-\eqref{eq:Dtozero}-\eqref{eq:rbounds},
$$
\EE\left[\left|r^{(1)}_{\eps}\right|^p\right] \leq \E^{p\overline{\varphi}_{\eps}}\EE\left[\|\Df^{\eps}\|_\infty^p\right]
\qquad\text{and}\qquad
\EE\left[\left|r^{(2)}_{\eps}\right|^p\right] \leq \left|\E^{\overline{\varphi}_{\eps}}-1\right|^p \EE\left[\|X\|_\infty^p\right],
$$
with both right-hand sides tending to zero,
and, since $\speed_{\eps}\leq \overline{\speed}$ on $(0,\eps_{0}]$,
$$
\EE\left[\left|r^{(3)}_{\eps}\right|^p\right] \leq \left(\frac{\speed_{\eps} \E^{\overline{\varphi}_{\eps}}}{2}\right)^p 
\EE\left[\|X^{\eps}\|_\infty^{2p}\E^{p\overline{\speed}\|X^{\eps}\|_\infty}\right]
\leq \left(\frac{\speed_{\eps} \E^{\overline{\varphi}_{\eps}}}{2}\right)^p \Psi\left(2p,p\overline{\speed},\overline{m},\vartheta\right) \to 0,
$$
with~$\vartheta$ in Lemma~\ref{lem:holder}.
Hence $r_{\eps} := r^{(1)}_{\eps} + r^{(2)}_{\eps} + r^{(3)}_{\eps}$ tends to zero in~$L^p$ for every $p$, 
hence in probability. 
Since $\Af^{(\eps)}\sim\Nn(0,\EE[\Af^2])$ for every~$\eps$ by~\eqref{eq:AT}, 
Slutsky's lemma~\cite[Lemma~2.8]{van2000asymptotic} applied to~\eqref{eq:decomp} implies that
$\frac{Y_{\eps} - \Phi_{\eps}}{\speed_{\eps}} = \Af^{(\eps)} + r_{\eps}$
converges in distribution to $\Nn\left(0,\EE[\Af^2]\right)$,
and~\eqref{eq:CLT} follows recalling~\eqref{eq:phi_notations}. 
For the moments,~\eqref{eq:decomp} together with the Gaussian moments of $\Af^{(\eps)}$ and the~$L^{2p}$ bounds above show that 
$$
\sup_{\eps\leq \eps_{0}}\EE\left[\left|\speed_{\eps}^{-1}(Y_{\eps} - \Phi_{\eps})\right|^{2p}\right]<\infty.
$$
Hence $\left\{|\speed_{\eps}^{-1}(Y_{\eps} - \Phi_{\eps})|^p\right\}_{\eps\leq \eps_{0}}$ is uniformly integrable and the $p$-th absolute moments converge~\cite[Theorem~2.20]{van2000asymptotic}. 

\subsubsection{Proof of Theorem~\ref{thm:MDP}}

Set $\kappa_{\eps} := \frac{\speed_{\eps}}{a_\eps}\downarrow 0$
so $\tau_{\eps} := \frac{1}{\kappa_{\eps}^2}\uparrow\infty$ and $\frac{a_\eps}{\speed_{\eps}^{2}} \uparrow\infty$. 
By~\eqref{eq:phi_notations} and Lemma~\ref{lem:decomp},
\begin{equation}\label{eq:MDPdecomp}
\frac{1}{a_\eps}\left(\frac{Z_\eps}{\eps}-\E^{\varphi(0)}\Phi_{\eps}\right)
= \E^{\varphi(0)} \kappa_{\eps} \Af^{(\eps)} + \E^{\varphi(0)} \kappa_{\eps}\left(r^{(1)}_\eps + r^{(2)}_\eps + r^{(3)}_{\eps}\right).
\end{equation}

Starting with the first term on the right-hand side, 
the law~$\mu$ of~$X$ is a centred Gaussian measure on the separable Banach space $\Cc_0([0,1])$, with Cameron-Martin space $\Kk(L^2([0,1]))$, 
with $\|\Kk\dot{h}\|_{\Hh_{1}^{X}} = \|\dot{h}\|_{L^2}$
(Theorem~\ref{thm:pathLDP} or~\cite{decreusefond1999stochastic}). 
Thus $(\theta X)_{\theta>0}$ 
satisfies an LDP~\cite[Section~3.4]{deuschel2001large} as $\theta\downarrow 0$ 
with speed~$\theta^{-2}$ and good rate function $\If:=\half\|\cdot\|^2_{\Hh_{1}^{X}}$.
Since~$\widetilde{X}^{\eps}$ and~$X$ have the same law for every~$\eps$, 
then $(\kappa_{\eps}\widetilde{X}^{\eps})_{\eps>0}$ satisfies an LDP on $\Cc_0([0,1])$ 
with speed~$\tau_{\eps}$ and good rate function~$\If$. 
The map $\Theta: f\mapsto \E^{\varphi(0)}\int_{0}^{1} f(s)\D s$ is continuous on~$\Cc_0([0,1])$, so the contraction principle~\cite[Theorem~4.2.1]{dembo1998large} yields an LDP for 
$(\E^{\varphi(0)} \kappa_{\eps} \Af^{(\eps)})_{\eps>0} = (\Theta(\kappa_{\eps} \widetilde{X}^{\eps}))_{\eps>0}$ with speed~$\tau_{\eps}$ and good rate function
$$
\inf\{\If(f): \Theta(f) = u\}
= \inf\left\{\half\|\dot{h}\|_{L^2}^2: \left\langle \Kk^{*}\ind,\dot{h}\right\rangle_{L^2} = u \E^{-\varphi(0)}\right\}
= \frac{u^2\E^{-2\varphi(0)}}{2\|\Kk^{*}\ind\|_{L^2}^2} =: \Jj(u).
$$
By~\eqref{eq:MDPdecomp} and~\cite[Theorem~4.2.13]{dembo1998large}, 
it suffices to prove that for each $i\in\{1,2,3\}$, $\delta>0$,
$$
\limsup_{\eps\downarrow 0}
\frac{1}{\tau_\eps}\log\PP\left(\kappa_{\eps}|r^{(i)}_\eps| > \delta\right) = -\infty.
$$
In that case, 
$$
\PP\left(\E^{\varphi(0)}\kappa_{\eps}\left|r^{(1)}_\eps + r^{(2)}_\eps + r^{(3)}_\eps\right| > \delta\right)\leq \sum_{i=1}^{3}\PP\left(\kappa_{\eps}|r^{(i)}_\eps|>\frac{\delta}{3} \E^{-\varphi(0)}\right),
$$
and the two sides of~\eqref{eq:MDPdecomp} are hence exponentially equivalent.

$\bullet$ By~\eqref{eq:rbounds}, $\kappa_{\eps}|r^{(1)}_\eps|>\delta$ 
so that 
$\|\Df^{\eps}\|_\infty > \delta \E^{-\overline{\varphi}_{\eps}}\frac{a_\eps}{\speed_{\eps}} \geq \frac{\delta}{2} \frac{a_\eps}{\speed_{\eps}}$.
Since $\frac{a_\eps}{\speed_{\eps}}\uparrow\infty$ 
and $\EE[\|\Df^{\eps}\|_\infty]\downarrow 0$, for~$\eps$ small enough then
$\EE[\|\Df^{\eps}\|_\infty]\leq \frac{\delta}{4} \frac{a_\eps}{\speed_{\eps}}$ 
and Borell-TIS inequality gives, 
as $\eps\downarrow 0$, since $|\varrho_\eps|\downarrow 0$,
$$
\frac{1}{\tau_{\eps}}\log\PP\left(\|\Df^{\eps}\|_\infty\geq \frac{\delta}{2}\frac{a_\eps}{\speed_{\eps}}\right)
\leq -\frac{\speed_{\eps}^2}{a_\eps^2}\frac{\left(\frac{\delta}{2}\frac{a_\eps}{\speed_{\eps}}-\EE[\|\Df^{\eps}\|_\infty]\right)^2}{2\varrho_\eps^2}
\leq -\frac{\speed_{\eps}^2}{a_\eps^2}\frac{\left(\frac{\delta}{4}\frac{a_\eps}{\speed_{\eps}}\right)^2}{2\varrho_\eps^2}
= -\frac{\delta^2}{32 \varrho_\eps^2} \downarrow-\infty.
$$

$\bullet$ If $\overline{\varphi}_{\eps}=0$ then $r^{(2)}_\eps=0$. 
Otherwise, by~\eqref{eq:rbounds}, $\kappa_{\eps}|r^{(2)}_\eps|>\delta$ implies 
$$
\|\widetilde{X}^{\eps}\|_\infty > M_\eps := \frac{\delta a_\eps}{\speed_{\eps}(\E^{\overline{\varphi}_{\eps}}-1)}\uparrow\infty,
$$
since $\frac{a_{\eps}}{\speed_{\eps}}\uparrow\infty$ and $\overline{\varphi}_{\eps}\downarrow 0$. 
For~$\eps$ small, $\EE\left[\|\widetilde{X}^{\eps}\|_\infty\right]\leq \frac{M_\eps}{2}$, 
and Borell-TIS inequality gives,
as $\eps\downarrow 0$, since $\overline{\varphi}_{\eps}\downarrow 0$,
$$
\frac{1}{\tau_{\eps}}\log\PP\left(\|\widetilde{X}^{\eps}\|_\infty > M_\eps\right)
\leq -\frac{\speed_{\eps}^2}{a_\eps^2}\frac{(\frac{M_\eps}{2})^2}{2\vartheta_{0}^{2}}
= -\frac{\delta^2}{8\vartheta_{0}^{2} (\E^{\overline{\varphi}_{\eps}}-1)^2}
\downarrow-\infty.
$$

$\bullet$ By~\eqref{eq:rbounds}, $\kappa_{\eps}|r^{(3)}_\eps|>\delta$ implies
$$
\|X^{\eps}\|_\infty^2 \E^{\speed_{\eps}\|X^{\eps}\|_\infty} > \frac{2\delta \E^{-\overline{\varphi}_{\eps}}a_\eps}{\speed_{\eps}^2}\geq y_\eps := \frac{\delta a_{\eps}}{\speed_{\eps}^2}.
$$
For~$\eps$ small, $\delta a_{\eps}\leq 1$ hence $\sqrt{y_\eps} = \frac{\sqrt{\delta a_{\eps}}}{\speed_{\eps}}\leq \frac{1}{\speed_{\eps}}$ and, for $x\geq 0$,
$x^2\E^{\speed_{\eps}x}\geq y_\eps$ implies $x\geq \E^{-\half}\sqrt{y_\eps}$,
thus $\|X^\eps\|_\infty\geq \E^{-\half}\sqrt{y_\eps}$. 
Since $y_\eps = \delta a_{\eps}\speed_{\eps}^{-2}\uparrow\infty$ for~$\eps$ small,
then $\overline{m}\geq \EE[\|X^\eps\|_\infty]$ 
and $\overline{m}\leq \half\E^{-\half}\sqrt{y_\eps}$, and Borell-TIS inequality gives, as $\eps\downarrow 0$,
$$
\frac{1}{\tau_\eps}\log\PP\left(\|X^\eps\|_\infty\geq \E^{-\half}\sqrt{y_{\eps}}\right)
\leq -\frac{\speed_{\eps}^2}{a_{\eps}^2}\frac{\left(\half\E^{-\half}\sqrt{y_{\eps}}\right)^2}{2\vartheta^2}
= -\frac{\speed_{\eps}^2}{a_{\eps}^2}\frac{y_{\eps}}{8\E \vartheta^2}
= -\frac{\delta}{8\E \vartheta^2 a_{\eps}}
\downarrow -\infty.
$$


\subsection{Proofs for Section~\ref{sec:options}}\label{sec:proofs:options}

\subsubsection{Proof of Theorem~\ref{thm:OTM}}
First, for any $p>1$, under Assumption~\ref{assump:1},
we can compute
\begin{align}
\EE\left[\left|\frac{1}{T}\int_{0}^{T}\sigma^{2}(Y_{t})\D t-K\right|^{p}\right]
&\leq
\EE\left[\left(\frac{1}{T}\int_{0}^{T}\sigma^{2}(Y_{t})\D t+K\right)^{p}\right]
\nonumber
\\
&\leq
2^{p-1}\left(\EE\left[\left(\frac{1}{T}\int_{0}^{T}\sigma^{2}(Y_{t})\D t\right)^{p}\right]+K^{p}\right)
\nonumber
\\
&\leq
2^{p-1}\left(\EE\left[\frac{1}{T}\int_{0}^{T}\sigma^{2p}(Y_{t})\D t\right]+K^{p}\right)
\nonumber
\\
&\leq
2^{p-1}\left(\frac{1}{T}\int_{0}^{T}\EE\left[\theta_{1}^{p}\E^{p\theta_{2}|Y_{t}|}\right]\D t+K^{p}\right),
\end{align}
where we applied Jensen's inequality and Assumption~\ref{assump:1}.
Furthermore, since $\theta_{1},\theta_{2}>0$,
\begin{align}
\frac{1}{T}\int_{0}^{T}\EE\left[\theta_{1}\E^{p\theta_{2}|Y_{t}|}\right]\D t
&\leq
\frac{\theta_{1}}{T}\int_{0}^{T}\EE\left[\E^{p\theta_{2}Y_{t}}+\E^{-p\theta_{2}Y_{t}}\right]\D t
\nonumber
\\
&=\frac{\theta_{1}}{T}\int_{0}^{T}\left(\E^{p\theta_{2}y_{0}}+\E^{-p\theta_{2}y_{0}}\right)\E^{\half p^{2}\theta_{2}^{2}t^{2H}}\D t
\nonumber
\\
&\leq\theta_{1}\left(\E^{p\theta_{2}y_{0}}+\E^{-p\theta_{2}y_{0}}\right)\E^{\half p^{2}\theta_{2}^{2}},
\end{align}
for any $T\leq 1$.
With $\Sigma_{T}:=\frac{1}{T}\int_{0}^{T}\sigma^{2}(Y_{t})\D t$, then $\Call(K,T)=\EE[(\Sigma_{T}-K)_{+}]$; the bounds above show that, for every $q>1$, 
$c_{q}:=\sup_{T\leq 1}\EE\left[|\Sigma_{T}-K|^{q}\right]$ is finite.
Since $(\Sigma_{T}-K)_{+}\leq |\Sigma_{T}-K|\ind_{\{\Sigma_{T}\geq K\}}$, H\"older's inequality yields the upper bound 
\[
\Call(K,T)\leq c_{q}^{\frac{1}{q}}\left(\PP(\Sigma_{T}\geq K)\right)^{1-\frac{1}{q}},
\]
while 
\[
\Call(K,T)\geq \delta\PP(\Sigma_{T}\geq K+\delta),
\]
for every $\delta>0$.
Therefore, for every $q>1$ and $\delta>0$,
\begin{align}
\liminf_{T\downarrow 0}T^{2H}\log\PP\left(\Sigma_{T}\geq K+\delta\right)
&\leq \liminf_{T\downarrow 0}T^{2H}\log \Call(K,T)
\nonumber
\\
&\leq \limsup_{T\downarrow 0}T^{2H}\log \Call(K,T)
\nonumber
\\
&\leq \left(1-\frac{1}{q}\right)\limsup_{T\downarrow 0}T^{2H}\log\PP\left(\Sigma_{T}\geq K\right),
\end{align}
and, letting $q\uparrow\infty$ and $\delta\downarrow 0$,
since continuity holds whenever~$\sigma^{2}$ is strictly monotone or continuously differentiable (Lemma~\ref{lem:Lambda_cont}(ii), or assumed in Theorem~\ref{thm:OTM}), 
it suffices to determine the limit of $T^{2H}\log\PP\left(\Sigma_{T}\geq K\right)$.
With the change of variables $t\mapsto tT$, this probability can be re-written as
\begin{equation}
\PP\left(\Sigma_{T}\geq K\right)
=\PP\left(\int_{0}^{1}\sigma^{2}(Y_{tT})\D t\geq K\right).
\end{equation}
By~\cite[Theorem~3.4.12]{deuschel2001large} $\PP\left(\left\{\eps^{H}W_{t}^{H},0\leq t\leq 1\right\}\in\cdot\right)$ satisfies
an LDP with good rate function $I_{H}(g)$
and speed $\eps^{-2H}$,
where
\begin{equation}
I_{H}(g):=\half\int_{0}^{1}|\dot{h}(t)|^{2}\D t,
\end{equation} 
with $g(t)=\int_{0}^{t}K_{H}(s,t)\dot{h}(s)\D s$ for any $0\leq t\leq 1$.
By self-similarity, 
$Y_{tT}=y_{0}+W_{tT}^{H}=y_{0}+T^{H}W_{t}^{H}$
in distribution.
Moreover,
$\PP\left(\left\{T^{H}W^{H}_{t},0\leq t\leq 1\right\}\in\cdot\right)$ satisfies
an LDP with good rate function $I_{H}(g)$
and speed~$T^{-2H}$. 
Since $\sigma:\RR\to[0,\infty)$ is continuous, by contraction principle, 
\begin{align*}
\lim_{T\downarrow 0}T^{2H}\log\PP\left(\int_{0}^{1}\sigma^{2}(Y_{tT})\D t\geq K\right)
&=-\inf_{\int_{0}^{1}\sigma^{2}(y_{0}+g(t))\D t=K}I_{H}(g)
\\
&=-\inf_{\int_{0}^{1}\sigma^{2}\left(y_{0}+(\Kk_{H}\dot{h})(t)\right)\D t=K}\half\|h\|^{2}_{\Hh_{1}}.
\end{align*}
More precisely, $F$ being continuous, the large deviations upper bound on the closed set $\{F\geq K\}$ and the lower bound on the open set $\{F>K+\delta\}$ give, by Lemma~\ref{lem:Lambda_cont}(i),
\begin{align*}
\limsup_{T\downarrow 0}T^{2H}\log\PP\left(\Sigma_{T}\geq K\right)
 & \leq -\inf\left\{I_{H}(g): F(g)\geq K\right\} = -\Lambda(K),\\
\liminf_{T\downarrow 0}T^{2H}\log\PP\left(\Sigma_{T}\geq K+\delta\right)
 & \geq -\inf\left\{I_{H}(g): F(g)> K+\delta\right\} \geq -\Lambda(K+2\delta),
\end{align*}
for every $\delta>0$, and the right-continuity of~$\Lambda$ at~$K$ concludes as~$\delta\downarrow 0$.
If $\Lambda(K)=+\infty$, the upper bound alone gives the result.
Similarly, for Put options, we can show that
\begin{equation}
\lim_{T\downarrow 0}T^{2H}\log \Put(K,T)
=\lim_{T\downarrow 0}T^{2H}\log\PP\left(\int_{0}^{1}\sigma^{2}(Y_{tT})\D t\leq K\right),
\end{equation}
and by the same large deviations argument as before, 
\begin{equation}
\lim_{T\downarrow 0}T^{2H}\log\PP\left(\int_{0}^{1}\sigma^{2}(Y_{tT})\D t\leq K\right)
=-\inf_{\int_{0}^{1}\sigma^{2}\left(y_{0}+(\Kk_{H}\dot{h})(t)\right)\D t=K}\half \|h\|^{2}_{\Hh_{1}}.
\end{equation}
Here again the upper bound is taken on the closed set $\{F\leq K\}$, the lower bound on the open set $\{F<K-\delta\}$, and the left-continuity of~$\Lambda$ at~$K$ concludes.

\subsubsection{Proof of Corollary~\ref{cor:OTMlognormal}}
From the proof of Lemma~\ref{lem:moments}, for every $q\geq 1$ and $\eps\leq 1$,
$$
\EE\left[\left(\frac{Z_{\eps}}{\eps}\right)^{q}\right]
\leq \frac{1}{\eps}\int_{0}^{\eps}\E^{q\varphi(u)+\half q^{2}\nu^{2}\VV[X^{(\ell)}_{u}]}\D u
\leq \sup_{u\leq 1}\exp\left\{q\varphi(u)+\half q^{2}\nu^{2}\VV\left[X^{(\ell)}_{u}\right]\right\},
$$
which is finite and independent of~$\eps$, so that $\sup_{\eps\leq 1}\EE[|\frac{Z_{\eps}}{\eps}-K|^{q}]<\infty$.
The H\"older sandwich from the proof of Theorem~\ref{thm:OTM} therefore reduces the claim to the tail asymptotics of~$Z_{\eps}/\eps$.
For Call options, let $K>\sigma_{0}^{2}$ and $u:=\log(\frac{K}{\sigma_{0}^{2}})>0$.
Theorem~\ref{thm:LDP} and the monotonicity of~$\underline{\Jj}$ on $[0,\infty)$ (Lemma~\ref{lem:prop_lambda}(i)) give
$$
\limsup_{\eps\downarrow 0}\speed_{\eps}^{2}\log\PP\left(\frac{Z_{\eps}}{\eps}\geq K\right)
\leq -\inf_{z\geq K}\underline{\Jj}\left(\log\frac{z}{\sigma_{0}^{2}}\right)
= -\underline{\Jj}(u).
$$
For the lower bound, fix $\varsigma>0$ and $f$ with $\Gg_{1}(f)=u$ and $\If(f)\leq \underline{\Jj}(u)+\varsigma$.
The map $\theta\mapsto\Gg_{1}(\theta f)$ is finite and convex on~$[0,\infty)$ (by H\"older's inequality), vanishes at $\theta=0$ and equals $u>0$ at $\theta=1$; it is therefore strictly increasing on~$[1,\infty)$, so that for every $u'>u$ sufficiently close to~$u$ there exists $\theta_{u'}>1$, with $\theta_{u'}\downarrow 1$ as $u'\downarrow u$, such that $\Gg_{1}(\theta_{u'}f)=u'$ and $\If(\theta_{u'}f)=\theta_{u'}^{2}\If(f)$.
Hence, 
$$
\limsup_{u'\downarrow u}\underline{\Jj}(u')\leq \underline{\Jj}(u)+\varsigma,
$$
and with $\varsigma>0$ being arbitrary, $u\mapsto\underline{\Jj}(u)$ is right-continuous on~$(0,\infty)$.
The open-set lower bound of Theorem~\ref{thm:LDP} then yields
$$
\liminf_{\eps\downarrow 0}\speed_{\eps}^{2}\log\PP\left(\frac{Z_{\eps}}{\eps}\geq K\right)
\geq -\inf_{z> K}\underline{\Jj}\left(\log\frac{z}{\sigma_{0}^{2}}\right)
= -\underline{\Jj}(u),
$$
and the same right-continuity handles the limit $\delta\downarrow0$ in the sandwich, proving the Call statement.
For Puts, let $0<K<\sigma_{0}^{2}$, so that $\log(\frac{K}{\sigma_{0}^{2}})<0$.
By Lemma~\ref{lem:prop_lambda}(ii), $\widetilde{\underline{\Jj}}$ is convex and finite on~$\RR$, hence continuous, and coincides with~$\underline{\Jj}$ on~$(-\infty,0]$; the same sandwich argument (with $\Put_{\eps}(K)\geq\delta\PP(Z_{\eps}/\eps\leq K-\delta)$), together with Corollary~\ref{cor:lefttail}, concludes.

\subsubsection{Proof of Proposition~\ref{prop:EL}}
Define
\begin{equation}
\Lambda[h]:=\half\int_{0}^{1}|\dot{h}(t)|^{2}\D t+\lambda\left(\int_{0}^{1}\sigma^{2}\left(y_{0}+(\Kk_{H}\dot{h})(t)\right)\D t-K\right),
\end{equation}
where $\lambda$ is the Lagrange multiplier.
We can compute
\begin{align}
&\Lambda[h+\delta\eta]-\Lambda[h]
\nonumber
\\
&=\half\int_{0}^{1}|\dot{h}(t)+\delta\eta'(t)|^{2}\D t+\lambda\int_{0}^{1}\sigma^{2}\left(y_{0}+\int_{0}^{t}K_{H}(s,t)(\dot{h}(s)+\delta\eta'(s))\D s\right)\D t
\nonumber
\\
&\qquad
-\half\int_{0}^{1}|\dot{h}(t)|^{2}\D t-\lambda\int_{0}^{1}\sigma^{2}\left(y_{0}+(\Kk_{H}\dot{h})(t)\right)\D t
\nonumber
\\
&=\delta\left(\int_{0}^{1}\dot{h}(t)\eta'(t)\D t+\lambda\int_{0}^{1}(\sigma^{2})'\left(y_{0}+(\Kk_{H}\dot{h})(t)\right)\int_{0}^{t}K_{H}(s,t)\eta'(s)\D s\D t\right) + \Oo\left(\delta^{2}\right),
\end{align}
as $\delta\downarrow 0$, which implies
that the optimal~$h$ satisfies
\begin{equation}
\int_{0}^{1}\dot{h}(t)\eta'(t)\D t+\lambda\int_{0}^{1}(\sigma^{2})'\left(y_{0}+(\Kk_{H}\dot{h})(t)\right)\int_{0}^{t}K_{H}(s,t)\eta'(s)\D s\D t=0,
\end{equation}
which holds for any $\eta$, which, by Fubini's theorem, is equivalent to 
\begin{equation}
\int_{0}^{1}\dot{h}(t)\eta'(t)\D t+\lambda\int_{0}^{1}\int_{s}^{1}(\sigma^{2})'\left(y_{0}+(\Kk_{H}\dot{h})(t)\right)K_{H}(s,t)\eta'(s)\D t \D s = 0,
\end{equation}
which can be rewritten as
\begin{equation}
\int_{0}^{1}\dot{h}(t)\eta'(t)\D t+\lambda\int_{0}^{1}\left[\int_{t}^{1}(\sigma^{2})'\left(y_{0}+(\Kk_{H}\dot{h})(s)\right)K_{H}(t,s)\D s\right]\eta'(t)\D t=0.
\end{equation}
Since it holds for any test function $\eta$, we conclude that
\begin{equation}
\dot{h}(t)+\lambda\int_{t}^{1}(\sigma^{2})'\left(y_{0}+(\Kk_{H}\dot{h})(s)\right)K_{H}(t,s)\D s=0,
\end{equation}
for any $0\leq t\leq 1$.
Hence, the optimal $h$ satisfies the above integral equation
with the constraint 
$\int_{0}^{1}\sigma^{2}\left(y_{0}+(\Kk_{H}\dot{h})(t)\right)\D t = K$, which completes the proof.

\subsubsection{Proof of Theorem~\ref{thm:ATM}}
When $K=\sigma^{2}(y_{0})$, by Taylor's theorem,
\begin{align}
\Call(K,T)&=\EE\left[\left(\frac{1}{T}\int_{0}^{T}\sigma^{2}(Y_{t})\D t-\sigma^{2}(y_{0})\right)_{+}\right]
\\
&=\EE\left[\left(\frac{1}{T}\int_{0}^{T}\left((\sigma^{2})'(y_{0})W_{t}^{H}+\half(\sigma^{2})''(R_{t})\left(W_{t}^{H}\right)^{2}\right)\D t\right)_{+}\right],
\end{align}
for some $R_{t}$, $0\leq t\leq T$, where $R_{t}$ lies between $y_{0}$ and $Y_{t}$
so that 
\begin{equation}\label{R:ineq}
|R_{t}|\leq\max(|y_{0}|,|Y_{t}|)\leq |y_{0}|+|Y_{t}|.
\end{equation}

Since $x\mapsto x_{+}$ is $1$-Lipschitz, by Assumption~\ref{assump:2}, 
\begin{align}
&\left|\Call(K,T)-\EE\left[\left(\frac{1}{T}\int_{0}^{T}(\sigma^{2})'(y_{0})W_{t}^{H}\D t\right)_{+}\right]\right|
\\
&\leq
\EE\left[\frac{1}{T}\int_{0}^{T}\half\left|(\sigma^{2})''(R_{t})\right|\left(W_{t}^{H}\right)^{2}\D t\right]
\nonumber
\\
&\leq
\EE\left[\frac{1}{T}\int_{0}^{T}\half\theta_{1}\E^{\theta_{2}|R_{t}|}\left(W_{t}^{H}\right)^{2}\D t\right]
\nonumber
\\
&\leq
\half\theta_{1}\E^{\theta_{2}|y_{0}|}
\EE\left[\frac{1}{T}\int_{0}^{T}\E^{\theta_{2}|Y_{t}|}\left(W_{t}^{H}\right)^{2}\D t\right]
\nonumber
\\
&\leq
\half\theta_{1}\E^{\theta_{2}|y_{0}|}
\frac{1}{T}\int_{0}^{T}\EE\left[\E^{2\theta_{2}|Y_{t}|}\right]^{\half}
\EE\left[\left(W_{t}^{H}\right)^{4}\right]^{\half}\D t,
\end{align}
using Cauchy-Schwarz inequality and~\eqref{R:ineq}.
Since $Y_{t}=y_{0}+W_{t}^{H}$
with $W_{t}^{H}\sim\Nn(0,t^{2H})$,
\begin{align}
\EE\left[\E^{2\theta_{2}|Y_{t}|}\right]
&\leq
\EE\left[\E^{2\theta_{2}Y_{t}}\right]
+\EE\left[\E^{-2\theta_{2}Y_{t}}\right]
\nonumber
\\
&=\left(\E^{2\theta_{2}y_{0}}+\E^{-2\theta_{2}y_{0}}\right)
\E^{2\theta_{2}^{2}t^{2H}}
\leq\left(\E^{2\theta_{2}y_{0}}+\E^{-2\theta_{2}y_{0}}\right)
\E^{2\theta_{2}^{2}},
\end{align}
for any $0\leq t\leq T\leq 1$ and
$\EE\left[\left(W_{t}^{H}\right)^{4}\right]
=3t^{4H}\leq 3T^{4H}$,
for any $0\leq t\leq T$ so that
\begin{equation}\label{error:ATM}
\left|\Call(K,T)-\EE\left[\left(\frac{1}{T}\int_{0}^{T}(\sigma^{2})'(y_{0})W_{t}^{H}\D t\right)_{+}\right]\right|
\leq
\frac{\theta_{1}}{2}\E^{\theta_{2}|y_{0}|}
\left(\E^{2\theta_{2}y_{0}}+\E^{-2\theta_{2}y_{0}}\right)^{\half}
\E^{\theta_{2}^{2}}\sqrt{3}T^{2H},
\end{equation}
for any $T\leq 1$.
Finally, by symmetry of the centred Gaussian law of $\int_{0}^{T}W_{t}^{H}\D t$, 
\begin{equation}\label{final:follow:1}
\EE\left[\left(\frac{1}{T}\int_{0}^{T}(\sigma^{2})'(y_{0})W_{t}^{H}\D t\right)_{+}\right]
=\frac{1}{T}\left|(\sigma^{2})'(y_{0})\right|\EE\left[\left(\int_{0}^{T}W_{t}^{H}\D t\right)_{+}\right],
\end{equation}
and $\int_{0}^{T}W_{t}^{H}\D t$
is centred Gaussian with variance
$\frac{T^{2H+2}}{2(H+1)}$.
Therefore, 
\begin{equation}\label{final:follow:2}
\EE\left[\left(\int_{0}^{T}W_{t}^{H}\D t\right)_{+}\right]
=\frac{T^{H+1}}{\sqrt{2(H+1)}}\EE\left[\Nn(0,1)_{+}\right]
=\frac{T^{H+1}}{2\sqrt{\pi(H+1)}}.
\end{equation}
It then follows from~\eqref{final:follow:1} and~\eqref{final:follow:2} that
\begin{equation}\label{final:follow:3}
\EE\left[\left(\frac{1}{T}\int_{0}^{T}(\sigma^{2})'(y_{0})W_{t}^{H}\D t\right)_{+}\right]
 = \frac{\left|(\sigma^{2})'(y_{0})\right|}{2\sqrt{\pi(H+1)}}T^{H}.
\end{equation}
Hence, we conclude from~\eqref{error:ATM} and~\eqref{final:follow:3} that
\begin{equation}
\lim_{T\downarrow 0}\frac{\Call(K,T)}{T^{H}}=\frac{\left|(\sigma^{2})'(y_{0})\right|}{2\sqrt{\pi(H+1)}},
\end{equation}
and similarly for ATM Put options, which completes the proof.


\subsubsection{Proof of Proposition~\ref{prop:Jexpansion}}
\label{sec:proof_prop:Jexpansion}
The solution of the variational problem for the rate function satisfies an Euler-Lagrange equation as seen in Proposition~\ref{prop:EL},
which depends only on $\dot{h}(t)$. 
Writing the integrations explicitly, 
\begin{equation}\label{EL}
\dot{h}(t) = \lambda \int_t^1 (\sigma^2)' \left( \int_0^s K_{H}(v, s) f(v) \D v\right) 
K_{H}(t,s)\D s,
\end{equation}
with $K_{H}(s,t) = \sqrt{2H}(t-s)^{H-\half}$
and $\lambda$ is a Lagrange multiplier.
In addition, $\dot{h}$ satisfies
\begin{equation}\label{constraint}
\int_{0}^{1} \sigma^2\left( \sqrt{2H} \int_0^t (t-s)^{H-\half} f(s)\D s \right) \D t = K = K_{\mathrm{ATM}} \E^y = v_0 \E^y,
\end{equation}
where $y=\log(\frac{K}{K_{\mathrm{ATM}}})$ is log-moneyness.
At the money $K=K_{\mathrm{ATM}}=v_0$ both the Euler-Lagrange equation and the constraint are satisfied with 
$\dot{h}(t) = 0$. 
We expand the optimiser function in log-moneyness $y=\log(\frac{K}{K_{\mathrm{ATM}}})$:
$$
\dot{h}(t) := y \dot{h}_1(t) + y^2 \dot{h}_2(t) + \Oo(y^3).
$$
We also expand the Lagrange multiplier as
$\lambda = \lambda_0 y + \lambda_1 y^2 + \Oo(y^3)$.
Substituting into the Euler-Lagrange equation and expanding in $y$, we obtain to first two orders 
\begin{align}\label{EL1}
\Oo(y): \dot{h}_1(t) &= \lambda_0 v_1 \int_t^1 K_{H}(t,s)\D s \\
\label{EL2}
\Oo\left(y^2\right): \dot{h}_2(t) &= 2v_2 \lambda_0 \int_t^1 \left( \int_0^s K_{H}(v,s) \dot{h}_1(v) \D v\right)
K_{H}(t,s)\D s \\
&\qquad\qquad+ \lambda_1 v_1 \int_t^1 K_{H}(t,s)\D s. \nonumber
\end{align}

Expanding also the 
constraint~\eqref{constraint} gives additional equations 
\begin{align}\label{C1}
\Oo(y): & \  v_1\int_{0}^{1} \left(\int_0^t K_{H}(s,t)  \dot{h}_1(s)\D s\right) \D t  = v_0,\\
\label{C2}
\Oo(y^2): &  \ v_2 \int_{0}^{1} \left|\int_0^t K_{H}(s,t)  \dot{h}_1(s)\D s\right|^2 \D t 
+
v_1 \int_{0}^{1} \int_0^t K_{H}(s,t)  \dot{h}_2(s)\D s \D t =\frac{v_0}{2}.
\end{align}

These equations are easily solved, order by order in $y$, and the rate function then reads
\begin{align}\label{rate:function:expansion}
J_H(K) &= \left(\half \int_{0}^{1} (\dot{h}_1(t))^2 \D t\right) y^2 + \left( \int_{0}^{1} \dot{h}_1(t) \dot{h}_2(t) \D t\right) y^3 \\
&\qquad\qquad + \left(\int_{0}^{1} \left( \half (\dot{h}_2(t))^2 + \dot{h}_1(t) \dot{h}_3(t) \right) \D t \right)  y^4 + \Oo\left(y^5\right).\nonumber
\end{align}

\textbf{Solution for $\dot{h}_1(t)$.}
From the $\Oo(y)$ term in the expansion
of~\eqref{EL1} we obtain
\begin{equation}\label{f1pre}
\dot{h}_1(t) = \lambda_0  v_1 \sqrt{2H} \int_t^1 (s-t)^{H-\half } \D s = 
\lambda_0 v_1 \frac{\sqrt{2H}}{H+\half } (1-t)^{H+\half }.
\end{equation}
The Lagrange multiplier $\lambda_0$ is determined by substituting this 
into~\eqref{C1}:
$\lambda_0 v_1 =  \frac{(H+\half )^2 (H+1)}{H} \frac{v_0}{v_1}$.
Substituting this into~\eqref{f1pre} gives
\begin{equation}\label{f1sol}
\dot{h}_1(t) = c_1 (1-t)^{H+\half },\qquad 
\text{with }c_1 := \frac{2(H+\half )(H+1)}{\sqrt{2H}} \frac{v_0}{v_1},
\end{equation}
and the $\Oo\left(y^2\right)$ term in the rate function expansion~\eqref{rate:function:expansion} reads
\begin{align}
j_2 &= \half \int_{0}^{1} (\dot{h}_1(t))^2 \D t = \frac{(H+\half )^2(H+1)}{2H} \frac{v_0^2}{v_1^2}.
\end{align}

\textbf{Solution for $\dot{h}_2(t)$.} Substituting the solution for $\dot{h}_1(t)$ in~\eqref{f1sol} into~\eqref{EL2} gives
\begin{align}\label{f2sol}
\dot{h}_2(t) 
&=  \frac{4 v_2 \lambda_0 H c_1}{H+\half } \int_t^1 (s-t)^{H-\half } s^{H+\half } 
{}_2 F_1\left(1,-H-\half,H+\frac{3}{2},s\right) \D s 
\\
&\qquad\qquad\qquad\qquad\qquad\qquad
+  \frac{\lambda_1 v_1\sqrt{2H}}{H+\half } (1-t)^{H+\half }.
\nonumber
\end{align}
The Lagrange coefficient $\lambda_1$ is found by substituting this into the $\Oo\left(y^2\right)$ constraint~\eqref{C2}:
$$
v_1 \lambda_1 = \frac{(H+\half )^2 (H+1)}{2H} \frac{v_0}{v_1} 
\left( 1 - \frac{24 v_0 v_2}{v_1^2} (H+1)^2 R(H) \right).
$$
The coefficient of the $\Oo\left(y^3\right)$ term in the rate function expansion~\eqref{rate:function:expansion} 
is $j_3 =\int_{0}^{1} \dot{h}_2(t) \dot{h}_1(t) \D t$.
Substituting here $\dot{h}_1(t)$ from~\eqref{f1sol} and $\dot{h}_2(t)$ from~\eqref{f2sol} gives the stated result.

\bibliographystyle{siam}
\bibliography{references}


\end{document}